\documentclass[12pt,letterpaper,oneside]{book}  
\usepackage[utf8]{inputenc}
\usepackage{cmap}
\usepackage[T1]{fontenc}
\usepackage[brazil]{babel}
\usepackage{indentfirst}
\usepackage[top=3cm,bottom=3cm,right=2cm,left=2cm]{geometry}
\usepackage{lmodern}  
\usepackage{slantsc}  
\usepackage{chngcntr}
\usepackage{mathrsfs}
\counterwithin{figure}{chapter}
\usepackage{enumitem}
\usepackage{capt-of}
\usepackage{commath}
\usepackage{multicol}

\usepackage{etoolbox}  

\usepackage{url}
\usepackage{breakurl}
\usepackage{hyperref}

\usepackage{amsmath}
\usepackage{amsfonts}
\usepackage{amssymb}
\usepackage{amsthm}
\usepackage{breqn}

\usepackage{multicol}
\usepackage{multirow}
\usepackage{array}
\usepackage{booktabs}

\usepackage[usenames,dvipsnames,svgnames,table]{xcolor}

\usepackage{pdfpages}
\usepackage{graphicx}
\usepackage{epstopdf}
\usepackage{tikz}
\usetikzlibrary{fit}
\usepackage{wrapfig}

\usepackage{algorithmic}
\usepackage[chapter]{algorithm}
\floatname{algorithm}{Algoritmo}

\usepackage{textcomp}
\usepackage{listings}

\usepackage{makeidx}
\makeindex

\usepackage{lmodern} \bfseries
\usepackage{mathrsfs}

\newcommand{\autor}{Ciro Javier Díaz Penedo}
\newcommand{\Trabalho}{Exame de qualificação}
\newcommand\numberthis{\addtocounter{equation}{1}\tag{\theequation}}
\newcommand{\titulopt}{Uso combinado do método de elementos finitos mistos híbridos com decomposição de domínio e de métodos espectrais para um estudo de renormalização do modelo KPZ}
\newcommand{\titulo}{Combined use of mixed and hybrid finite elements method with domain decomposition and spectral methods for a study of renormalization for the KPZ model}
\def\aplicada{}
\newcommand{\orientador}{Eduardo Cardoso de Abreu}
\newcommand{\ano}{2015}
 
\appto\frontmatter{\pagestyle{plain}}  

\numberwithin{equation}{section}
\numberwithin{section}{chapter}

\theoremstyle{definition}
\newtheorem{thm}{Teorema}[section]

\newtheorem{dfn}[thm]{Definição}
\newtheorem{spn}[thm]{Suposição}

\newtheorem{pps}[thm]{Proposição}

\newcommand{\annexname}{Anexo}
\makeatletter
\newcommand\annex{\par
\setcounter{chapter}{0}%
\setcounter{section}{0}%
\gdef\@chapapp{\annexname}%
\gdef\thechapter{\@Roman\c@chapter}}
\makeatother

\newcommand{\subf}[2]{%
  {\small\begin{tabular}[t]{@{}c@{}}
  #1\\#2
  \end{tabular}}%
}
\renewcommand{\arraystretch}{2.0}

\begin{document}
 
\frontmatter

\thispagestyle{plain}
\includegraphics[width=.94in, height=1in,
keepaspectratio=true]{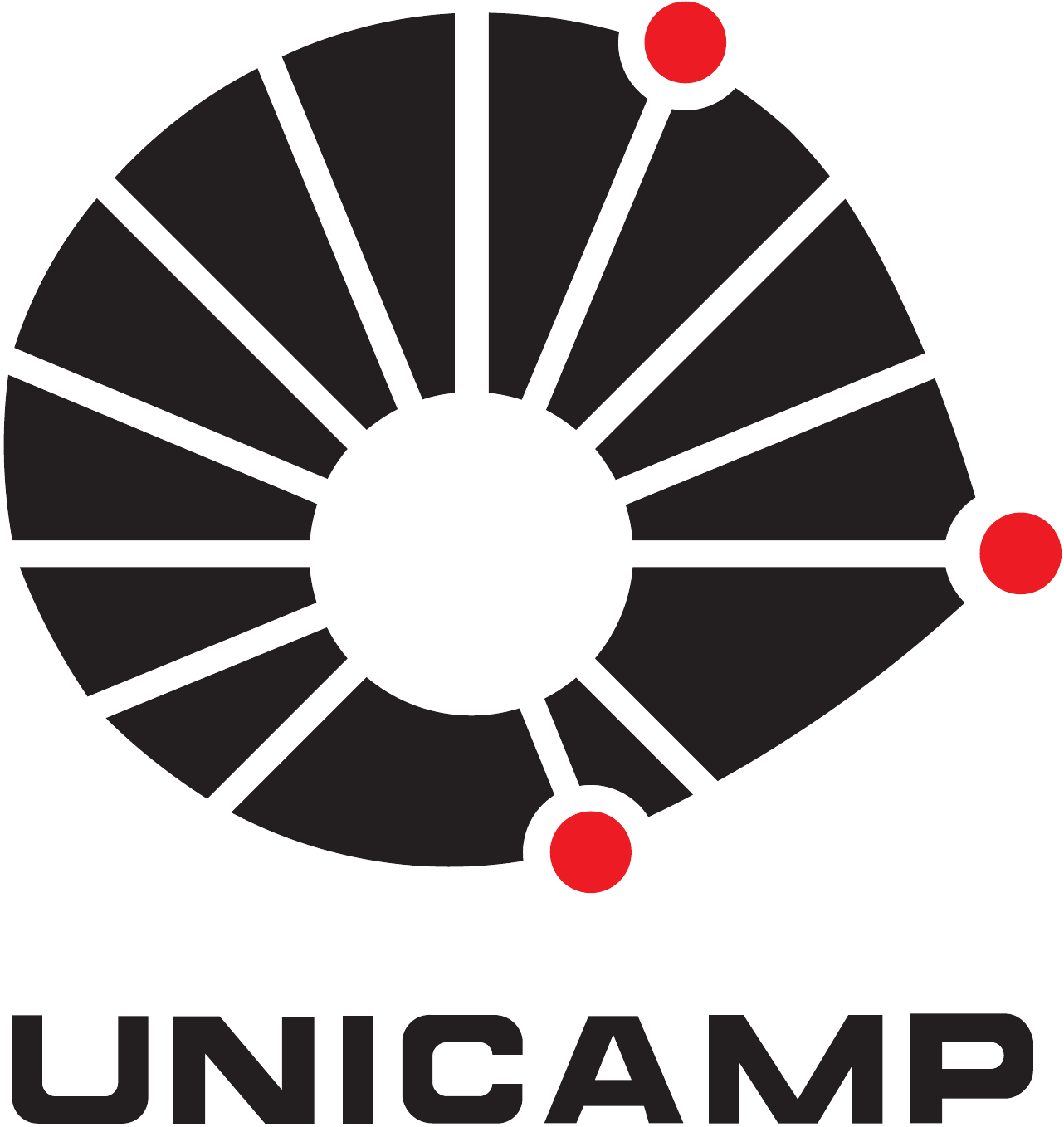}
\vspace*{.8cm}
\begin{center}
  {\Large\textsc{\autor}}
\end{center}
\begin{center}
  {\small\textsc{\Trabalho}}
\end{center}
\vspace{1.9cm}
\begin{center}
  {\Large\textbf{\textsc{\titulopt}}}
\end{center}
\vspace{2.2cm}
\begin{center}
  {\Large\textsl{\textsc{\titulo}}}
\end{center}
\vspace{3cm}
\begin{center}
  {\Large\textsc{orientador : \orientador}}
\end{center}
\vfill
\begin{center}
  \textbf{CAMPINAS \\ \ano}
\end{center}
  
\newpage\mbox{}\thispagestyle{plain}\newpage  

\thispagestyle{plain}
\includegraphics[width=.94in, height=1in,
keepaspectratio=true]{unicamp-logo}
\begin{center}
  {\large\textbf{\textsc{Universidade Estadual de Campinas}}
  \vspace{.4cm}

  Instituto de Matemática, Estatística \\
  e Computação Científica}
\end{center}
\vspace{.01cm}
\begin{center}
  {\large\textbf{\textsc{\autor}}}
\end{center}
\vspace{.01cm}
\begin{center}
  {\Large\textbf{\textsc{\titulo}}}
\end{center}
\vspace{.01cm}
\begin{center}
  {\Large\textsl{\textsc{\titulopt}}}
\end{center}
\vspace{.01cm}

\begin{flushright}
  \begin{minipage}[c]{.5\textwidth}
    \ifx\mestrado\undefined
    Thesis
    \else
    Dissertation
    \fi
    presented to the Institute of Mathematics, Statistics and Scientific Computing
    of the University of Campinas in partial
    fulfillment of the requirements for the degree of
    \ifx\mestrado\undefined
    \ifx\femaleAuthor\undefined
    Mestre
    \else
    Mestre
    \fi
    \else
    \ifx\femaleAuthor\undefined
    Master
    \else
    Master
    \fi
    \fi
    in
    \ifx\matematica\undefined
    \else
    mathematics.
    \fi
    \ifx\aplicada\undefined
    \else
    applied mathematics.
    \fi
    \ifx\estatistica\undefined
    \else
    statistics.
    \fi
  \end{minipage}
\end{flushright}

\vfill
\begin{flushright}
  \begin{minipage}[c]{.5\textwidth}
    \textit{\ifx\mestrado\undefined
    Tese
    \else
    Dissertação
    \fi
    apresentada ao Instituto de Matemática,
    Estatística e Computação Científica da Universidade
    Estadual de Campinas como parte dos requisitos exigidos
    para a obtenção do título de
    \ifx\mestrado\undefined
    \ifx\femaleAuthor\undefined
    Mestre
    \else
    Mestre
    \fi
    \else
    \ifx\femaleAuthor\undefined
    Mestre
    \else
    Mestre
    \fi
    \fi
    em
    \ifx\matematica\undefined
    \else
    matemática.
    \fi
    \ifx\aplicada\undefined
    \else
    matemática aplicada.
    \fi
    \ifx\estatistica\undefined
    \else
    estatística.
    \fi
    }
  \end{minipage}
\end{flushright}
\vspace{.5cm}

\noindent
\textbf{Orientador\ifx\femaleOrientador\undefined
\else
a\fi: \orientador
}
\vspace{.25cm}

\ifx\coorientador\undefined
\else
\noindent
\textbf{Coorientador\ifx\femaleCoorientador\undefined
\else
a\fi: \coorientador
}
\vspace{.5cm}
\fi

\noindent
\begin{minipage}[c]{.5\textwidth}
  {\footnotesize\textsc{Este exemplar n\~ao corresponde à versão final da
  \ifx\mestrado\undefined
  tese
  \else
  dissertação
  \fi
  defendida
  \ifx\femaleAuthor\undefined
  pelo aluno
  \else
  pela aluna
  \fi
  \autor,
  e orientada pel\ifx\femaleOrientador\undefined
  o\else
  a\fi{} Prof\ifx\femaleOrientador\undefined
  \else
  a\fi. Dr\ifx\femaleOrientador\undefined
  \else
  a\fi. \orientador.
  }}
\end{minipage}
\vspace{.5cm}

\noindent
{\small\textbf{Assinatura
\ifx\femaleOrientador\undefined
do Orientador
\else
da Orientadora
\fi
}

\vspace{.3cm}
\noindent
\rule[1pt]{7cm}{.5pt}  
}
\vspace{.4cm}

\ifx\coorientador\undefined
\else
\noindent
{\small\textbf{Assinatura
\ifx\femaleCoorientador\undefined
do Coorientador
\else
da Coorientadora
\fi
}

\vspace{.3cm}
\noindent
\rule[1pt]{7cm}{.5pt}  
}
\fi
\vfill
\begin{center}
  {\small\textbf{\textsc{ Campinas \\ \ano}}}
\end{center}

\includepdf{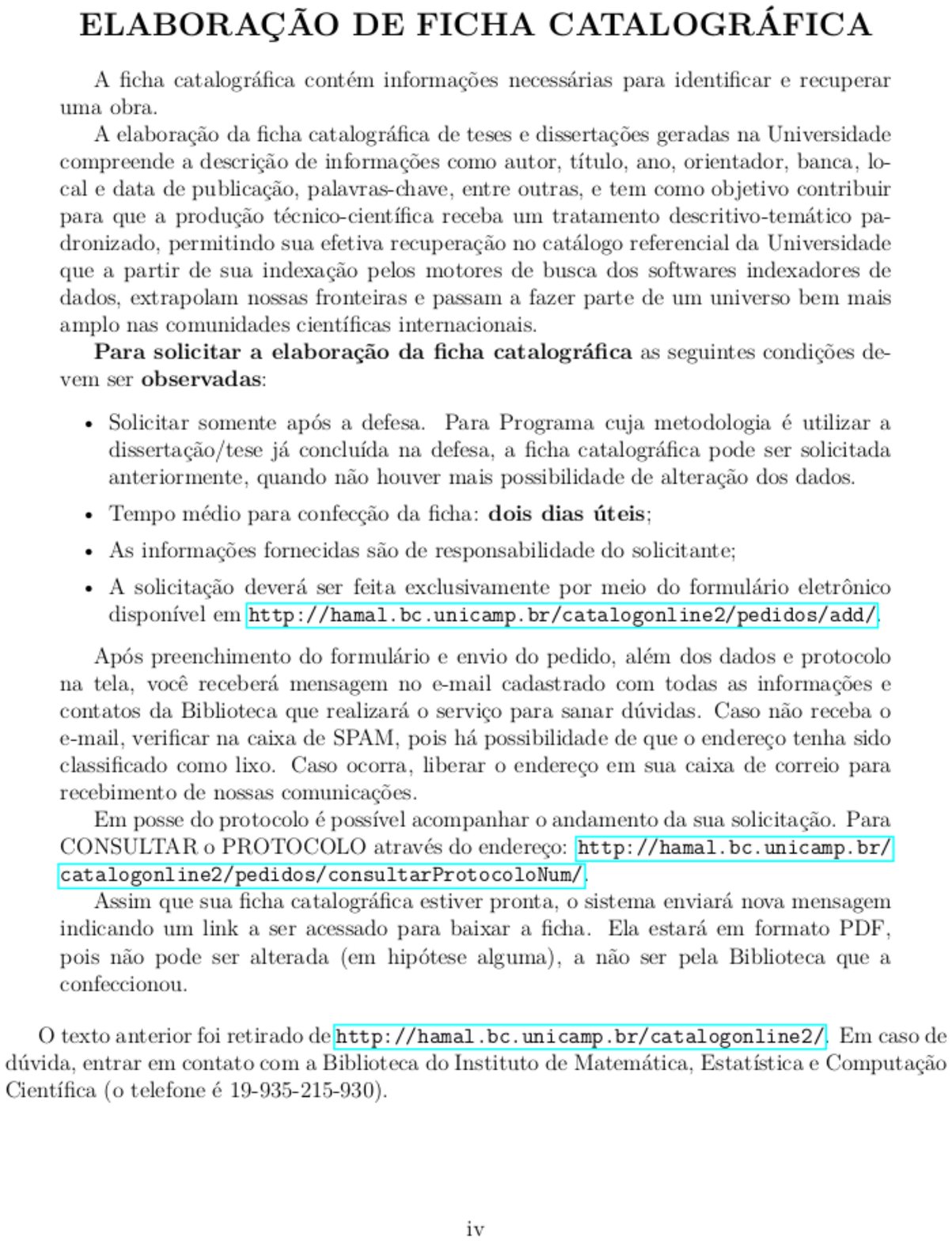}
\includepdf{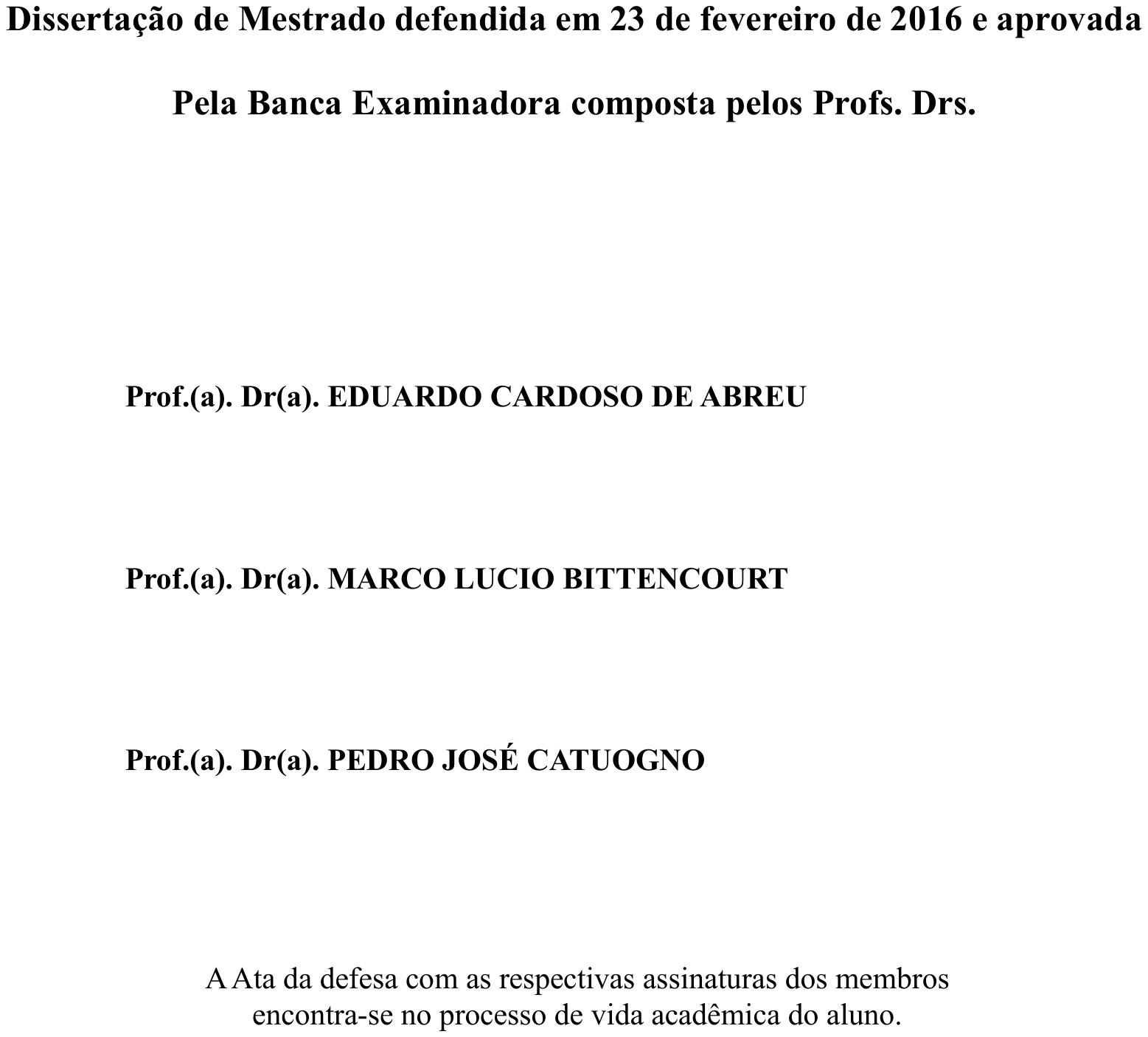}
\chapter*{}  
\begin{center}
  \large{\textbf{Abstract}}
\end{center}

The focus of this work is the numerical approximation of time-dependent partial differential equations associated to initial-boundary value problems. This master dissertation is mostly concerned with the actual computation of the solution to nonlinear stochastic evolution problems governed by Kardar-Parisi-Zhang (KPZ) models. In addition, the dissertation aims to contribute to corroborate, by means of a large set of numerical experiments, that the initial-boundary value problem with periodic boundary conditions for the equation KPZ is ill-posed and that such equation needs to be renormalized. The approach to discretization of KPZ equation perfomed by means of the use of hybrid and mixed finite elements with a domain decomposition procedure along with a pertinent mollification of the noise. The obtained solution is compared with the well known solution given by the Cole-Hopf transformation of the stochastic heat equation with multiplicative noise. We were able to verify that both solutions exhibit a good agreement, but there is a shift that grows as the support of the mollifier decreases. For the numerical aproximation of the stochastic heat equation we use a state-of-the-art numerical method for evaluating semilinear stochastic PDE , which in turn combine spectral techniques, Taylor’s expantions and particular numerical treatment to the underlying noise. Furthermore, a state-of-the-art renormalization procedure introduced by Martin Hairer is used to renormalize KPZ equation that is validated with nontrivial numerical experiments.

\vspace{.2cm}
\textbf{Keywords}:

Stochastic partial diferential equations, Mixed and hybrid finite elements, Domain decomposition, Expectral methods, Discretization of the KPZ model, Renormalization
\vspace{.1\textheight}
\begin{center}
  \large{\textbf{Resumo}}
\end{center}

O foco deste trabalho é realizar aproximação numérica de equações diferenciais parciais dependentes do tempo associadas a problemas de valor inicial e de contorno, em particular,  problemas de evolução estocásticos não lineares regidos por modelos da classe Kardar-Parisi-Zhang (KPZ). Além disso, o trabalho visa corroborar, por meio de um grande conjunto de experimentos numéricos, que o problema de valor inicial com condições de contorno periódicas para a equação KPZ, é mal posto e que  a equação precisa ser renormalizada. A discretização da equação KPZ é realizada por meio do uso de elementos finitos mistos e híbridos, juntamente com um procedimento de decomposição do domínio e um pertinente amolecimento do ruído. Por sua vez, a solução obtida é comparada com a bem conhecida transformação de Cole-Hopf da solução da equação estocástica do calor com ruído multiplicativo. Ao longo do desenvolvimento deste trabalho, foi verificado que os perfis de ambas soluções exibem uma boa concordância, porém há um crescente deslocamento à medida que o suporte do \textit{mollifier} diminui. Para a aproximação numérica da equação estocástica do calor utilizamos métodos numéricos recentemente desenvolvidos para equações estocásticas semilineares, que por sua vez, combinam técnicas espectrais, expansão de Taylor e uma abordagem particular do ruído. Além disso, um procedimento de renormalização  introduzido por Martin Hairer é usado para renormalizar a equação KPZ e a sua eficácia é validada com experimentos numéricos não triviais.

\vspace{.2cm}
\textbf{Palavras-chave}:

Equações diferenciais estocásticas, Métodos computacionais para o modelo KPZ, Métodos espectrais, Elementos finitos mistos híbridos, Equações estocásticas semilineares, Decomposição de domínio, Renormalização.
\tableofcontents

\mainmatter
\chapter{Introdução}\label{cap1}

  A equação KPZ foi introduzida em (\cite{KPZ86}) 
  
  \begin{equation*}
  \partial_t h = \nu\partial^2_x h + \lambda (\partial_x h)^2 + \xi,
  \end{equation*}
\\
  com a intenção de modelar o processo de deposição balística, neste caso $\xi$ representa um ruído branco 
  no espaço e no tempo. Em pouco tempo muitos
  cientistas voltaram sua atenção a esta equação que tentava capturar o crescimento lateral de uma 
  superfície em evolução, propriedade que escapava aos modelos anteriores (Deposição aleatória e 
  Edward-Wilkinson). As primeiras tentativas de 
  aproximar possíveis soluções de 
  problemas iniciais e de contorno envolvendo a equação KPZ fracassaram. Paralelamente, foi ficando
  claro que tais problemas eram mal postos e que algum processo de
   renormalização poderia dar sentido a esta equação.\\
   
   Do ponto de vista de métodos de aproximação no contexto da análise numérica, desde 
  a introdução do modelo KPZ até hoje, foram discutidos diversos esquemas numéricos 
   com a ``ilusão'' de obter aproximações de uma solução que descrevesse o 
  processo de deposição balística. Por exemplo, esquemas em diferenças 
  finitas \cite{Tiago}, \cite{NuS}, métodos pseudo-espectrais \cite{GGR08} além do 
  princípio de mínima ação \cite{MAM09}. Muitos destes métodos devolviam resultados que 
  recuperavam algumas das propriedades que exibem as superfícies obtidas por processos 
  de crescimento porém, estes resultados ficavam longe de serem interpretados como soluções. Assim, 
  a interpretação das soluções da equação KPZ ficou como um problema aberto nas últimas duas 
  décadas.

\section{Motivação da pesquisa}\label{sec1}
  
  Nos últimos anos muitos progressos foram obtidos para uma melhor compreensão da equação KPZ. 
  Um exemplo notável foi o trabalho seminal de Martin Hairer (2013), 
  no mesmo, o autor introduziu novos conceitos de aproximação de 
  quantidades estocásticas que culminou em uma nova noção de solução
  da equação KPZ. Desta maneira é de se esperar que 
  novos procedimentos numéricos devam ser igualmente investigados. Resultados neste sentido 
  podem ser encontrados nos trabalhos de Arnulf Jentzen \cite{JER12,JEK11,JEN11,JEK08} 
  onde uma inovadora estratégia para conectar métodos espectrais e expansões 
  em serie de Taylor é acompanhada de uma rigorosa teoria de aproximação. 
  
  Existem poucos trabalhos apresentando experimentos numéricos  que mostrem a eficácia 
  dos processos de renormalização introduzidos em (\cite{Hai11}) ou 
  em (\cite{FrH14}). É por isso que decidimos enfatizar neste trabalho a apresentação de um conjunto 
  representativo de experimentos numéricos em concordância com os resultados teóricos reportados na
  literatura. Em particular os métodos numéricos empregados serão também descritos.\\

  Com base na revisão bibliográfica realizada até o presente momento, métodos de 
  elementos finitos clássicos ou  métodos de elementos finitos 
  mistos e híbridos, não aparecem relacionados à resolução numérica da equação KPZ.
  Por isso, consideramos pertinente a proposta de explorar 
  o desempenho desse tipo de metodologia na aproximação da solução da equação KPZ, acompanhada 
  de um processo de renormalização que permita uma reinterpretação adequada das soluções.\\

\section{Objetivos específicos e proposta da dissertação}\label{sec2}

    Este trabalho visa mostrar experimentalmente que a equação 
    KPZ (clássica) deve ser renormalizada para fazer sentido, e que as soluções desta equação 
    renormalizada se aproximam à transformada de Hopf-Cole da solução da equação estocástica do calor
    com ruído branco multiplicativo. Também experimentalmente, mostraremos que a escolha 
    do \textit{mollifier} não altera o resultado anterior. Isso permite concluir que o uso combinado 
    de teoria de aproximação com métodos numéricos pertinentes via analise numérica é 
    uma ferramenta matemática relevante para a compreensão da equação KPZ e de outras equações
     diferenciais estocásticas relacionadas. Para tal fim procedemos na seguinte ordem.\\
    
\begin{enumerate}[label=(\alph*)]

\item Estudo e implementação de algoritmos de aproximação para equações estocásticas 
      semilineares baseados 
      em métodos espectrais combinados com expansões de Taylor. Estes métodos serão utilizados na 
      aproximação da solução da equação estocástica do calor com ruído multiplicativo.

\item Adaptar e implementar o método de elementos finitos mistos e híbridos com decomposição de domínio      
      (EFMH-DD) para ser aplicado no problema de valor inicial e de condições de contorno para a 
      equação KPZ com ruído branco amolecido. 
      
\item Medir experimentalmente o desempenho dos métodos discutidos para obter indícios de que 
      os algoritmos estejam aproximando corretamente as soluções desejadas, fazendo uso dos itens 
      \textit{(a)} e \textit{(b)}.

\item Validar experimentalmente os resultados teóricos referentes à conexão entre a 
      transformada de Hopf-Cole da solução da equação estocástica do calor com ruído branco 
      multiplicativo e o limite de um processo de renormalização aplicado à equação KPZ, com base 
      nos itens \textit{(a)}, \textit{(b)} e \textit{(c)}. 

\end{enumerate}

\section{Resultados}\label{sec3}

Lista-se o conjunto de resultados atingidos, face aos objetivos desta dissertação de mestrado:

\begin{enumerate}[label=(\alph*)]      

\item    Foram implementados os códigos\textit{ Milstein.m}, \textit{LordRougemont.m} e
         \textit{EulerGalerkinSemimplícito.m}  para modelar computacionalmente os métodos numéricos 
         de Milstein, Lord-Rougemont e Euler-Galerkin-semi-implícito, respetivamente.  Estes 
         códigos foram utilizados na aproximação das soluções de equações estocásticas semilineares. 
          
\item   Foi feito um estudo numérico para validar os códigos implementados em problemas onde
        tem-se a solução conhecida (ver \cite{BBHW00}). 
        Um estudo de erro revelou que  este decresce à medida que aumentamos o número de 
        funções base na expansão da solução. Além disso, foi recuperado via transformada de Cole-Hopf,
        o comportamento previsto da rugosidade para o processo de deposição balística.

\item    Foi feita uma construção formal de uma nova formulação numérica via método de elementos finitos
         mistos e  híbridos com decomposição de domínio (EFMH-DD) para um modelo KPZ com ruído
         branco amolecido.   
             
\item    Foi implementado o código \textit{EFMH_KPZ.m} para modelar computacionalmente o método EFMH-DD
         e foi usado em um estudo numérico  para reprodução de resultados apresentados  em \cite{BBHW00},
         onde foi considerado um modelo KPZ determinístico.  
         Nesse contexto, também foram realizados alguns experimentos numéricos com o objetivo de 
         medir o desempenho do algoritmo.
          
\item   Foi implementado um código que compara a solução aproximada da equação KPZ 
        (amolecendo o ruído) do item \textit{(c)} e \textit{(d)} com a transformada de Hoph-Cole 
        da solução aproximada da equação estocástica do calor do item \textit{(a)}. 
        Os experimentos mostraram a  similaridade entre os dois perfis quando consideramos a 
        mesma realização do ruído branco.

\item   Foi feito um estudo numérico representativo para verificar que a conexão entre as duas soluções 
        do item \textit{(c)} não depende do \textit{mollifier} utilizado.
        
\item   Foi aplicado o processo de renormalização proposto em \cite{Hai11} para o modelo KPZ 
        e foram realizados experimentos numéricos nos quais aproximamos a solução como em \textit{(d)}. 
        Estes experimentos forneceram evidências de que tal processo  atenua a divergência provocada 
        pelo termo $(\partial_x h)^2$ e que o processo limite independe da escolha do \textit{mollifier}.   
            
\end{enumerate}

\section{Organização do trabalho}\label{sec4}

  No Capítulo \ref{cap11}, apresentamos algumas definições e resultados a serem 
  utilizados ao longo deste trabalho. No Capítulo \ref{cap2}, fazemos um breve estudo da equação KPZ,
  detalhando aspectos teóricos necessários para a construção dos métodos propostos. 
  No Capítulo \ref{cap33}, fazemos um estudo 
  de alguns métodos numéricos para equações diferenciais estocásticas semilineares onde o caso
  de maior interesse é a equação estocástica do calor com ruído branco multiplicativo.
  No Capítulo \ref{cap3}, introduzimos o método de elementos finitos mistos e híbridos com 
  decomposição de domínio adaptado para ser aplicado na
  equação KPZ determinística. Também apresentamos os resultados obtidos nas simulações onde utilizamos 
  o código \textit{EFMH_KPZ.m} e comparamos estes com outros 
  reportados na literatura considerando uma equação KPZ determinística. 
  No Capítulo \ref{cap4},  mostramos experimentalmente  que a transformada de Hopf-Cole 
  das soluções numéricas da equação estocástica do 
  calor obtidas usando os métodos introduzidos no Capítulo 
  \ref{cap33}, se aproximam à solução da equação KPZ renormalizada com ruído branco 
  amolecido. Para amolecer o ruído, utilizamos um \textit{mollifier}, que em nosso caso será
   uma função real, de classe $\mathcal{C}^{\infty}$ e com suporte compacto.
   Além disso, verificamos experimentalmente que o limite do processo de renormalização existe quando 
  o diâmetro do suporte vá para zero e que independe da escolha do \textit{mollifier}. 
  No final deste capitulo implementamos o processo de renormalização proposto em \cite{Hai11} 
  e apresentamos alguns resultados numéricos obtidos aplicando o  método de elementos finitos 
  mistos e híbridos adaptado para a equação renormalizada. Por fim, no Capítulo \ref{cap5}, 
  são reportadas as conclusões e perspectivas desta dissertação de mestrado.
\chapter{Preliminares}\label{cap11}

A continuação mostram-se algumas notações, definições e teoremas a serem utilizados nesta dissertação de mestrado. Estamos seguindo, fundamentalmente, a notação e o estilo do texto (\cite{LPS14}).

\section{Notações, definições e teoremas}

\begin{dfn}\label{HSc} \textbf{(Norma Hilbert-Schmidt)}  Sejam $H$ e $U$ espaços de Hilbert separáveis com normas $\parallel \cdot \parallel$ e $\parallel \cdot \parallel_U$. Para uma base ortonormal $\lbrace\varphi_j : j \in \mathbb{N} \rbrace \subset U$ definimos a norma de Hilbert-Schmidt como,

\begin{equation}
  \parallel L \parallel_{HS(U,H)} := \left( \sum_{j=1}^{\infty} \parallel L\varphi_j \parallel_U^2 \right)^{1/2}\!\!\!\!\!\!.
\end{equation}
O conjunto $HS(U,H):=\lbrace L \in \mathcal{L}(U,H): \parallel L \parallel_{HS(U,H)} < +\infty \rbrace$ é um espaço de Banach munido da norma de Hilbert-Schmidt. Um operador $L \in HS(U,H)$ é conhecido como um operador de Hilbert-Schmidt.
\end{dfn} 

  Ao longo deste trabalho vamos considerar somente espaços de probabilidade filtrados $(\Omega, \mathcal{F}, \mathcal{F}_t, \mathbb{P})$ que vamos abreviar escrevendo $(\Omega, \mathcal{F}, \mathbb{P})$.

%
%

\begin{dfn}\textbf{(variáveis aleatórias)} Seja $\left(\Omega,\mathcal{F},\mathbb{P} \right)$  um espaço de probabilidade e $(\Psi, \mathcal{G})$ um espaço de medida. Então, $X$ é uma variável aleatória que toma valores em $\Psi$ se $X$ é uma função medível  de $\left(\Omega,\mathcal{F}\right)$ em $(\Psi, \mathcal{G})$. Para enfatizar a $\sigma$-álgebra sobre $\Omega$, podemos escrever que $X$ é uma variável aleatória $\mathcal{F}$-mensurável. O valor observado $X(\omega)$  para um dado $\omega \in \Omega$ é chamado de uma realização de $X$.
\end{dfn}

Neste trabalho vamos trabalhar com variáveis aleatórias reais, ou seja $(\Psi, \mathcal{G}) = (\mathbb{R},\mathcal{B}(\mathbb{R}))$ onde $\mathcal{B}(\mathbb{R})$ é a $\sigma$-algebra de Borel.


\begin{dfn}(valor esperado)\label{ValEsp}
 Seja $X$ uma variável aleatória que toma valores num espaço  de Banach sobre o espaço de probabilidade $\left(\Omega,\mathcal{F},\mathbb{P} \right)$. Se $X$ é integrável, a expectativa de $X$ é
\begin{equation}
\mathbb{E}(X) = \int_{\Omega}X(\omega)d\mathbb{P}(\omega) \label{Int},
\end{equation}
a integral de $X$ respeito da medida de probabilidade $\mathbb{P}$.
\end{dfn}

\begin{dfn}(covariância) A covariância entre duas variáveis aleatórias reais $X$ e $Y$ define-se como
\begin{equation}
Cov(X,Y)=\mathbb{E}\left[ (X-\mu_X)(Y - \mu_Y)  \right] = \mathbb{E}\left[XY\right] - \mu_X\mu_Y.
\end{equation}
\end{dfn}

\begin{dfn}
(\textbf{processo estocástico}). Seja um conjunto $T \subset \mathbb{R}$, um espaço mensurável $(H,\mathcal{H})$, e um espaço de probabilidade $(\Omega, \mathcal{F}, \mathbb{P})$. Chamaremos de processo estocástico que toma valores em $H$, um conjunto de variáveis aleatórias $\lbrace X (t): t \in T \rbrace$ que tomam valores em $H$.
\end{dfn}

\begin{dfn}(\textbf{segunda ordem}). Um processo estocástico ${X(t) : t \in T}$  é de segunda ordem se $X(t) \in L^2(\Omega)$ para cada $t \in T$, a função de meia é definida como $\mu(t) = \mathbb{E}[X(t)]$ e a função de co-variância é definida por $C(s, t) = Cov(X(s), X(t))$  para todo $s,t \in T$.
\end{dfn}

\begin{dfn} (\textbf{processo Gaussiano}).  Um processo estocástico de segunda ordem ${X(t) : t \in T}$  é Gaussiano  se $X = [X(t_1), . . . , X(t_M)]$ segue uma distribuição gaussiana multivariada para cada $t_1, \cdots,t_M \in T$ e cada $M \in \mathbb{N}$.
\end{dfn}

\begin{dfn} (\textbf {movimento Browniano}). Dizemos que  $\lbrace W(t) : t \in \mathbb{R} \rbrace$ é um Movimento Browniano se é um processo gaussiano  com trajetórias continuas, meia $\mu(t) = 0$ e função de 
covariância $C(t,x) = min\lbrace t,x \rbrace$.
\end{dfn}

\begin{dfn} \label{PB} (\textbf {ponte Browniano}). Dado um movimento Browniano $W(t)$ chamamos de Ponte Browniano ao processo $B(t)$ em $[0,T]$ cuja distribuição é obtida condicionando a distribuição de $W(t)$ por condições de contorno em $t = 0$ e $t=T$.
\end{dfn}

\begin{dfn} (\textbf{ruído branco escalar}). Chamamos de \textit{ruído branco escalar} ao processo estocastico $\xi$ onde cada $\xi(t)$ é uma variável aleatória com distribuição normal. Além disso a função de covariância é: 
\begin{equation}
\langle \xi(t),\xi(t')\rangle = \delta(t-t').
\end{equation}
Uma maneira simples de descrever um ruído branco é
\begin{equation}\label{RB1D}
\zeta(t) = \sum_{j=1}^{\infty} \xi_j \varphi_j(t),
\end{equation}
onde $\xi_j$ são variáveis aleatórias com distribuição normal e $\lbrace \varphi_j \rbrace_{j\in\mathbb{R}}$ é uma base ortonormal de $L^2(\mathbb{R})$.
\end{dfn}

\begin{dfn} (\textbf{ruído branco bidimensional}). Chamamos de \textit{ruído branco} à distribuição $\zeta$ que toma valores em um campo gaussiano bidimensional (i.e. $ \zeta(t,x) \thicksim N(0,1)$)  com função de correlação,
\begin{equation}
\langle \xi(t,x),\xi(t',x')\rangle = \delta(t-t')\delta(x-x').
\end{equation}
Uma maneira simples de descrever um ruído branco é
\begin{equation}\label{RB2D}
\zeta(t,x) = \sum_{j=1}^{\infty} \xi_j(t) \varphi_j(x),
\end{equation}
onde $\xi_j$ são ruídos brancos escalares e $\lbrace \varphi_j \rbrace_{j\in\mathbb{R}}$ é uma base ortonormal de $L^2(\mathbb{R})$.
\end{dfn}

\begin{dfn} (\textbf{processo $Q$-Wiener}) Seja $(\Omega,\mathcal{F},\mathcal{F}_t,\mathbb{P})$ um espaço de probabilidade filtrado. Um processo estocástico $\lbrace W(t) : t\geq 0\rbrace$ que toma valores em $U$ é um processo $Q$-Wiener se:
\begin{itemize}
\item $W(0) = 0$ a.s. \footnote{a.s. significa \textit{almost sure}. Escolhemos neste trabalho manter a abreviatura a.s. ao longo do todo o texto.}
\item $W(t): \mathbb{R}_+ \times \Omega \longrightarrow U$ é uma função continua para cada $\omega \in \Omega$
\item $W(t)$ é $\mathcal{F}_t$ - adaptado e $W(t) - W(s)$ é independente de $\mathcal{F}_s$ para cada $s\leq t$ e 
\item $W(t) - W(s) \sim N(0,(t-s)Q)$ para todo $0 \leq s \leq t$.
\end{itemize}\label{QWdef}
\end{dfn}

 Podemos provar que o processo $Q$-Wiener $W$ pode ser escrito como:
 \begin{equation}
 W_t(x) = \sum_{j=1}^{\infty} \sqrt{q_j} \beta_j(t) \chi_j(x),\label{QW}
 \end{equation}
onde $\beta_j$ são movimentos brownianos escalares \cite{JEK11,LPS14}. Dizemos que um processo $Q$-Wiener é \textbf{cilíndrico} quando $Q \equiv I$. 

\begin{dfn}
(\textbf{equação semilinear}). É uma equação estocástica da forma:
\begin{equation}
dX_t = \left[ AX_t + F(X_t) \right]dt + G(X_t)dW_t, \label{semiL}
\end{equation}
onde $W_t$ é um processo $Q$-Wiener e os operadores $F:H \longrightarrow H$ e $G: H \longrightarrow H$ são, em geral, não lineares.
\end{dfn}

\begin{dfn}
(\textbf{solução forte}) Um processo previsível $X:[0,T] \longrightarrow H$ é chamado solução forte de (\ref{semiL}) quando $X_t$ satisfaz:
\begin{align}
dX_t = \int_0^t\left[ A X_s + F(X_s) \right]ds + \int_0^t G(X_s)dW_s, && \forall t\geq 0.\label{forte}
\end{align}
\end{dfn}

\begin{dfn}
(\textbf{solução fraca}) Um processo previsível $X:[0,T] \longrightarrow H$ é chamado solução fraca de (\ref{semiL}) quando $X_t$ satisfaz:
\begin{align}
\langle X_t,v \rangle = \langle X_0,v \rangle + \int_0^t\left[ \langle A X_s,v \rangle + \langle F(X_s),v\rangle \right]ds + \int_0^t \langle G(X_s)dW_s, v\rangle, && \forall v\in H.
\end{align}
e de (\ref{QW}) temos que a segunda integral pode ser escrita como:
\begin{equation}
\int_0^t \langle G(X_s)dW_s,v\rangle = \sum_{j=1}^{\infty} \int_0^t \langle  G(X_s)\sqrt{q_j}\chi_j,v \rangle d\beta_j(s).\label{fraca}
\end{equation}
\end{dfn}

\begin{dfn}
(\textbf{solução \textit{mild}}) Um processo previsível $X:[0,T] \longrightarrow H$ é chamado solução \textit{mild} de (\ref{semiL}) quando $X_t$ satisfaz:
\begin{align}
X_t = e^{tA} X_0 + \int_0^t e^{(t-s)A} F(X_s)ds + \int_0^t e^{(t-s)A} G(X_s)dW_s.
\end{align}
\end{dfn}\label{leve}
%
%

\chapter{Uma revisão do modelo Kardar-Parisi-Zhang}\label{cap2}

\section{Modelos de crescimento}

Os fenómenos crescimento de superfícies são influenciados por muitos fatores, a maioria deles indistinguíveis se levadas em conta apropriadas escalas de tempo e espaço. Porém, com base na física (observação e entendimento do modelo real físico que ocorre) os cientistas sempre esperam que exista um pequeno número de leis básicas fundamentais que os determinam, ou melhor, que caracterizam de forma única e geral a morfologia e a dinâmica de crescimento. Com este fim, o estudo do modelo de deposição balística tem ajudado a encontrar as propriedades essenciais de várias classes de fen\^omenos de crescimento, ver e.g., \cite{LBGG97,MRS86}.\\
\subsection{Deposição balística (\textit{balistic deposition} (BD))}

\begin{figure}[h]
\centering
\includegraphics[scale=.6]{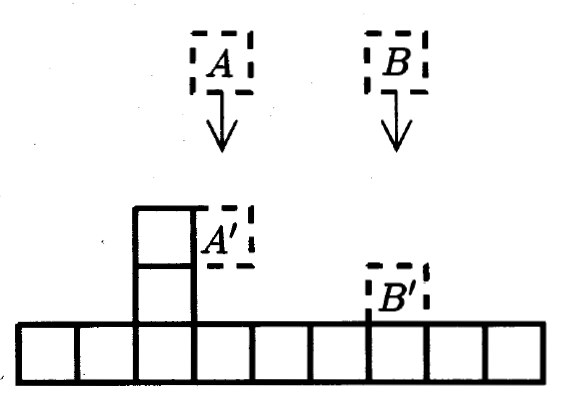} 
\caption{Modelo discreto de deposição balística. Extraído de \cite{BS95}.} \label{fig:BD1}
\end{figure}

   No modelo de deposição balística (BD), uma partícula cai de algum ponto aleatório de uma altura maior que a altura máxima da superfície e segue uma trajetória vertical até chegar à superfície onde fica aderida. Nesta versão simples do modelo, as partículas que chegam se fixam à primeira partícula que tocam. Assim, em um primeiro momento, estamos considerando a superfície inicial como sendo plana e com longitude $L$. Definimos a superfície como sendo o conjunto de partículas que ocupam a maior altura em cada coluna. A Figura \ref{fig:BD1} fornece uma ideia geométrica dessa modelagem.
  
Duas quantidades de interesse em nosso modelo são {\it a altura média} e {\it a rugosidade} \cite{Qua12,KS92,BS95}. A altura média pode ser considerada como a média aritmética de todas as alturas, enquanto a rugosidade é a soma dos desvios médios quadráticos entre as alturas e a altura média, dividida pelo comprimento do intervalo,
\begin{align}
&\hat{h}(t) = \frac{1}{L} \sum_{i=1}^n h(i,t),  \\
&w(L,t)    = \sqrt{\frac{1}{L} \sum_{i=1}^n \left[ h(i,t) - \hat{h} (t) \right]^2}. 
\end{align}

\begin{figure}[h!]
\centering
\includegraphics[scale=.5]{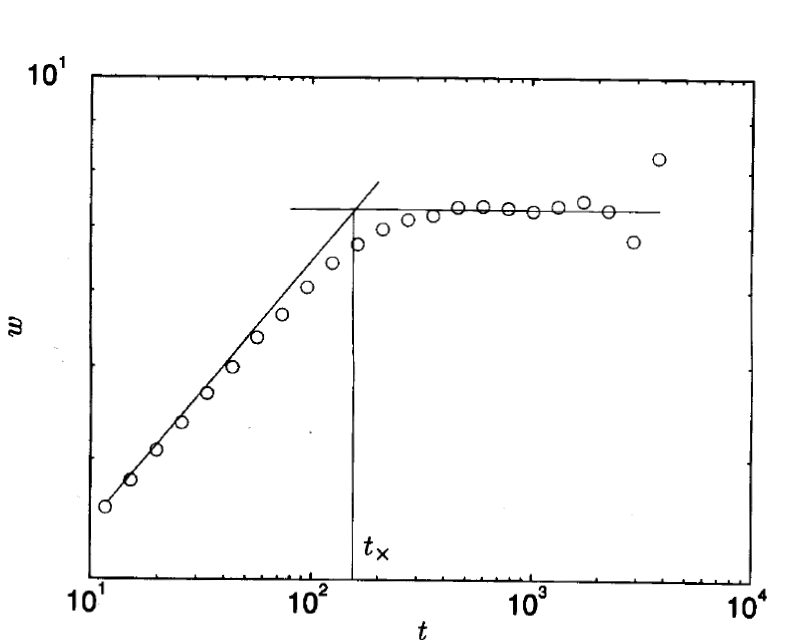}
\caption{Modelo simples de deposição balística. Extraído de \cite{BS95}.} \label{fig:Rug1}
\end{figure}

Se o ritmo de deposição (chegada de partículas à superfície por unidade de tempo) é constante, não é difícil demonstrar (ver \cite{BS95}) que o ritmo médio de crescimento é,
\begin{equation}
\hat{h}(t) \sim t.
\end{equation}
Os resultados experimentais indicam que existe um tempo \textit{crossover}. A saber, um tempo \textit{crossover} $t_x$ que separa dois regimes (ver Figura \ref{fig:Rug1}) onde a rugosidade tem comportamentos distintos dados pelas equações (\ref{Rug_Exp1}\, e\, \ref{Rug_Exp2}),
\begin{align}
w(L,t) \sim t^\beta, && t \ll t_x \label{Rug_Exp1}, \\
w(L,t) \sim L^\alpha, && t \gg t_x \label{Rug_Exp2}.
\end{align}

\begin{figure}[h!]
\centering
\includegraphics[scale=.5]{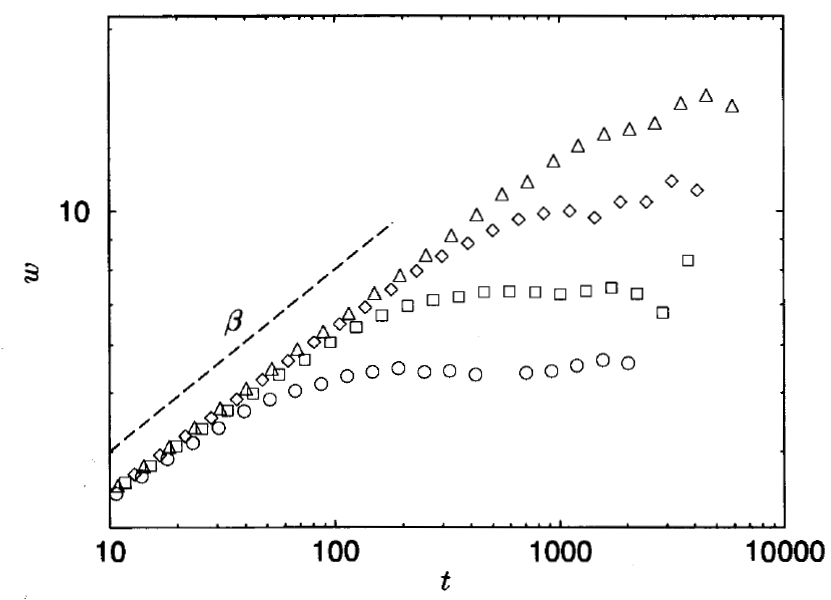} 
\caption{Efeito de crossover e rugosidade. Resultados experimentais para o modelo BD. Extraído de \cite{BS95}.} \label{fig:Rug2}
\end{figure}

Quando $L$ varia, também varia o tempo de {\it crossover}  ao estado de saturação e então tem-se,
\begin{align*}
&t_x \sim L^z \numberthis. \label{Tx}
\end{align*}
Na Figura \ref{fig:Rug2} se mostra a variação da rugosidade no tempo para distintos valores de $L$.
Os expoentes $\alpha$, $\beta$ e $z$ não são independentes e podemos provar que a rugosidade cumpre a chamada relação de escala de \cite{FV85}, ou seja, 
\begin{equation}\label{RelEsc}
w(L,t) \sim L^\alpha f\left(  \frac{t}{L^z} \right),
\end{equation}
onde a função $f(u) \sim u^\beta$ para $u \ll 1$ e $f(u) = constante$ para $u \gg 1$ (o parâmetro $u$ representa o termo $t/t_x$). Aproximando o ponto \textit{crossover} (ver também e.g., \cite{GJ12,BFRT93,KAT13}) $(t_x,w(L,t_x))$ pela esquerda, segundo (\ref{Rug_Exp1}), temos que $w(t_x) \sim t^\beta_x$ e, aproximando pela direita e usando (\ref{Rug_Exp2}), temos que $w(t_x) \sim L^\alpha$. Portanto, temos que $t_x^\beta = L^\alpha$ e utilizando (\ref{Tx}) obtemos a lei de escala,
\begin{equation}\label{LeiEsc}
z = \frac{\alpha}{\beta}.
\end{equation}

\begin{figure}[h!]
\centering
\includegraphics[scale=.6]{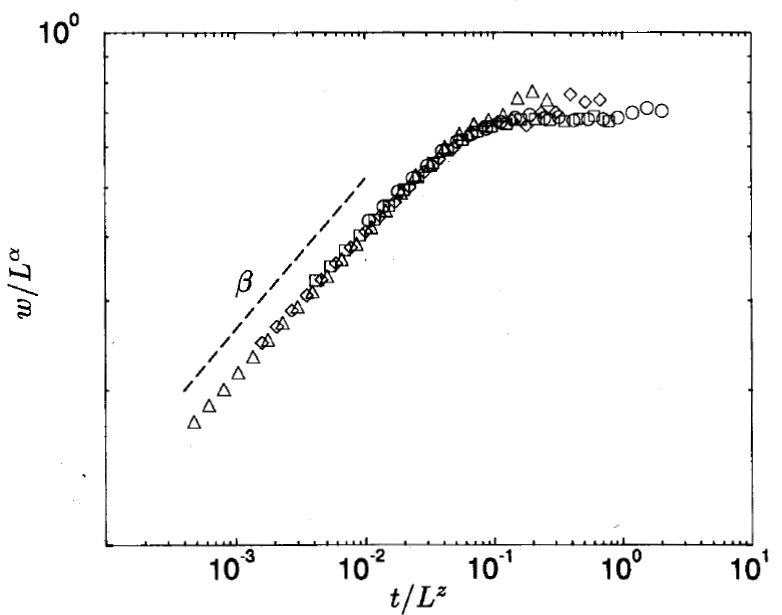}
\caption{Reescalamento que leva à colisão das curvas de rugosidade obtidas para quatro escolhas diferentes de $L$ no modelo BD. Extraído de \cite{BS95}.} \label{fig:Rug3}
\end{figure}

As curvas de rugosidade e os tempos de \textit{crossover}, mostrados na Figura \ref{fig:Rug2}, agora colapsam na curva $f$, considerando o reescalamento (\ref{RelEsc}), como mostra a Figura \ref{fig:Rug3}.

A lei de escala (\ref{LeiEsc}) é válida para todo processo de crescimento que obedece a relação de escala (\ref{RelEsc}). A Tabela \ref{tab:Somario} sumariza os conceitos expostos nesta seção. Cumpre mencionar que assim temos em mãos uma forma simples e efetiva de testar o desempenho do método proposto, ou seja, se as aproximações estiverem, de fato, corretas então seremos capazes de recuperar essas curvas de crescimento a partir do pós-processamento das soluções numéricas calculadas. É claro que isso não caracteriza qualquer forma rigorosa de demonstração matemática. Porém teremos um bom argumento formal para dar um suporte de motivação para perseguir estudos mais avançados.

\begin{table}
\centering
\begin{tabular}{ccc}
\hline
\hline
Média & $\bar h(t) = \frac{1}{L}\sum_{i=1}^L h(i,t)$ \\
Rugosidade  & $w(L,t) = \sqrt{\frac{1}{L} \sum_{i=1}^L [h(i,t) - \bar{h} (t)]^2}$\\
\hline
Expoente de crescimento & $w(L,t)\sim t^{\beta}$ & $[t \ll t_x]$\\
Expoente de rugosidade & $w_{sat}(L)\sim L^{\alpha}$ & $[t \gg t_x]$\\
Expoente dinâmico & $t_x \sim L^z$\\
\hline
Relação de escala & $w(L,t) \sim L^{\alpha} f\left(t/L^z \right)$\\
lei de escala & $z = \alpha/\beta$\\
\hline 
\hline
\end{tabular}
\caption{Resumo das principais grandezas associadas ao modelo de deposição balística que serão utilizadas neste trabalho. Extraído de \cite{BS95}.}
\label{tab:Somario}
\end{table}

\subsection{Correlações}

  Uma propriedade a ressaltar dos modelos \textbf{BD} é a presença de  correlações sobre a superfície em crescimento. Isso significa que os locais na superfície não evoluem de forma independente, mas dependem das suas vizinhanças. Cada nova partícula que chega à interface se fixa na primeira partícula que encontra e isto faz com  que o crescimento tenha uma componente lateral como mostra a Figura \ref{fig:BD1}. Embora este processo de crescimento seja local, a informação sobre a altura em cada local propaga-se globalmente. A distância sobre a qual as alturas estão relacionadas é chamada de comprimento de correlação e denota-se por $\xi_\Vert$. Sendo que os pontos da superfície no início do processo não estão correlacionados, mas começam a ficar mais e mais relacionados, e dado que $\xi_\Vert$ é limitado por $L$, podemos deduzir que quando $\xi_\Vert = L$, então toda a superfície está correlacionada e o processo chegou ao equilíbrio, i.e.,
\begin{align}
\xi_\Vert \sim L, && t\gg t_x.
\end{align}

Nesse contexto, sabendo que no equilíbrio $t_x \sim L_z$, e substituindo $L$ por $\xi_\Vert$ obtemos,
\begin{align}
\xi_\Vert \sim t^{1/z} && t\ll t_x.
\end{align}  

\subsection{Deposição Aleatória (\textit{random deposition} (RD))}

\begin{figure}
\centering
\includegraphics[scale=.7]{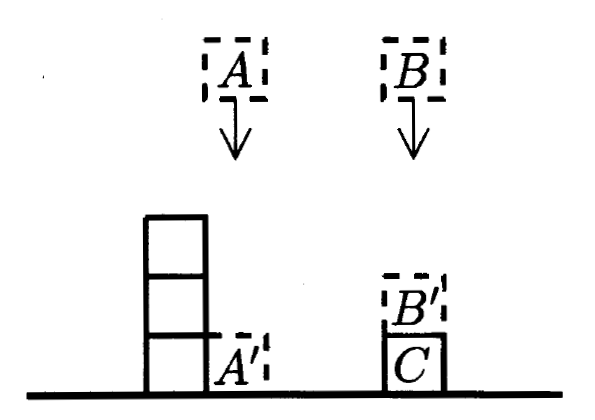} 
\caption{Modelo discreto de deposição aleatória. Extraído de \cite{BS95}.}
\label{fig:RD}
\end{figure}

O modelo de deposição aleatória (RD) é muito simples, mas é muito útil para introduzir a ideia de associar a um modelo discreto uma equação no contínuo. O ideia mecânica do modelo RD é similar do modelo \textbf{BD} só que neste caso as partículas descem até fixar na maior altura na direção vertical de descida, como mostra a Figura \ref{fig:RD}.

  Neste caso todas as posições na superfície são não correlacionadas e cada altura cresce com probabilidade $p = 1/L$. A probabilidade de uma coluna ter altura $h$, após a deposição aleatória de $N$ partículas, pode ser calculada explicitamente. Assim, isso explica o crescimento da rugosidade para este modelo, i.e., que cresce indefinidamente dado por $w(t) \sim t^{\beta}$ com $\beta = 1/2$. A independência entre os crescimentos das alturas em cada posição são a causa de que o processo não apresenta um estado de saturação e, portanto, os coeficientes $\alpha$ e $z$ ficam indefinidos.\\
  Vamos associar uma {\it equação estocástica contínua} que descreve o modelo. Como o modelo é discreto, vamos considerar um amolecimento da superfície de tal forma que tenha sentido esta aproximação para escalas pequenas. 
  Atendendo às caraterísticas do modelo RD teremos equação geral,
\begin{equation}
\partial_t h(x,t) = \Phi(x,t),
\end{equation}  
onde $\Phi$ é o número de partículas que chegam na posição $x$ no tempo $t$, levando em conta que o processo é aleatório, podemos reescrever o modelo, 
\begin{equation}\label{Eq:RD}
\partial_t h(x,t) = F + \xi(x,t),
\end{equation}  
onde $F$  é o número médio de partículas que chegam no ponto $x$ e $\xi(x,t)$ é a aleatoriedade. Em geral $\xi$ é um ruído branco, ou seja, uma função que para cada ponto $(x,t)$ tem distribuição normal com média zero e variância,
\begin{equation*}
\mathbb{E}\left[\xi(x,t)\xi(y,s) \right] = \delta\left(x-y \right) \delta\left(t-s \right),
\end{equation*}
A solução da equação (\ref{Eq:RD}) pode ser calculada explicitamente por,
\begin{align}
h(t,x) &= Ft + \int_0^t \xi(\tau,x) d\tau,  \\
       & = Ft.
\end{align}  
Os momentos de $h(t,x)$ são $\langle h(t,x)\rangle = Ft$ e $\langle h^2(t,x)\rangle = F^2 t^2 + 2Dt$ de onde temos que, 
\begin{equation}
w^2(t,x) = \langle h^2 \rangle - \langle h\rangle^2 = 2Dt.
\end{equation}  
  Assim, chegamos ao mesmo expoente de escala $\beta = 1/2$.
  
\subsection{Deposição aleatória com difusão}  

  No modelo RD as partículas que chegavam à  superfície se fixavam no ponto de máxima altura na direção de decrescimento. Neste novo modelo permitiremos que as partículas se difundam  até atingir o ponto de menor altura, onde permanecerão fixadas. Este modelo com difusão faz com que a superfície seja mais suave. Além disso, é claro que os pontos da superfície estarão mais correlacionados. É valido mencionar, ainda, que 
resultados experimentais reportados em  (\cite{Fam86}) mostram que os coeficientes de crescimento para este modelo são,
\begin{align}
\alpha = 0.48 \pm 0.02, && \beta = 0.24 \pm 0.01. \label{Exp_EW_exper}
\end{align}

  Da mesma forma que fizemos no exemplo anterior, vamos associar o modelo discreto com uma equação estocástica contínua. Neste caso a equação geral toma a forma,
\begin{equation}
\partial_t h(x,t) = G(t,x,h) + \xi(x,t).
\end{equation}  

  Para a derivação da equação, vamos utilizar algumas leis do crescimento que deve cumprir a solução (ver \cite{BS95}),
como por exemplo,
\begin{enumerate}
\item Invariância por translações no tempo.
\item Invariância por translações na direção do crescimento.
\item Invariância por translações na direção perpendicular ao crescimento.
\item Simetria de inversão e rotação com respeito ao eixo do crescimento.
\end{enumerate}  

  Fazendo expansão em série de Taylor de $h(t,x)$ e levando em conta as condições de simetria antes mencionadas, podemos eliminar alguns termos da série. Também não levamos em conta os termos de ordem inferior, ficando apenas com os mais relevantes para obter a chamada equação de Edward-Wilkinson (EW) introduzida em (\cite{EdW82} e \cite{CW78}), i.e., 
\begin{equation}\label{Eq:EW}
\partial_t h(t,x)= \nu \triangledown^2 h(t,x) + \xi(t,x).
\end{equation}

  A equação (\ref{Eq:EW}) apenas tem sentido para valores pequenos de $\triangledown h$. O termo difusivo  $\triangledown^2 h$ provoca a reorganização das partículas na interface compensando de certo modo os efeitos da aleatoriedade provocada pelo ruído $\xi$. No caso em que a superfície estiver se movendo com uma certa velocidade $v$, o termo difusivo deve ser agregado à equação,
\begin{equation}\label{Eq:EW_vel}
\partial_t h(t,x) = v + \nu \triangledown^2 h(t,x) + \xi(t,x).
\end{equation}

  Para achar os coeficientes $\alpha$, $\beta$ e $z$ deste processo podemos utilizar o argumento de escala ou simplesmente resolver a equação via transformada de Fourier. Vamos optar pela primeira forma. Consideremos a transformação de escala,
\begin{subequations}  
\begin{align*}
x \longrightarrow x' = bx, \\
h \longrightarrow h' = b^\alpha h, \\
t \longrightarrow t' = b^z t.
\end{align*}
\end{subequations} 

  Substituindo as relações anteriores em (\ref{Eq:EW}), e levando em conta a forma que tais relações ``escalam'' com o ruído branco, obtemos:
\begin{equation}\label{EW_esc}
\partial_t h = \nu b^{z-2} \triangledown^2 h + b^{\frac{z-d}{2}-\alpha}\xi.
\end{equation}  

  Para garantir a invariância, a equação (\ref{EW_esc}) não pode depender de $b$. Isto define a escolha dos coeficientes,
\begin{align}\label{Exp_EW_Cal}
\alpha = 1/2,  && \beta = 1/4,  && z = 2.
\end{align}  

A similaridade dos mecanismos de difusão entre o modelo e a solução da equação de Edward- Wilkinson (\ref{Eq:EW}) mostram, junto à  similaridade dos expoentes achados experimentalmente (\ref{Exp_EW_exper}) e os calculados (\ref{Exp_EW_Cal}),  que o modelo e a equação EW pertencem à  mesma {\it classe de universalidade}, que é diferente da classe do modelo RD.\\

\subsection{Kardar-Parisi-Zhang}

  Até agora identificamos duas classes de universalidade diferentes para modelos de crescimento: a classe a que pertence o modelo RD e a função estocástica asociada (\ref{Eq:RD}), e a classe à que pertence o modelo de deposição aleatória com difusão e a equação estocástica associada (\ref{Eq:EW}). Para o modelo balístico as simulações numéricas realizadas em \cite{CV92} e \cite{MRS86} sugerem as seguintes predições para os expoentes:
\begin{align}\label{Exp_KPZ_exp}
\alpha = 0.47 \pm 0.02, && \beta = 0.33 \pm 0.006.
\end{align}

\begin{figure}[h!]
\centering
\includegraphics[scale=.5]{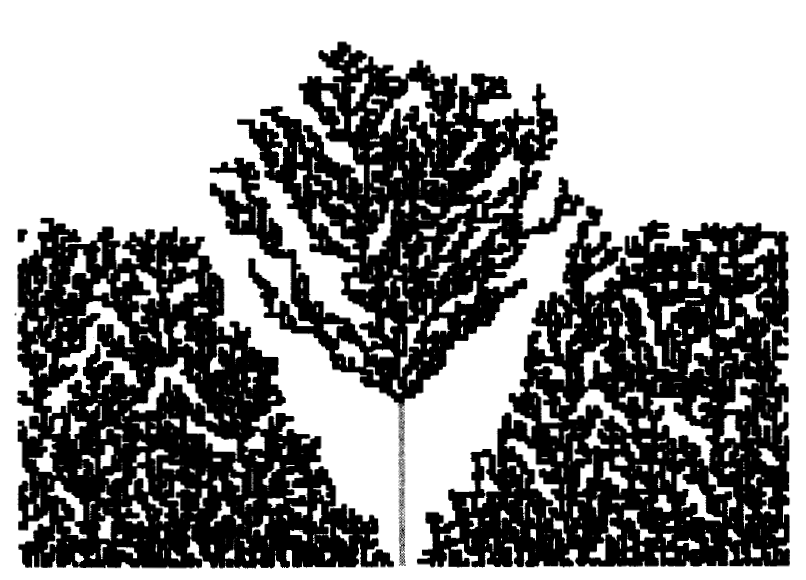}
\caption{Crescimento lateral em deposição balística. Extraído de \cite{BS95}.} 
\label{fig:BD_est}
\end{figure}

As  predições anteriores indicam que este modelo pertence a uma terceira classe de universalidade e que, portanto, não deve ser descrito pelas equações associadas aos modelos anteriores RD ou EW. O modelo de deposição balística apresenta, como característica que o distingue dos anteriores, o crescimento lateral. Na Figura \ref{fig:BD_est} se mostra o efeito do crescimento lateral de uma superfície que cresce por deposição balística.

\begin{figure}[h!]
\centering
\includegraphics[scale=1.2]{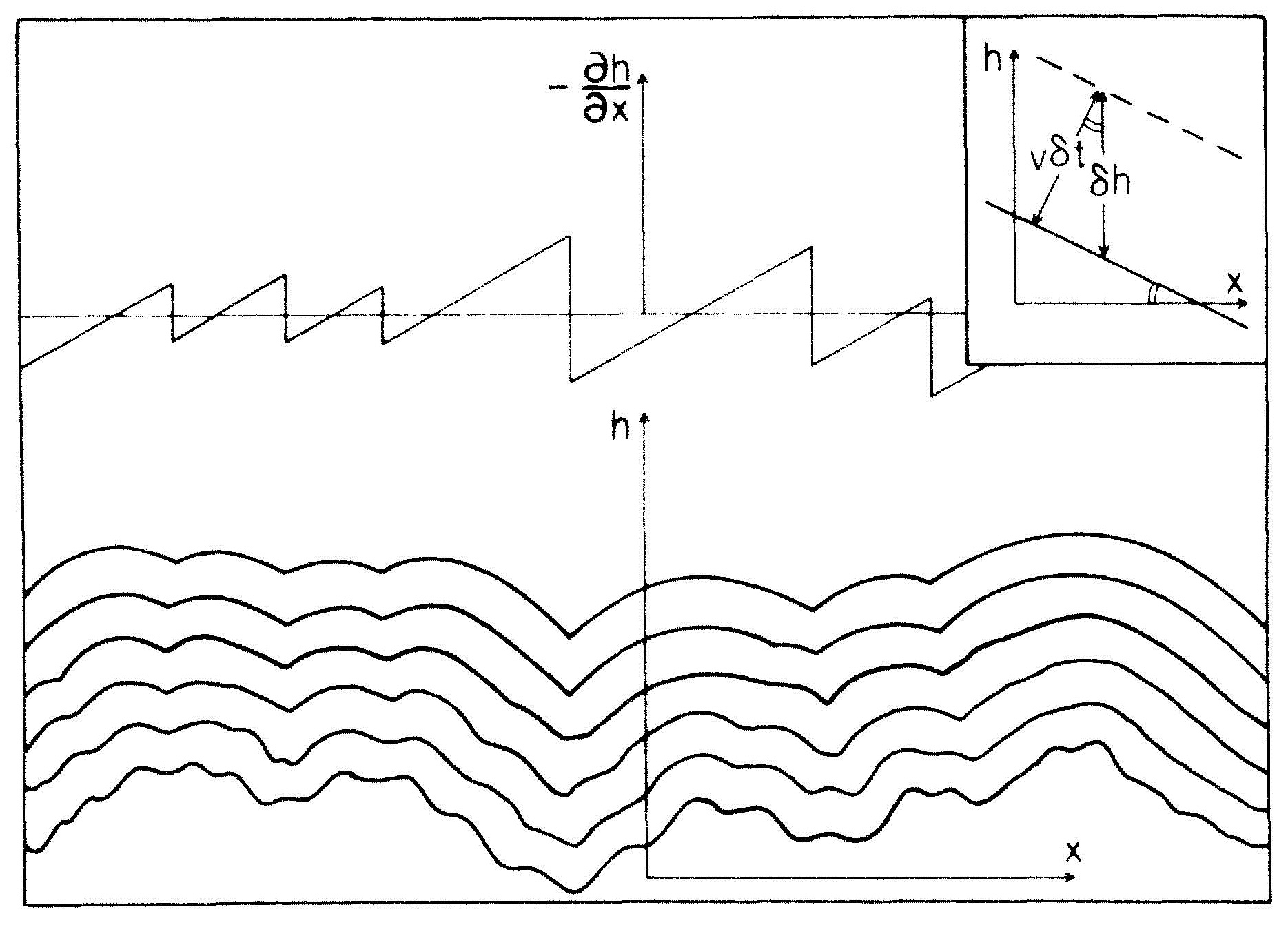}
\caption{Ideia geométrica da derivação do termo não linear. Extraído de \cite{KPZ86}.}
\label{fig:1} 
\end{figure}

No caso de escalas muito grandes, onde o modelo pode ser analisado no contínuo, o crescimento acontece localmente na direção do vetor normal à superfície. Precisaremos de uma nova equação estocástica associada a este modelo. Para tal fim vamos generalizar a equação EW. Uma forma de derivar o termo que determina o crescimento lateral é considerando uma linearização da superfície $h$ em torno do ponto $x$. Dada uma pequena variação no tempo $\Delta t$ poderíamos aproximar a variação da altura $\delta h$ no ponto $x$ sabendo que a variação no sentido do vetor normal é $\lambda\delta t$. A Figura \ref{fig:1} mostra as considerações geométricas para obter a relação
\begin{equation*}
\delta h = \left[(\lambda\delta t)^2 + (\lambda\delta t \partial_x h)^2 \right]^{1/2} = \delta t\left[(\lambda)^2 + (\lambda \partial_x h)^2 \right]^{1/2}\!\!\!\!\!.
\end{equation*}

Dividindo a equação anterior por $\delta t$ e tomando o limite quando $\delta t \longrightarrow 0$ chegamos em $\partial_t h = \lambda \left( 1 + \partial_x h \right)^{1/2}$. Agora fazemos expansão em série de Taylor da função $(1+u)^{1/2}$ em torno do ponto zero e a avaliamos  em $ \partial_x h $. Assim temos que $\partial_t h = (1 + (\partial_x h)^2)^{1/2} = \lambda + (\lambda/2)\cdot (\partial_x h)^2$.\\
  
\begin{figure}[h!]
\centering
\includegraphics[scale=.5]{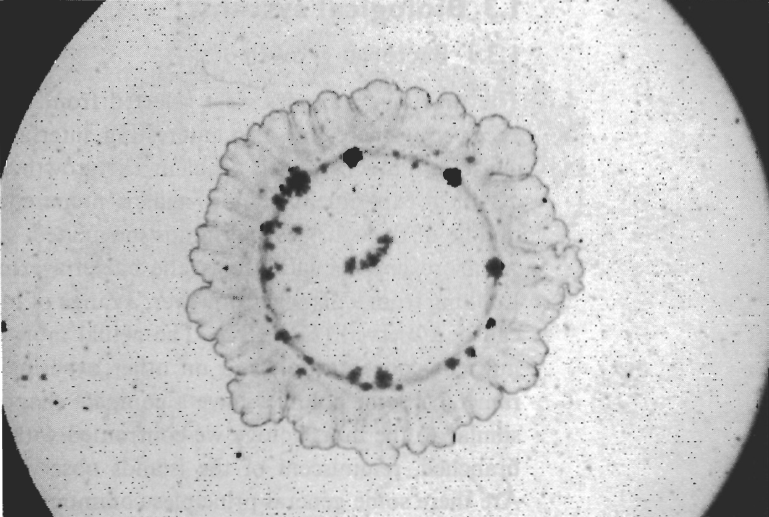}
\caption{Crescimento de uma colónia de bactérias com rugosidade compacta. Extraído de \cite{BSC94}.}
\label{fig:bact}
\end{figure}

Isto sugere a presença do termo $(\partial_x h)^2$ na nova equação refletindo o crescimento lateral. Acrescentando este termo na equação EW obtemos a equação KPZ \cite{KPZ86}, \\
\begin{equation}\label{KPZ}
\partial_t h = \nu \partial_x^2 h + \frac{\lambda}{2}\left( \partial_x h \right)^2 + \sqrt{D} \xi.
\end{equation} 
Claro que há um problema com esta derivação. Para que a equação faça sentido, o
valor de $\partial_x h$ tem que estar perto de zero, mas na realidade este valor é muito
grande. Então, nós teríamos que subtrair um termo muito grande que reflete as pequenas 
variações de escala. 
\newpage
Mesmo assim, usando esta ``inocente derivação'', encontramos campos 
não triviais \cite{KS92}, \cite{HHZ95} e \cite{BS95}. \\
\\
    A equação KPZ rapidamente tornou-se um modelo protótipo 
    fundamental para a modelagem matemática da dinâmica de crescimento de diversos processos em 
    diferentes áreas \cite{LBGG97,FFTV85,Hai11}. Em Biologia, a KPZ é útil para modelar o crescimento 
    de colónias de bactérias como foi mostrado nos trabalhos \cite{BSC94} e \cite{MF90}. 
    A Figura \ref{fig:bact} ilustra um exemplo do crescimento dessas colônias.\\
    
\begin{figure}[h!]
\centering
\includegraphics[scale=.5]{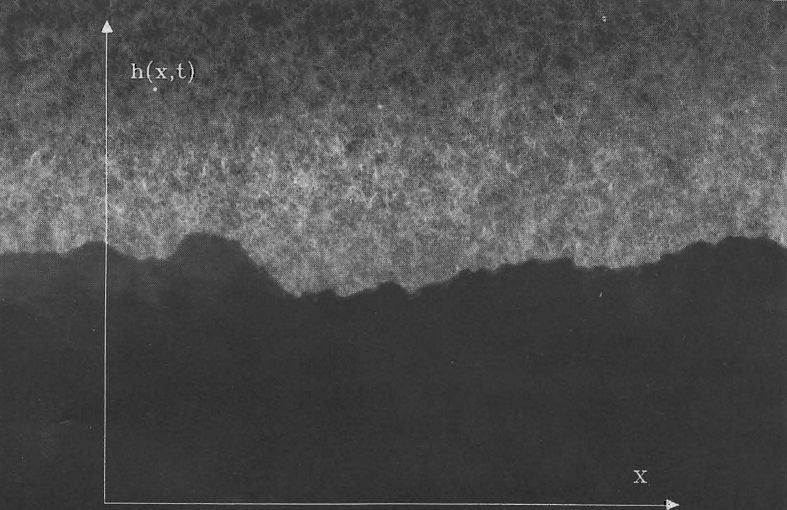} 
\caption{Foto ilustrando um segmento de uma folha de papel queimando. A parte transversal mede cerca de 9 cm. Extraído de \cite{ZZ92}.}
\label{fig:burn}
\end{figure}

   A equação KPZ também é útil na modelagem matemática de frentes de fogo. Em \cite{ZZ92}, 
   os resultados de um experimento físico em laboratório, baseado na queima de uma folha de papel, 
   ajudam a predizer o comportamento de um incêndio florestal. A Figura \ref{fig:burn} ilustra a
    evolução da frente de fogo em experimento de queima de uma folha de papel.\\
  
\begin{figure}[h!]
\centering
\includegraphics[scale=.5]{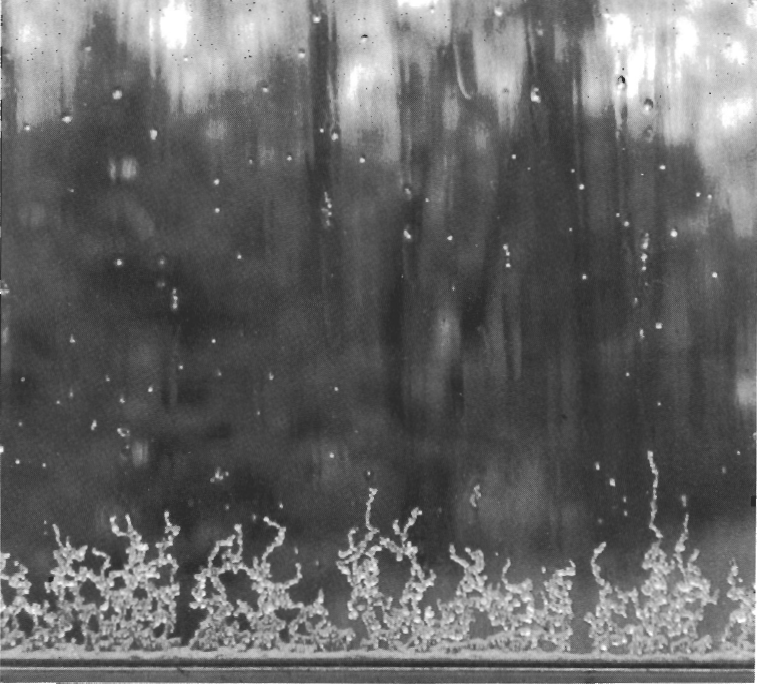}
\caption{Partículas de neve que caem sobre uma janela de cristal. Extraído de \cite{CV92}.} 
\label{fig:Atom1}
\end{figure}

\begin{figure}[h!]
\centering
\includegraphics[scale=.5]{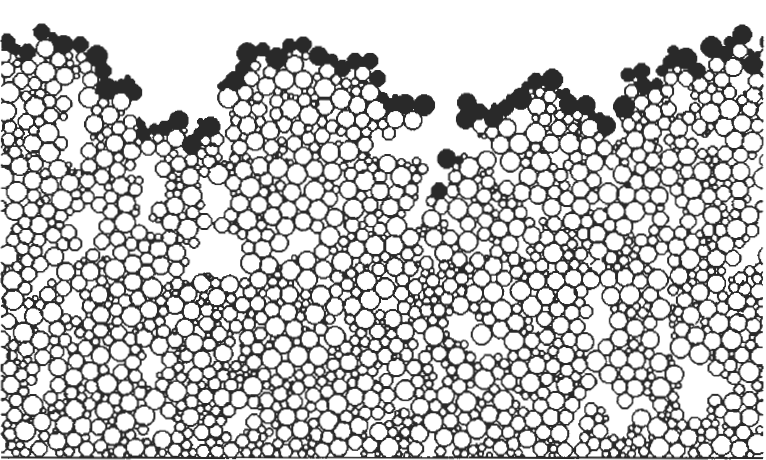} 
\caption{Partículas esféricas com diâmetros aleatórios, normalmente distribuídos, que chegam a uma superfície e rolam até a colisão com outras duas partículas. Extraído de \cite{Leo91}.} 
\label{fig:Atom2}
\end{figure}

  Na Física também existem diversos fen\^omenos de crescimento que podem ser modelados pela equação KPZ, 
  por exemplo, a deposição de partículas sobre uma superfície com certa geometria. 
  Em \cite{CV92} e \cite{Leo91}, encontramos resultados experimentais que corroboram com o anterior.
   As Figuras \ref{fig:Atom1} e \ref{fig:Atom2} ilustram alguns destes fenómenos físicos.\\
  
\newpage
\begin{center}
\textbf{Coeficientes da equação KPZ}
\end{center}

Pode ser demonstrado (\cite{BS95}) que os coeficientes $\alpha$, $\beta$ e $z$ são coerentes com os valores obtidos nos experimentos (\ref{Exp_KPZ_exp}). Mostraremos que podemos reduzir a analise à escolha particular  $\lambda = 2$, $\nu = 1$ e $D = 1$ para a qual a equação KPZ toma a forma 
\begin{equation}\label{KPZ_Canon}
\partial_t h =  \left( \partial_x h \right)^2 + \partial_x^2 h + \xi,
\end{equation} 
a qual chamaremos de equação KPZ canônica. Consideremos a mudança de escala,
\begin{equation}
h_{\epsilon}(x,t) = \epsilon^{\beta} h(\epsilon^{-z} t,\epsilon^{-1} x),
\end{equation} 
de onde temos $\partial_t h = \epsilon^{z-\beta}\partial_t h_{\epsilon}$, 
$\partial_x h = \epsilon^{1-\beta} h_{\epsilon}$ e
$\partial_x^2 = \epsilon^{2-\beta} h_{\epsilon}$. O ruído branco também rescala,
\begin{equation}
\xi(t,x) \overset{dist}{=} \epsilon^{\frac{z+1}{2}}\xi(\epsilon^{-z} t, \epsilon^{-1} x),
\end{equation} 
onde a igualdade significa que os dois campos aleatórios têm a mesma distribuição. Substituindo as relações anteriores em (\ref{KPZ_Canon}) e dividimos por $\epsilon^{z-\beta}$ obtemos,
\begin{equation}\label{KPZ_Scaling}
\partial_t h_\epsilon =  \epsilon^{2-z-\beta} \left( \partial_x h \right)^2 + \epsilon^{2-z} \partial_x^2 h + \epsilon^{\beta-\frac{1}{2}z - \frac{1}{2}} \xi.
\end{equation}

Claramente podemos tomar agora $\lambda = 2\cdot\epsilon^{2-z-\beta}$, $\nu = \epsilon^{2-z}$ e
$D = \epsilon^{\beta-\frac{1}{2} z - \frac{1}{2}}$ e assim recuperamos a equação KPZ (\ref{KPZ}) a partir da equação KPZ canônica (\ref{KPZ_Scaling}).\\

Podemos agora então tentar achar um reescalamento,
\begin{equation}
h_{\epsilon}(t,x) = \epsilon^\beta h(\epsilon^{-z}t, \epsilon^{-1} x),
\end{equation}
mediante o qual esperamos observar algum comportamento não trivial para a equação (\ref{KPZ_Scaling}) quando $\epsilon\longrightarrow 0$ em escalas grandes para o tempo e o espaço. \\

Se fixamos $t = 0$, o fato de que a solução seja localmente Browniana \cite{Qua12}, faz com que,
\begin{equation}
\beta = 1/2.
\end{equation}
  
Então a equação fica,
\begin{equation}
\partial_t h_\epsilon = \epsilon^{3/2-z} \left( \partial_x h \right)^2 + \epsilon^{2-z} \partial_x^2 h + \epsilon^{-\frac{1}{2}z} \xi.
\end{equation}
 
Para evitar a divergência no termo não linear tomamos,
\begin{equation}
 z = \frac{3}{2}.
\end{equation}

\subsection{Solução de Hopf-Cole}\label{HC}

Percebemos que a equação (\ref{KPZ_Canon}) tem um grande problema. O termo não-linear pode não fazer sentido. Vemos que o termo não-linear precisa de um tipo de renormalização não finita. Portanto, seria mais honesto escrever a equação como,
\begin{equation}\label{KPZ_R}
\partial_t h = \left[\left( \partial_x h \right)^2 - \infty \right]+ \partial_x^2 h + \xi.
\end{equation}  

    Em \cite{BG97} os autores propuseram que a solução correta da equação KPZ poderia ser obtida da seguinte forma: A equação estocástica do calor com ruído branco multiplicativo é:
\begin{equation}\label{Heat}
\partial_t z = \partial^2_x z + \xi z.
\end{equation}

Esta equação (\ref{Heat}) deve ser interpretada no sentido de Itô, em tal caso o problema de valor inicial e de contorno é bem posto e para um dado inicial razoável $z_0(x)>0$ teremos que
$z(t,x)>0$ $\forall x$. Bertini e Giacomi propuseram que,
\begin{equation}\label{HC_KPZ}
h(t,x) = \log z(t,x),
\end{equation}  
é a solução correta da equação KPZ. Existem várias razões que apoiam a anterior:
  
\begin{enumerate}
\item[(1)]  Se $\xi$ for uma função suave, (\ref{HC_KPZ}) seria a solução de (\ref{KPZ_Canon}). Para verificar isto basta substituir (\ref{HC_KPZ}) em (\ref{Heat}). Isto é chamado, simplesmente, a transformação de Hopf-Cole.

\item[(2)]  No caso que $\xi$ for um ruído branco, poderíamos amolecer a solução $z(t,x)$ de (\ref{Heat}) utilizando como \textit{mollifier} $\mathcal{G}_k(x) = \frac{1}{\sqrt{2\pi k^2}}\exp\lbrace -x^2/2k^2\rbrace$, ou seja, 
\begin{equation}
z_k(t,x) = \langle\mathcal{G}_k,z(t)\rangle= \int z(t,y) \mathcal{G}_k(x-y)\mathbb{d}y.
\end{equation}

Definamos $h_k(t,x) = -\log  z_k(t,x)$.  Então pela formula de Itô temos,
\begin{equation}
\partial_t h_k + \left(\partial_x h_k\right)^2 - \partial^2_x h_k - \xi = \{ z_k^{-1}\langle \mathcal{G}_k z,\xi\rangle - \xi\} + \frac{1}{2}z_k^{-2}\langle\mathcal{G}_k^2,z^2\rangle.
\end{equation}

Para o primeiro termo, podemos calcular $\mathbb{E}\left[\left(\int \int\varphi(t,x)  \{ z_k^{-1}\langle \mathcal{G}_k z,\xi\rangle - \xi\}\mathbb{d}x\mathbb{d}t \right)^2\right]$ para uma função $\varphi$ suave com suporte compacto utilizando a isometria de Itô,
\begin{equation*}
\int \int \mathbb{E} \left[ \left( \int \frac{\varphi(t,y)\mathcal{G}_k(x-y) d y}{\int \mathcal{G}_k(y-y')z(t,y') d y'}z(t,s) - \varphi(t,x) \right)^2 \right]dx\,dt,
\end{equation*}
que tende a zero quando $k \searrow  0$ pela continuidade de $z(t,x)$.
Agora calculamos o último termo. Definindo $\mathcal{J}_k(x) = 2k\sqrt{\pi}\mathcal{G}_k(x)$, temos que $\mathcal{J}$, $k>0$ é uma nova aproximação da identidade. O último termo é,
\begin{equation*}
\frac{1}{4}k^{-1}\pi^{-1/2}\int\mathcal{J}_k(x-y)e^{2\left(h(y) - h(x) \right)}dy \, \left( \int \mathcal{G}_k(x-y)e^{h(y)-h(x)}\,dy\right)^{-2}.
\end{equation*}

Como estamos no equilíbrio $h(y) - h(x)$ são incrementos Brownianos e podemos fazer um cálculo de variação quadrática para obter,
\begin{equation*}
\frac{1}{2}z_k^{-2}\langle\mathcal{G}_k^2,z^2\rangle \sim \frac{1}{4}k^{-1}\pi^{-1/2}.
\end{equation*}
Então
\begin{equation}\label{KPZ_Mol}
\partial_t h_k = \frac{1}{2}\left[\left(\partial_x h_k\right)^2 + \frac{1}{2}k^{-1}\pi^{-1/2}  \right] + \frac{1}{2} \partial^2_x h_k + \xi + o(1).
\end{equation}
Assim obtemos nossa primeira forma precisa de (\ref{KPZ_R}).

\item[(3)] Suponha que, no lugar  de amolecer a solução de (\ref{HC_KPZ}), amolecemos o ruído branco no espaço, utilizando $\mathcal{G}_k(x)$,

\begin{equation*}
\xi_k(t,x) = \int \mathcal{G}_k(x-y)\xi(y)\,dy.
\end{equation*}

Como tal operação é linear, e sendo $\xi_k(t,x)$ Gaussiana com média zero e covariância,
\begin{align*}
\mathbb{E}\left[\xi_k(t,x),\xi_k(s,y) \right] = C_k(x-y)\delta(t-s),
\end{align*}
onde,
\begin{align*}
C_k(x-y) = \int \mathcal{G}_k(x-u)\mathcal{G}_k(y-u)\,du,
\end{align*}
assim, em particular temos que,
\begin{align*}
C_k(0) = \frac{1}{2}k^{-1}\pi^{-1/2}.
\end{align*}

Seja $z_k(t,x)$ a solução da equação estocástica do calor com ruído amolecido,
\begin{align*}
\partial_t z_k = \partial_x^2 z_k - z_k \xi_k, && t>0, && x \in \mathbb{R}.
\end{align*}
Não é difícil provar que $z_k \rightarrow z$ uniformemente em conjuntos compactos, e como $z(t,x)>0$ para $t>0$ podemos definir, 
\begin{align*}
h_k(t,x) = \log z_k(t,x),
\end{align*} 
e $h_k(t,x)$ converge para $z(t,x) = log h(t,x)$. Pela fórmula de Itô obtemos,
\begin{equation}\label{RB_Mol}
\partial_t h_k = \frac{1}{2}\left[\left(\partial_x h_k\right)^2 - C_k(0)  \right] + \frac{1}{2} \partial^2_x h_k + \xi_k.
\end{equation}

Comparando (\ref{KPZ_Mol}) com (\ref{RB_Mol}) teria sentido pensar que no limite 
$k \longrightarrow 0$ as soluções poderiam coincidir.
   
\item[(4)] A solução de Hopf-Cole é obtida aproximando a equação KPZ pela energia livre de polímeros aleatórios dirigidos, e pela função de alturas de exclusão assimétrica \cite[veja Section 3.12]{BG97}. Este limite fracamente assimétrico é esperado para uma ampla classe de sistemas com uma simetria ajustável.
  
\item[(5)] A solução de Hopf-Cole tem os expoentes que foram preditos em \cite{BQS11}. Podemos também obter algumas das flutuações preditas \cite{ACQ11}, \cite{SS10a}, \cite{SH10}.
\end{enumerate}  

Diante do exposto, as evidências de que a transformada de Hopf-Cole da solução da equação estocástica do calor é a solução da equação KPZ são esmagadoras. A dificuldade está em encontrar uma definição apropriada para  (\ref{KPZ}) para fazer com que as soluções coincidam. Este problema foi resolvido por Martin Hairer; para os leitores interessados em mais detalhes sobre esse assunto, indica-se a referência \cite{Hai11}.

\chapter{Revisão de alguns esquemas numéricos para SPDE semilineares}\label{cap33}

   Como veremos, muitos esquemas numéricos usados para aproximar soluções de equações em derivadas parciais podem ser adaptados para resolver numericamente equações estocásticas em derivadas parciais. São os casos dos métodos de diferenças finitas, elementos finitos e métodos espectrais. Neste capítulo ficaremos concentrados na solução de equações do tipo semilinear (\ref{semiL}). A chave para aproximar a solução de este tipo de equações vai ser a aproximação do ruído a partir de um processo $Q$-Wiener (\ref{QWdef}) (Em particular vamos trabalhar com processos cilíndricos cujas trajetórias vão ser pontes Brownianos (\ref{PB})). Alguns destes métodos foram desenvolvidos em anos recentes. Resultados destacados nesta área foram obtidos por cientistas como Arnulf Jentzen, Peter Klodem, Gabriel J. Lord, Catherine E. Powel, Tony Shardlow entre outros em trabalhos como \cite{JER12,JEK11,JEN11,JEK08,LOR04,LPS14}.\\
\\
\section{Derivação dos métodos}
De aqui em diante vamos supor que sejam válidas todas as hipóteses para que as definições seguintes façam sentido:

\begin{spn}\label{A1}
(\textbf{Operador Linear A}) Seja $\lbrace \lambda_i \rbrace_{i\in \mathbb{N}}$ um conjunto de números reais positivos e seja $\lbrace e_j \rbrace$ uma base ortonormal de $H$. O operador $A$ é dado por:
\begin{equation}
A\cdot X = \sum_{j=1}^{\infty} -\lambda_j \langle X,e_j\rangle_H \cdot e_j,
\end{equation} 
para cada $X \in D(A)$ onde $D(A) = \lbrace X \in H : \sum_{j=1}^{\infty} \vert\lambda_j \vert^2 \langle X,e_j \rangle^2 < \infty \rbrace \subset H$. O operador linear $A$ é fechado em $H$ e define um semigrupo analítico $e^{tA}:H\longrightarrow H$ para $t > 0$.
\end{spn}

\begin{spn}\label{A2}
(\textbf{Desvio F}) O operador $F:H \longrightarrow H$ é um operador infinitamente Frechet-diferenciável e 
$inf_{v\in H}\lbrace \parallel F^{(n)}(v)\parallel_{L(H,H)} \rbrace < \infty$ para cada $n\in\mathbb{N}$.
\end{spn}

\begin{spn}\label{A3}
(\textbf{Difusão G}) O operador $G: H \longrightarrow L(U,D)$ é infinitamente Frechet-diferenciável e os operadores $ e^{tA} G^{(n)}(v)(\omega_1,\cdots,\omega_n)$ e $(-A)^\gamma e^{tA} G(v)$ são da classe Hilbert-Schmidt (\ref{HSc}). Além disso existe uma família de números reais $\lbrace (L_n)\rbrace_{n \in \mathbb{N}}$ e números reais $\theta$, $\rho \in (0,\frac{1}{2}]$ e $\gamma \in (0,1)$ tais que:
\begin{align*}
\parallel e^{tA}G^{(n)}(v)(\omega_1,\cdots,\omega_n) \parallel_{HS(H,U)} &\leq L_n \cdot (1+\parallel v \parallel_H)\cdot \parallel\omega_1\parallel\cdots\parallel\omega_n\parallel_H\cdot t^{\theta-\frac{1}{2}}\\
\parallel e^{tA}(G(v)-G(w))\parallel_{HS(H,U)} &\leq L_0\cdot\parallel v-w\parallel_H t^{\rho - \frac{1}{2}}\\
\parallel(-A)^\gamma e^{tA} G(v)\parallel_{HS(H,U)} &\leq L_0\cdot (1+\parallel v\parallel_H)\cdot t^{\rho - \frac{1}{2}},
\end{align*}
para cada $v,\omega,\omega_1,\cdots\omega_n \in H$, $n = \lbrace 1,2,\cdots\rbrace$ a cada $t \in (0,T]$.
\end{spn}

\begin{spn}\label{A4}
\textbf{(Valor inicial)} Seja $X_0 :\Omega \longrightarrow D(-A)$ um mapeamento $\mathcal{F}_0/\mathcal{B} (D((-A)^\gamma))$ mensurável com $\mathbb{E}\parallel (-A)^\gamma X_0\parallel_H^p \leq \infty$ para algum $p \in [1,\infty)$.
\end{spn}
\begin{pps}
Se as suposições (\ref{A1}),(\ref{A2}),(\ref{A3}) e (\ref{A4}) são satisfeitas então existe um único processo estocástico previsível $X:[0,T]\times \Omega \longrightarrow D((-A)^\gamma)$ com $\mathbb{E}\parallel (-A)^\gamma X_t\parallel_H^p < \infty$ para $p\in [0,\infty)$ tal que:
\begin{equation}
\mathbb{P}\left[ X_t = e^{tA} x_0 + \int_0^t e^{(t-s)A} F(X_s)ds + \int_0^t e^{(t-s)A} G(X_s))dW_s \right] = 1.
\end{equation}
\end{pps}
A demonstração pode ser encontrada em \cite{JEN11}.\\

Vamos enfatizar que estamos trabalhando com um problema de valor inicial $X_0(x) \in H$ com condições de contorno periódicas $X_t(0)\equiv X_t(1)$ onde $t \in [0,T]$ e $x \in [0,1]$, para a equação semilinear:
\begin{equation}
dX_t = \left[AX_t + F(X_t) \right]dt + G(X_t)dW_t, 
\end{equation}
onde $H$ é um espaço de Hilbert que em nosso caso tomaremos como sendo $L^2([0,1])$, os operadores $A$, $F$ e $G$ satisfazem as suposições (\ref{A1}), (\ref{A2}), (\ref{A3}) e $W_t$ é um processo $Q$-Weiner cilíndrico. As aproximações são entendidas no sentido da seguinte norma:

\begin{equation}
\Vert X_t - Y_t\Vert = \left(\mathbb{E}\left[\int_0^1 \vert X_t(x) - Y_t(x)\vert^2 dx\right]\right)^{1/2}.
\end{equation}\label{erro}

É claro que a equação estocástica do calor com ruído branco multiplicativo é um caso particular da equação (\ref{semi1}) quando fazemos $A=\nu\partial_{x}^2$,\, $F \equiv 0$ e $G(X,t) = \lambda X_t$:
\begin{equation}
dX_t = \nu A X_t dt + \lambda X_t dW_t, \label{ESC}
\end{equation}

Vamos construir o processo $Q$-Weiner $W_t$ como em (\ref{QWdef}). Um processo $Q$-Weiner cilíndrico $\left( Q \equiv I \right)$ é também chamado de ruído branco já que $\dot{W}_t$ possui as propriedades deste último. Esta vai ser de fato a caracterização de ruído branco que vamos considerar ao longo deste trabalho. Então, temos que $q_j = 1$ para cada $\forall j\in \mathbb{N}$ o processo $Q$- Weiner vai ter a forma:
 \begin{equation}\label{QW}
 W_t(x) = \sum_{j=1}^{\infty} \beta_j(t) \chi_j(x), 
 \end{equation}
onde $\beta_j(t)$ são movimentos Brownianos escalares. As funções  $\chi_j$ vão ser as autofunções do operador $A$. Para discretizar $W_t$ truncamos a série (\ref{QW}) no valor $N$ que vai coincidir com o número de funções a serem consideradas na expansão da aproximação da solução e com o número de elementos na partição no espaço (dependendo do tipo de solução que estejamos levando em consideração, forte, fraca ou \textit{mild}).
 
\subsection{Esquemas baseados na solução forte}

O método de diferenças finitas pode ser facilmente adaptado a SPDE. Para ilustrar isto consideremos a equação semiliniar (\ref{semiL}) com condições de contorno de Dirichlet homogéneas. Onde estamos considerando parâmetros $\varepsilon$, $\sigma> 0$, um termo da reação $F: R \longrightarrow R$, e $W_t$ é um processo $Q$-Wiener em $L^2(0,1)$. Introduzimos os pontos da malha $x_j = jh$ para $h = a / J$ e $j = 0,. . . , J$. Seja $Y_t^J$ a aproximação em diferenças finitas de $[X_t(x_1),. . . , X_t(x_{J-1})]^T$, resultante da aproximação em diferenças centradas $A^D$  para o Laplaciano. Isto é $Y_t^J$ e a solução de,
\begin{equation}
dY_t^J = \left[ A^D Y_t^J + F(Y_t^J) \right]dt + G(Y_t^J)dW^J_t,
\end{equation}
com $Y^J_0 = [X_0(x_1),\cdots,X_0(x_{J+1})]^T$ e  $W^J_t = [W_t(x_1),\cdots,W_t(x_{J+1})]^T$. O método das diferenças finitas é adaptável a outras condições de contorno, por exemplo periódicas, alterando $A^D$ e $W^J$. Agora podemos escolher uma discretização no tempo para obter o sistema de equações algébrico associado. Assim por exemplo, derivamos o esquema de Euler Maruyama semi-implícito,
\begin{equation}
Y^J_{n+1} = (I+\Delta t A^D)^{-1}\left[Y^J_n + F(Y^J_n)\Delta t + G(Y^J_n)\Delta W^J_n\right],
\end{equation}
com $Y^J_0 = X^J_0$ e $\Delta W^J_n = W^J_{t_{n+1}} - W^J_{{t_n}}$.

\subsection{Esquemas baseados na solução fraca}

Para a aproximação de Galerkin da equação estocástica de evolução (\ref{semiL}) devemos introduzir um espaço de dimensão finita $\tilde{V}$ e consideramos $\tilde{P}:H \longrightarrow \tilde{V}$ a projeção ortogonal sobre $\tilde{V}$. Então estamos procurando um processo $\tilde{X}_t \in \tilde{V}$ definido por,
\begin{equation}
\begin{aligned}
\langle \tilde{X}_t,v\rangle = \langle \tilde{X}_0,v\rangle + \int_0^t \left[ \langle A\tilde{X}_s,v \rangle + \langle f(\tilde{X}_s),v\rangle \right]ds \\
+ \langle\int_0^t  G(\tilde{X}_s\cdot dW_s)ds,v\rangle, && t\in [0,T], && v\in \tilde{V},
\end{aligned}\label{Gal_Fraca}
\end{equation}
com dado inicial $\tilde{X}_0 = \tilde{P}X_0$. Ou seja, estamos considerando a  solução fraca (\ref{fraca}) da equação (\ref{semiL}). A aproximação de Galerkin $\tilde{X}_t$ também satisfaz,
\begin{equation}
d\tilde{X}_t = \left[ \tilde{A}\tilde{X}_t + \tilde{P} F(\tilde{X}_t) \right]dt + \tilde{P} G(\tilde{X}_t)dW(t),\label{Gal_Int_Form}
\end{equation}
onde $\tilde{A} = \tilde{P} A$. Para discretizar no tempo podemos aproximar $d\tilde{X}_{t_n}$ por $X_{n+1} - X_n$ para $t_n = n\Delta t$ onde estamos chamando de $X_n$ a $X_{t_n}$. Dado que $\tilde{V}$ é de dimensão finita temos que a equação (\ref{Gal_Fraca}) é uma SODE e podemos aplicar o método de Euler semi-implícito com passo de tempo $\Delta t > 0 $ para obter o algoritmo iterativo, 
\begin{equation}
\tilde{X}_{n+1} = \left(I + \Delta t \tilde{A} \right)^{-1}\left[\tilde{X}_n + \tilde{P} F(\tilde{X}_n)\Delta t + \tilde{P} G(\tilde{X}_n) \Delta W_n \right],
\end{equation}
com $\Delta W(s) = \int_{t_n}^{t_{n+1}} dW(s)$. Na prática é necessário aproximar $G$ por algum $\mathcal{G}: \mathbb{R}_+ \times H \longrightarrow L^2_0$ e estudamos a aproximação definida por, 
\begin{equation}
\tilde{X}_{n+1} = \left(I + \Delta t \tilde{A} \right)^{-1}\left( \tilde{X}_n + \tilde{P} F(\tilde{X}_n)\Delta t + \tilde{P} \int_{t_n}^{t_{n+1}}\mathcal{G}(s,\tilde{X}_n) dW(s) \right),\label{GalIF}
\end{equation}
para dado inicial $\tilde{X}_0 = \tilde{P}X_0$. Dada uma base ortonormal $\lbrace \chi_1,\cdots \chi_J \rbrace$ uma forma usual de tomar dita aproximação é $\mathcal{G}(s,u)= G(u) P_J$ onde $P_j$ e a projeção ortogonal sobre $V_J = span\lbrace \chi_1,\cdots \chi_J \rbrace$.

\begin{center}
\textbf{Euler Galerkin semi-implícito}
\end{center}

\begin{figure}[h!]\label{fig:EECE}
\centering
\includegraphics[scale=.25,trim = 5cm 0 0 0]{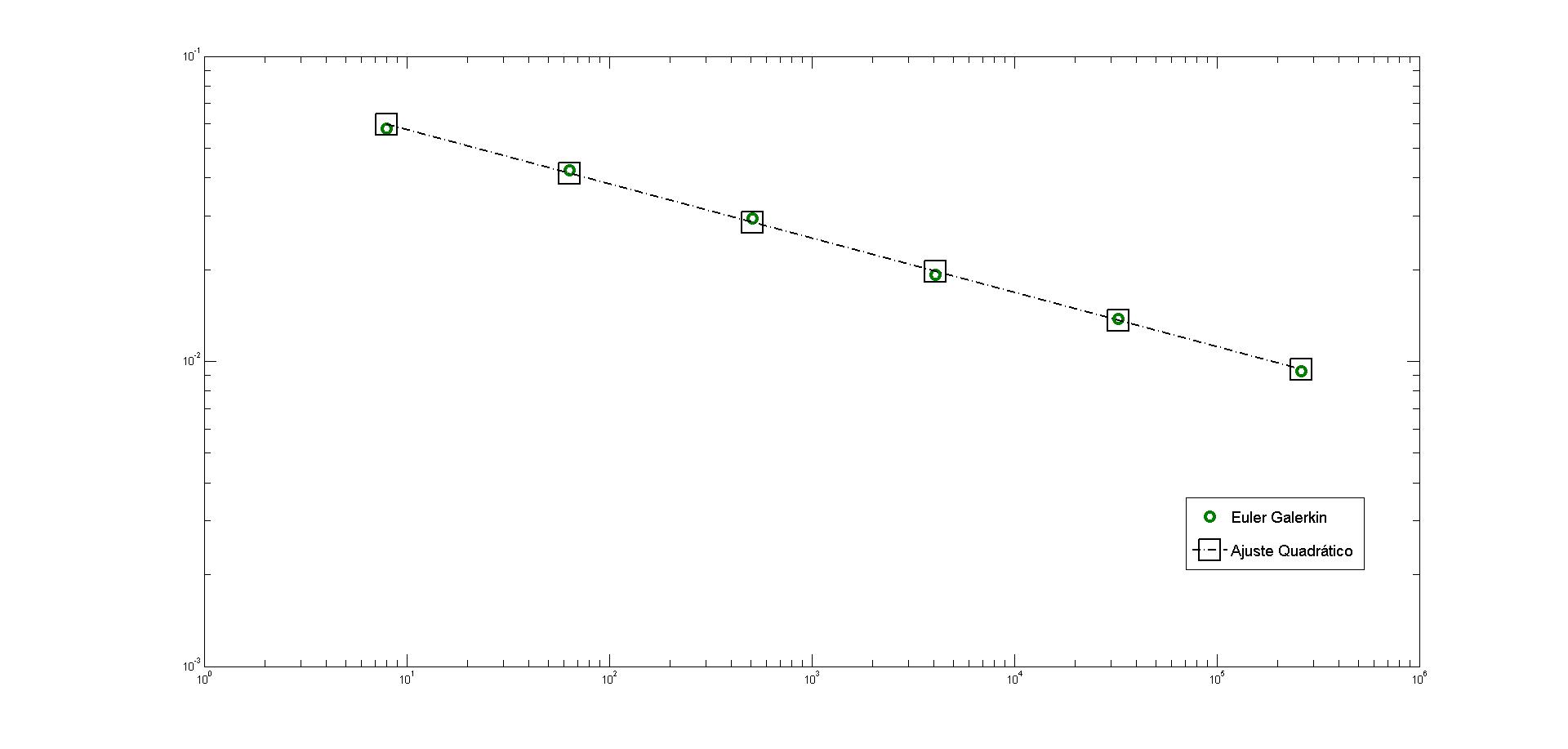} 
\caption{Erros do método Euler Galerkin semi-implícito (\ref{EGSI}) sobre $250$ realizações do ruído branco aplicado na equação estocástica do calor com ruído branco multiplicativo (\ref{ESC}) com parâmetros $\nu = 1$ e $\lambda = 1$ para $N = 2, 4, 8, 16, 32, 64$.}
\end{figure}

Suponhamos que o operador $A$ satisfaz a suposição (\ref{A1}) e denotemos por  $\varphi_j$ as autofunções de $A$ com autovalores $\lambda_j$ para $j\in\mathbb{N}$. Para a aproximação espectral de Galerkin nos escolhemos $\tilde{V} = V_J := span\lbrace \varphi_1\cdots\varphi_J \rbrace$ e escrevemos $X^J$ para a aproximação do método das linhas $\tilde{X}$, $P_J:H\longrightarrow V_J$ é a projeção ortogonal $\tilde{P}$ e $A_J = P_J A$ vai ser  $\tilde{A}$. Então de (\ref{GalIF}) segue:
\begin{align}
d X^J_t = \left[ A_J X^J_t + P_J F(X^J_t) \right]dt + P_J G(X^J_t)dW(t) && X^J_0 = P_J X_0,
\end{align}
assim, chamando de $Y^J_n$ à aproximação de $[X_n(x_1),\cdots,X_n(x_{J-1})]$, obtemos o método de Euler Galerkin semi-implícito,
\begin{equation}
Y^J_{n+1} = \left(I + \Delta t A_J \right)^{-1}\left( Y^J_{n} + F_n(Y^J_{n})\Delta t +  \int_{t_n}^{t_{n+1}}\left(\mathcal{G}(s,Y^J_{n}) dW_{J}\right)_n \right),
\end{equation}\label{EGSI}
onde $F_n(Y^J_{n})$ e $\left(\mathcal{G}(s,Y^J_{n})
dW_{J}\right)_n$ representam os $n$-ésimos coeficientes de 
Fourier de ditas funções. Aqui temos escolhido como aproximação 
de $G$ o operador $\mathcal{G}(s,u) = G(u)P_J$ e o dado inicial
$X^J_0 = P_J X_0$. Esta aproximação é particularmente útil para
problemas com ruído aditivo onde $U=H$ e as autofunções de $Q$
coincidem com as autofunções de $A$. Nesse caso 
$P_J \int_{t_n}^{t_{n+1}}\mathcal{G}(s,X^J_n) dW^J_n$ 
pode ser calculada de forma explícita aproveitando a ortogonalidade da base. A Figura \ref{fig:EECE} mostram o desempenho do método quando é aplicado na equação estocástica do calor com ruído branco multiplicativo (\ref{ESC}) com parâmetros $\nu = 1$ e $\lambda = 1$. Os experimentos foram realizados sobre $250$ realizações do ruído branco e o erro é tomado no sentido de (\ref{erro}).\\

\newpage
\subsection{Esquemas baseados na solução \textit{mild}}

Até agora temos derivado esquemas baseados na solução forte e na solução fraca da equação (\ref{semiL}). Agora veremos esquemas numéricos que derivam da solução \textit{mild} (\ref{leve}). Para tal fim consideremos um processo $Q$-Weiner como em (\ref{QW}) e suponhamos satisfeitas as suposições (\ref{A1}),(\ref{A2}),(\ref{A3}) e (\ref{A4}). Suponhamos também que as autofunções do processo $Q$-Wiener coincidem com as autofunções do operador $A$ as quais vamos a denotar de  $\lbrace \varphi_j \rbrace_{j\in\mathbb{N}}$.\\

\begin{center}
\textbf{Lord Rougemont modificado}
\end{center}

 Vamos desenvolver uma variante do método de Lord-Rougemont introduzido em \cite{LOR04} no ano 2003. Para este fim vamos considerar primeiro uma discretização no passo de tempo $t_n = n\cdot\Delta t$ e escrever nossa aproximação como,
 \begin{align}
X_{n+1}(x) \approx e^{\Delta t A} X_n(x) + \int_{t_n}^{t_{n+1}} e^{(t_{n+1}-s)A} F(X_s(x))ds + \int_{t_n}^{t_{n+1}} e^{(t_{n+1}-s)A} G(X_s(x))dW_s(x), \label{LR1}
\end{align}
onde $X_n(x) = X_{t_n}(x)$ e as integrais são integrais de Itô. Definamos  $V_J =span\lbrace \varphi_1,\cdots,\varphi_J \rbrace$. Seja uma aproximação do operador $G$ dada por $\mathcal{G}(s,u) = G(u)P_N$ onde $P_J$ é a projeção sobre $V_J$ e vamos aproximar as integrais em (\ref{LR1}) aproximando o integrando pelo limite a esquerda, ou seja, $ e^{(t_{n+1}-s)A} F(X_s(x)) \approx  e^{\Delta t A}F(X_n(x))$ e $ e^{(t_{n+1}-s)A} G(X_s(x)) \approx  e^{\Delta t A} G(X_n(x))$, e ainda que $dW^J_s(x):=P_J dW_s(x) = \sum_{j=1}^J \sqrt{q_j} \cdot d\beta_j(s) \cdot \varphi_j(x)$. Assim obtemos,

\begin{equation}
\begin{aligned}
X_{n+1}(x) &\approx e^{\Delta t \cdot A} \left[ X_n(x) + F(X_n(x))\int_{t_n}^{t_{n+1}} ds + G(X_n(x))\int_{t_n}^{t_{n+1}} dW^J_s(x) \right]\\
&\approx e^{\Delta t \cdot A} \left[ X_n(x) + F(X_n(x))\Delta t + G(X_n(x))\Delta W^J(x) \right].
\end{aligned}\label{LR2}
\end{equation}

Levando em conta a representação do processo $Q$-Weiner (\ref{QW}) podemos reescrever,
\begin{equation}
\Delta W^J(x) = \sum_{j=1}^J \sqrt{q_j} \xi_j \varphi_j(x),
\end{equation}
onde cada $\xi_j$ é uma variável aleatória independente com distribuição gaussiana com meia $0$ e variância $\Delta t$ \textit{i.e.} $\xi_j \thicksim N(0,\Delta t)$. Se agora consideramos a projeção de (\ref{LR2}) sobre o espaço $V_J$ e chamamos de $Y^J_n$ os coeficientes de Fourier da projeção $P_J X_n$, poderíamos escrever o esquema numérico como sendo,

 \begin{align}
Y^J_{n+1} &\approx e^{\Delta t \cdot A_J} \left[ Y_n^J + F_J(Y^J_n)\Delta t + \left( G(Y^J_n)\Delta W^J\right)_J \right], \label{LRM}
\end{align}

onde $A_J = diag\lbrace \lambda_1,\cdots,\lambda_J \rbrace$ e $F_J(Y^J_n)$ e $\left( G(Y^J_n)\Delta W^J\right)_J$ representam os vetores dos coeficientes de Fourier das  funções $F(X_n(x))$ e $G(X_n(x))\Delta W^J(x)$, respetivamente. Os coeficientes de Fourier de $G(X_n(x))\Delta W^J(x)$ não são difíceis de calcular uma vez que os coeficientes de Fourier de $\Delta W^J$ são conhecidos. Nos apoiando na transformada de Fourier podemos fazer o cálculo de forma relativamente eficiente.

\begin{center}
\textbf{Esquema Milstein para SPDE}
\end{center}

 Davie $\&$ Gaines  demonstraram em \cite{DAG00}  que qualquer esquema numérico aplicado a uma SPDE semilinear (\ref{semiL}) com $F \equiv 0$, que utiliza apenas valores equidistantes do ruído, não pode convergir mais rapidamente do que a taxa $1/6$ com respeito ao esforço computacional. Muller-Gronbach e Ritter mostraran que este é também um limite inferior para a taxa de convergência no caso do ruído multiplicativo, (veja \cite{MUG07}). Eles também provaram que a taxa de convergência geral $1/6$ não pode ser melhorada tomando passos de tempo não uniformes, \cite{MUR07,MRW08}.\\

  Taxas de ordem superior a $1/6$ para SPDEs não-lineares da forma (\ref{semiL}) foram obtidos para  tipos  de ruído mais suaves. Por exemplo, em \cite{Hau02} aplicado o regime de Euler-linear implícito e explícito e o esquema de Crank-Nicholson ao SPDE (\ref{semiL}) com um processo ruído de dimensão infinita. Para o ruído de classe de traço, ele obteve a ordem $1/4$ com respeito ao esforço computacional, mas no caso do ruído branco espaço-tempo  a taxa de convergência não foi melhor do que a taxa barreira de Davie-Gaines $1/6$. Da mesma forma, em \cite{LOR04}, os autores propuseram um esquema numérico que mostrou ser útil quando o ruído é muito suave no espaço, em particular com  regularidade de tipo Gevrey (\cite{LOR04}). No caso do ruído branco espaço-tempo a convergência também não  supera a taxa  barreira de Davie-Gaines  $1/6$.\\

O seguinte algoritmo é um exemplo de que a barreira de Davie-Gaines pode ser superada chegando numa taxa de convergência de $1/4$ com respeito ao esforço computacional. A demonstração pode ser encontrada em \cite{LPS14}. Um estudo mas geral destes esquemas baseados em expansões em séries de Taylor dos operadores envolvidos na equação (\ref{semiL}) foi desenvolvido por Arnulf Jentzen e é apresentado em \cite{JEK11}.\\

Vamos considerar novamente a formulação \textit{mild} (\ref{leve}). Se supomos que $G(u)$ é pelo menos uma vez Frechet diferenciável poderíamos considerar a aproximação $G(X_s) \approx G(X_0) + G'(X_0)(X_s-X_0)$ e ao substituir na formulação \textit{mild} obtemos,
\begin{align*}
X_t &\approx e^{tA} X_0 + \int_0^t e^{(t-s)A} F(X_s)ds + \int_0^t e^{(t-s)A} G(X_0)dW_s + \int_0^t e^{(t-s)A} G'(X_0)(X_s-X_0)dW_s\\
&\approx e^{tA} \left( X_0 + t\cdot F(X_s) + \int_0^t G(X_0)dW_s + \int_0^t G'(X_0)(X_s-X_0)dW_s \right),\numberthis, \label{app1}
\end{align*}
para $t\in[0,T]$. A aproximação anterior ainda não se adequa para a implementação de um esquema de aproximação numérico devido à presença do termo $X_s$ na segunda integral. Agora, nos apoiando na formulação fraca podemos considerar a seguinte estimativa $X_s \approx X_0 + \int_0^s G(X_0)dW_s$. Substituindo na aproximação (\ref{app1}) obtemos,

\begin{equation}
X_t \approx e^{tA} \left( X_0 + t\cdot F(X_s) + \int_0^t G(X_0)dW_s + \int_0^t G'(X_0)\left( \int_0^s G(X_0)dW_u\right) dW_s \right).
\end{equation}

Combinar um esquema de discretização temporal de tipo Euler com uma discretização no espaço de tipo Galerkin espectral resulta no seguinte esquema numérico,
\begin{equation}
\begin{aligned}
Y^J_{n+1} = e^{\Delta t A_J} \left( Y_n + F(Y^J_n)\Delta t + G(Y^J_n)\left(W^J_{n+1} -W^J_n \right) +  \right.\\
\left. + \int_0^t G'(Y^J_n) \left( \int_0^s G(Y^J_n) dW^J_u \right) dW^J_s  \right),
\end{aligned} \label{Milstein}
\end{equation}
onde  $Y^J_n:\Omega \longrightarrow H$ são mapeamentos $\mathcal{F}/\mathcal{B}(H)$-mensuráveis e  $Y^J_0 = P_J X_0$. A última integral na equação (\ref{Milstein}) sugere um inconveniente computacional para a utilização deste método devido ao cálculo das integrais iteradas. Este problema é superado com a ajuda da formula, 
\begin{equation}
 \int_0^t G'(Y^J_n) \left( \int_0^s G(Y^J_n) dW^J_u \right) dW^J_s = G'(Y^J_n)G(Y^J_n)\left( \left( W^J_{n+1} - W^J_n \right)^2 + \Delta t \sum_{j=1}^J  \lambda_j \cdot (\varphi_j)^2 \right),\label{form}
\end{equation}
(veja \cite{JER12}) onde o somatório em (\ref{form}) e calculado uma vez somente no algoritmo.

\section{Algoritmos e códigos}

Foram programados em Matlab códigos para aproximar as soluções de equações semilineares utilizando os esquemas, \textit{Euler Galerkin semi-implícito (\ref{EGSI}), Lord Rougemont (\ref{LRM}) e Milstein (\ref{Milstein})}. Os códigos foram testados para a equação estocástica do calor (\ref{ESC}) que é o caso de equação semilinear de interesse neste trabalho já que a transformada de Cole-Hopf de dita solução modela corretamente o processo de crescimento de interfaces por deposição balística.  Na continuação apresentamos alguns dos resultados obtidos que mostram como os códigos programados fornecem evidencias de convergência. \\

Vamos considerar o problema de valor inicial e de contorno da secção anterior para a equação semilinear com $F \equiv 0 $, $G$ é um operador de tipo Nemyiski, i.e. $G(X_t(x)) = (g \circ X)(t,x)$ onde $g:\mathbb{R}^2 \longrightarrow \mathbb{R}^2$ e $A$  é o operador Laplaciano $\frac{\partial^2}{\partial x^2}$. Assim, a equação estocástica do calor escrita como equação integral abreviada fica:

\begin{equation}
dX_t = \nu\frac{\partial^2 X_t}{\partial x^2} dt + \lambda X_t \cdot dW_t \label{semi1}
\end{equation}

\begin{figure}[h!]
\centering
\includegraphics[scale=.4,trim = 4cm 0 0 0]{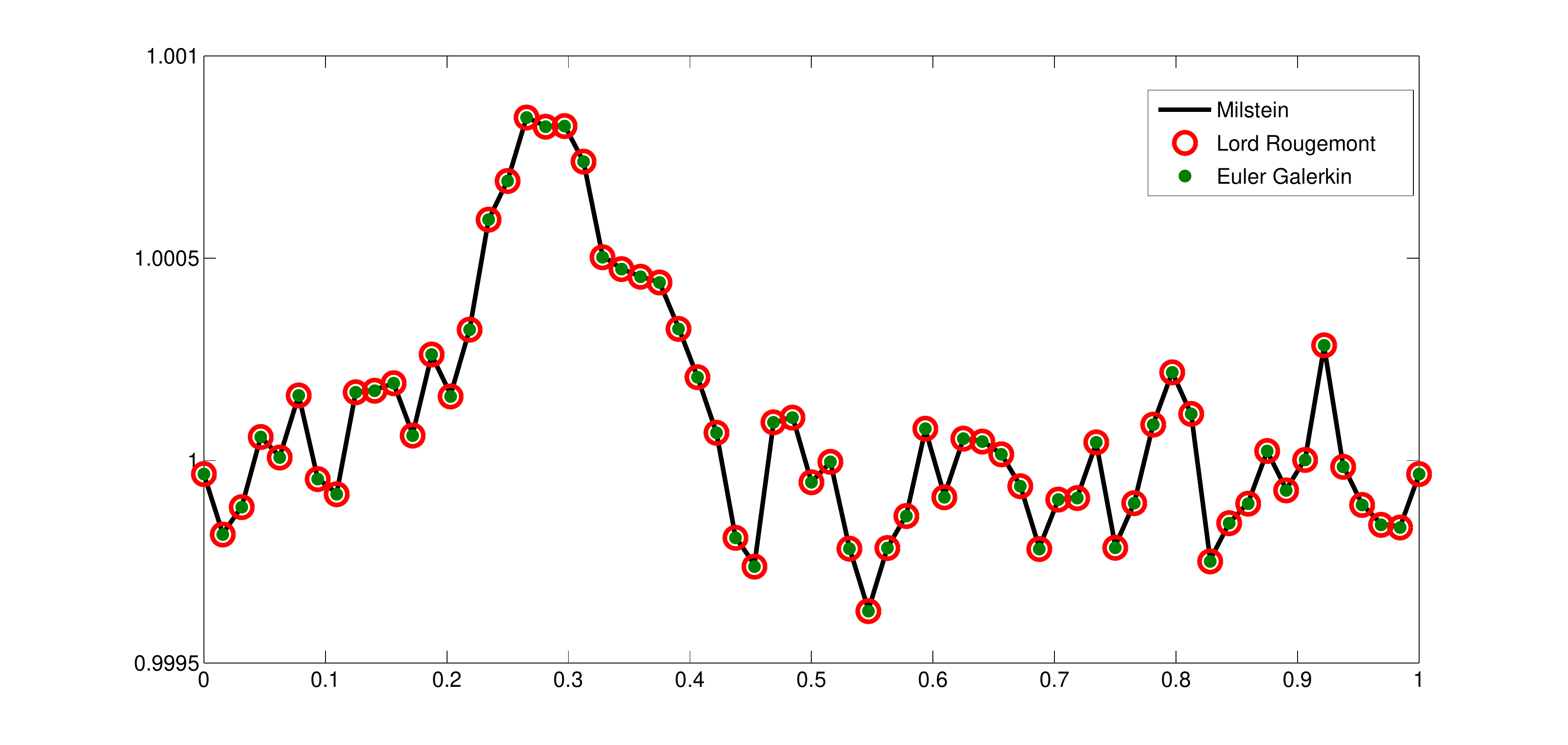} 
\caption{Trajetórias obtidas aplicando os métodos \textit{Euler Galerkin semi-implícito},\textit{ Lord Rougemont} e \textit{Milstein} na equação estocástica do calor até o tempo $T = 1$ para uma mesma realização do ruído branco para valores dos parâmetros $\nu = 1$ e $\lambda = 1$.}\label{aprox}
\end{figure}

Vamos comparar as aproximações calculadas  até o tempo $T=1$ para $N = 128$ (128 funções na base para aproximar o ruído branco e a solução) usando os  algoritmos de \textit{Euler Galerkin semi-implícito}, \textit{Lord Rougemont} e \textit{Milstein}. Vamos partir de uma superfície plana $X_0 \equiv 1$. As trajetórias mostradas foram obtidas para a mesma realização do ruído branco. Os valores dos parâmetros foram $\nu = 1$ e $\lambda = 1$. Os tempos reais de processamento em um computador Core I7 6700HQ foram: \textit{Euler GS} ($13min:07seg$.), Lord-Rougemont ($13min:20seg$.) e \textit{Milstein} ($1min:20seg$.) As trajetórias são mostradas na Figura \ref{aprox}.

Este experimento numérico foi repetido 50 vezes e, em todos os
casos, observamos um boa concordância entre as aproximações
calculadas como mostrado na Figura \ref{aprox}.\\

\begin{figure}[h!]
\centering
\includegraphics[scale=.5,trim = 6cm 0 0 0]{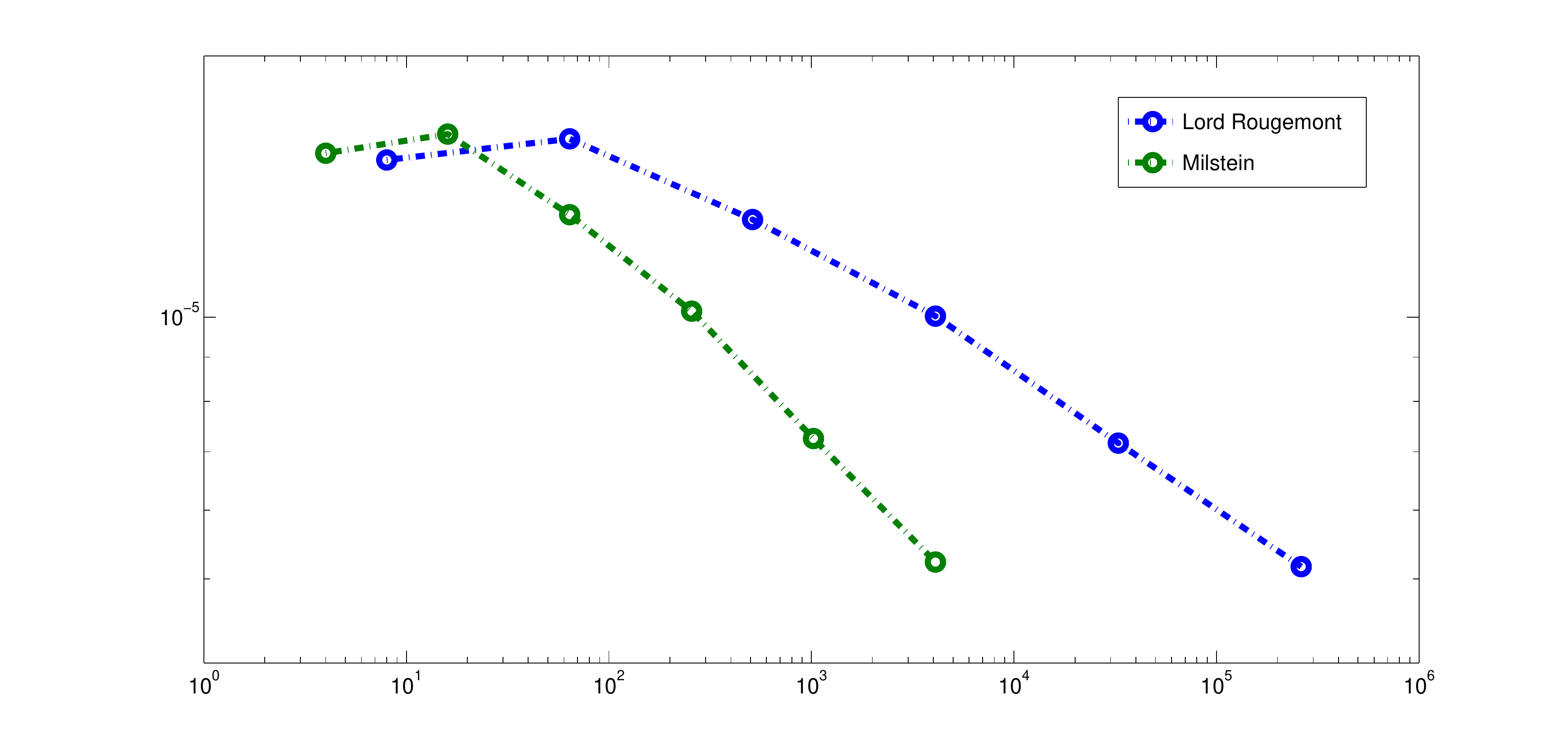} 
\caption{Erro de aproximação no sentido de (\ref{erro}) dos métodos de \textit{ Milstein} e \textit{Lord Rougemont}, para valores de $N = 2,4,8,16,32,64,128$ calculado ao longo de 50 realizações do ruído branco.}\label{ordem}
\end{figure}

Agora queremos saber como é o comportamento do erro cometido por cada algoritmo no sentido da definição (\ref{erro}) com respeito a $N$, ou seja, como varia o erro com respeito ao esforço computacional. Para isso vamos repetir o experimento anterior mas desta vez para várias realizações do ruído branco variando a quantidade de funções na expansão. Usaremos os valores $N = 2,4,8,16,32,64,128$ e estimaremos o erro cometido comparando os resultados da solução $Y_{N+1}$ e a anterior $Y_N$. A Figura (\ref{ordem}) mostra como o erro decai na medida em que aumentamos o número de funções na expansão com o qual temos mais uma evidencia de convergência. Além disso podemos observar como o método de \textit{ Milstein} apresenta uma ordem de aproximação melhor do que o método de \textit{Lord Rougemont}, como era esperado.

\begin{figure}[h!]
\centering
\includegraphics[scale=.4,trim = 4cm 0 0 0]{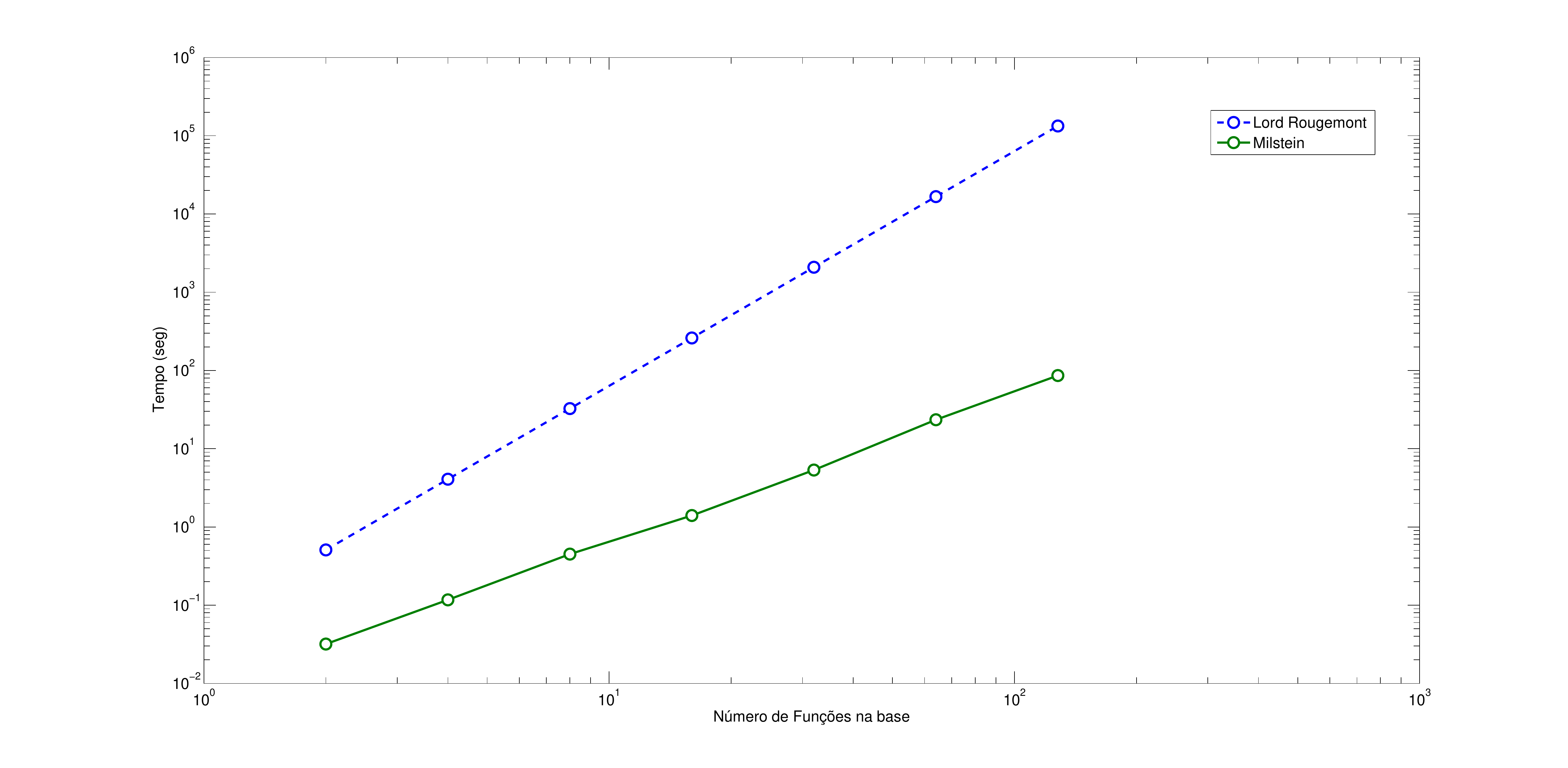} 
\caption{Tempos reais de procesamento dos métodos de \textit{ Milstein} e \textit{Lord Rougemont}, para valores de $N = 2,4,8,16,32,64,128$ calculado ao longo de 50 realizações do ruído branco.}\label{tempos}
\end{figure}

A Figura \ref{tempos} mostra os tempos reais de processamento dos algoritmos de \textit{Milstein} e \textit{Lord-Rougemont} onde podemos observar uma marcada diferencia entre os dois métodos sendo o primeiro deles mais rápido. Os experimentos foram realizados num computador Core I7 6700HQ.\\

\begin{figure}[h!]
\includegraphics[scale=.4]{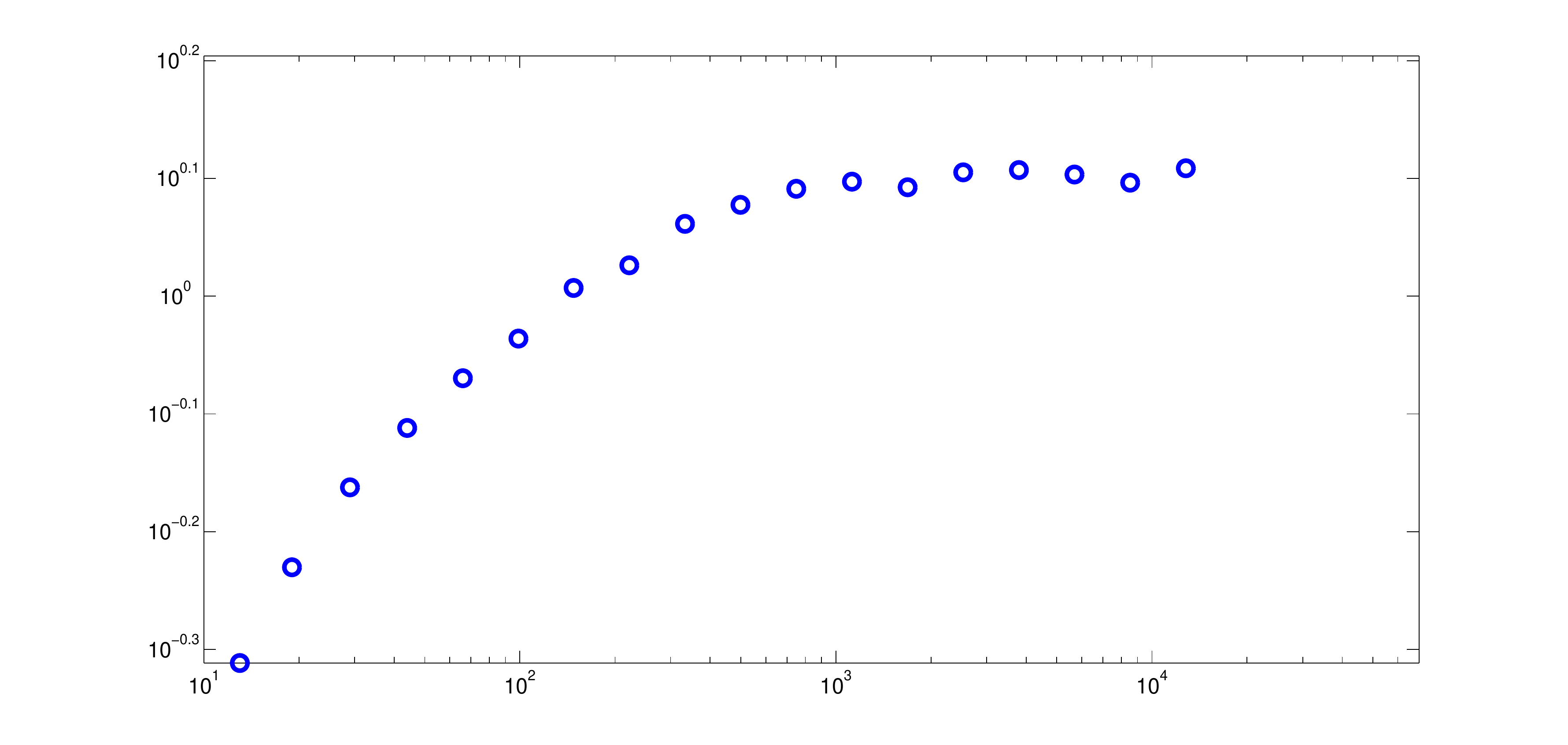} 
\caption{Curva de rugosidade para a transformada de Hopf-Cole da solução da equação estocástica do calor com ruído branco multiplicativo.}\label{Rugo}
\end{figure}

  Pergunta: O que acontece com os coeficientes de crescimento do processo de deposição balística? Se os algoritmos estão trabalhando corretamente então devemos esperar que a transformada de Hopf-Cole da solução obtida corrobore com os valores teóricos previstos para os ditos coeficientes. Para conferir com o caso anterior, fizemos um experimento numérico aplicando o método de \textit{Milstein} (\ref{Milstein}) na equação estocástica do calor. Tomamos $N = 128$ funções na expanção do ruído branco e calculamos a média, sobre 50 realizações, das transformadas de Hopf-Cole das soluções aproximadas obtidas. A Figura (\ref{Rugo}) mostra o gráfico da evolução da rugosidade no tempo.
  
\newpage
Podemos conferir na Figura \ref{Rugo} que de fato temos um tempo de \textit{crossover} que separa os dos regimes transiente (``pré-assintótico'')  e  estacionário (``assintótico'') onde o comportamento da rugosidade varia de um crescimento da ordem de $t^3$ para um comportamento constante.

\chapter{Proposta de modelagem computacional e aproximação numérica para modelos KPZ}\label{cap3}

  Neste capítulo abordaremos a aproximação numérica da solução de uma variante não 
  estocástica da equação KPZ, onde vamos substituir o ruído branco por uma função suave. 
  Especificamente, vamos propor, no nível formal, uma formulação do método de elementos 
  finitos mistos e híbridos com decomposição de domínio para o tratamento numérico desta 
  equação. Tal metodologia
  é baseada no trabalho de  \cite{DPLW93}, onde tal procedimento foi bem sucedido
  para a aproximação de um problema parabólico proveniente da área da dinâmica
  de fluidos em meios porosos. Uma generalização desse método para uma 
  classe mais ampla de problemas parabólicos, assim como rigorosas 
  estimativas de erro e de convergência, pode-se encontrar no trabalho 
  \cite{KPP00}. De fato, aqui estaremos também bem próximos dos trabalhos
  \cite{JV15,EAJV15}, onde uma formulação semelhante foi utilizada e bem sucedida para
  aproximar um modelo da classe pseudo-parabólica em problemas de fluxo
  bifásico em meios porosos. O método de elementos finitos mistos e híbridos 
  teve suas origens nos trabalhos de Fraeijs de Veubeke \cite{Veu65,Veu75} 
  e foi logo detalhadamente analisado por \cite{AB85} (veja-se também 
  \cite{BDF87,BDM85,BDD87}). Provas de convergência para problemas 
  diferenciais em forma forte podem ser encontradas em \cite{Des90,DJR90}, 
  ver também \cite{JBBD97}.

Nesse contexto, é crucial destacar que o método dos elementos finitos é uma 
técnica geral para construir subespaços de dimensão finita de um apropriado 
espaço de Hilbert $V$, a fim de aplicar o método fundamental Ritz-Galerkin 
para um problema variacional \cite{FBMR91}. Esta técnica baseia-se em algumas 
ideias simples. O fundamental é a partição do domínio $\Omega$ (contínuo,
tipicamente em dimensão 1, 2 ou 3) em que o problema diferencial é posto em 
um conjunto (discreto) de ``subdomínios simples'', chamados de {\it elementos}. 
Estes elementos podem ser {\it intervalos}, {\it triângulos}, {\it quadriláteros}, 
{\it tetraedros}, e assim por diante (ver mas detalhes em \cite{FBMR91}). Um espaço $V$ de funções 
definidas em $\Omega$ é então aproximado por tais funções ``simples'' ou 
elementares, definidas em cada subdomínio com condições adequadas de 
compatibilidade, ou condições de transmissibilidades, nas interfaces entre os 
elementos da partição. Essas funções simples são geralmente {\it polin\^omios} 
ou {\it funções obtidas a partir de polin\^omios} por uma conveniente mudança 
de variáveis. O ponto fundamental que queremos destacar aqui 
é o seguinte: um método de aproximação via elementos finitos só pode ser 
considerado em relação a {\it um princípio variacional} e {\it um espaço 
funcional}. Assim, alterando o princípio variacional associado ao modelo 
diferencial e o espaço associado de aproximação no qual o mesmo é posto 
conduzirá também para uma aproximação distinta de elementos finitos, mesmo 
que a solução para o problema original posto no contínuo permaneça inalterada.

Desta forma, fica evidente que é importante a existência de um princípio 
variacional para o modelo diferencial KPZ sob investigação. Nessa direção 
entendemos que é relevante mencionar o trabalho \cite{HSW09}, onde uma 
formulação variacional para a equação Kardar-Parisi-Zhang (KPZ) foi 
apresentada, que por sua vez leva a um potencial termodinâmico, como o
potencial estacionário para o modelo KPZ, e também como para outras equações 
cinéticas relacionados \cite{WERDL11}. Com efeito, no trabalho \cite{HSW09} 
uma formulação variacional teórica da equação KPZ é introduzida, mas 
também dedicado a questões de discretização do ponto de vista teórico, 
ressaltando o valor da consistência em integrações numéricas da KPZ. 
Mesmo reconhecendo as  propriedades analíticas 
da equação KPZ \cite{Qua12,Hai11}, é evidente que investigar o comportamento 
das suas soluções via uma integração numérica direta pode ser uma alternativa 
cientificamente válida para o entendimento e compreensão de modelos KPZ,
e.g., \cite{BS95,BG97,CW78,CV92,Tiago,Fam86,ACQ11,MAM09,KS92,LS98,NuS,MRS86}. 
Essa abordagem foi também utilizada, por exemplo, para obter os expoentes 
críticos para modelos KPZ em uma e mais dimensões espaciais 
\cite{ACQ11,KS92,BQS11,Qua12}. Chamamos a atenção que vários outros 
aspetos são relevantes para uma formulação variacional para KPZ, e 
que estão muito além do escopo do presente trabalho. Por exemplo, a 
formulação variacional discreta leva naturalmente a uma consistente 
discretização da equação KPZ, mas também levando em conta (i) sua 
relação com o funcional de Lyapunov, (ii) a relação flutuação-dissipação, 
e (iii) a invariância de Galileu (simetrias); para maiores detalhes 
sobre tais aspetos, e que vão muito além do escopo do presente trabalho, 
ver, e.g., \cite{HSW09,WERDL11,WDER14}.

\section{Modelagem computacional da KPZ}

\subsection{Formulação mista e formulação fraca}
 
Métodos de elementos finitos em que dois espaços são utilizados para aproximar duas variáveis recebem a denominação geral de métodos mistos. Em alguns casos, a segunda variável é introduzida na formulação do problema devido a seu interesse físico e geralmente está relacionada com alguma derivada da variável original. Por exemplo, vimos no Capitulo precedente a relevância da equação do calor para KPZ. Ou seja, é crucial então uma ``boa'' aproximação do núcleo do calor. Assim, a formulação mista parece ser também adequada para a aproximação do fluxo difusivo. É o caso, por exemplo, das equações de elasticidade, em que o esforço pode ser introduzido para ser aproximado ao mesmo tempo que o deslocamento. As vezes, existem duas variáveis naturais independentes, nestes casos a formulação mista é uma formulação natural. Um exemplo são as equações de Stokes, onde as magnitudes do fen\^omeno físico analisado são a pressão e a velocidade, representadas por duas variáveis independentes.\\

  Embora pareça artificial, na verdade, a formulação mista é muito natural e é usada na matemática há muito tempo. Um exemplo simples é o seguinte: considere a equação linear de segunda ordem como em \cite{Ger54},
\begin{equation} \label{TFEq}
  L(\psi) = k(\sigma)\frac{\partial^2\psi}{\partial\vartheta^2} + \frac{\partial^2\psi}{\partial\sigma^2} = 0,
\end{equation}
onde $\psi$ é uma função de $(\sigma,\vartheta)$. Se introduzimos a variável $\phi$, que é uma função de $(\sigma,\vartheta)$, podemos escrever a equação anterior como,
\begin{align}\label{sisTF}
\frac{\partial\phi}{\partial\vartheta} = - \frac{\partial\psi}{\partial\sigma},
&& \frac{\partial\phi}{\partial\sigma} = k(\sigma)\frac{\partial\psi}{\partial\vartheta}.
\end{align}

Aqui a função potencial $\phi$ passa a ser uma incógnita a ser calculada junto com $\phi$ no sistema (\ref{sisTF}) equivalente a (\ref{TFEq}).\\
 
  A análise matemática e as aplicações de métodos de elementos finitos mistos  têm sido amplamente desenvolvidas desde os anos setenta. Uma análise geral para este tipo de métodos foi desenvolvido pela primeira vez por \cite{Bre74}. Podemos também  mencionar os trabalhos de \cite{Bab73} e de Crouzeix e Raviart \cite{CR73}. Nessos trabalhos foram considerados problemas particulares mas algumas das ideias fundamentais para a análise de métodos mistos foram discutidas. Outros textos importantes nessa área são \cite{FO80} e \cite{For77}, onde são apresentados resultados gerais; ver também \cite{Veu65,Veu75,AB85,BDF87,BDM85,BDD87,Des90,DJR90,JBBD97}.\\

Em uma primeira abordagem de uma possível formulação mista para a equação KPZ consideraremos o problema de valor de contorno unidimensional com condições de contorno de Dirichlet, dado por,
\begin{equation}
\begin{aligned}
h_t &= \frac{\lambda}{2}\left( h_x \right)^2 + \nu h_{xx} + \sqrt{D} \xi  && (t,x) \in J\times \Omega, \\
h(0,x) &= h_0(x), &&  x \in \Omega,  \\ 
h(t,a) &= h_a(t), &&  t \in J, \\
h(t,b) &= h_b(t), &&  t \in J.  
\end{aligned}\label{PVC}
\end{equation}
Em nosso caso $\Omega = [a,b]$ e $J=[0,T]$. Consideremos o espaço $L^2(\Omega)$ (espaço das funções quadrado integrável munido de produto interno e normas induzidas usuais). 
\begin{align*}
\left( f,g \right)_{L^2} &= \int^a_b f(x) \,\, g(x) \mathrm{d}x, \\
\parallel f \parallel_2 &= \sqrt{\left( f,f \right)_{L^2}},
\end{align*}
e consideremos também o espaço $H_1(\Omega)$ das funções de $L^2(\Omega)$ com derivada também  em $L^2(\Omega)$ munido com o mesmo produto interno e mesma norma,
\begin{align*}
 L^2(\Omega) &=  \{\, u : \Omega \rightarrow \mathbb{R} : \| u \|_2 < \infty \,\}, \\
 H_1(\Omega) &= \{\, u\in L^2(\Omega): u_x \in L^2(\Omega) \,\}. 
\end{align*}  

Como já foi dito a função $\xi$ vai ser alguma função suave em $L^2([a,b])$. No caso da equação (\ref{PVC}) identificamos a variável de fluxo $u$,
\begin{equation}\label{varU}
u = -\nu h_x.
\end{equation}

Substituindo (\ref{varU}) no PVC (\ref{PVC}) obtemos a formulação mista do problema,
\begin{equation}
\left\{ \begin{array}{c}
u = -\nu h_x, \\
h_t + u_x - \frac{\lambda}{2\nu^2} u^2 = \xi.
\end{array}
\right.
\end{equation}

Simplificamos a notação,
\begin{align*}
h(0) &= h(0,x),   \\
\alpha &= 1/\nu,   \\
\beta &= -\frac{\alpha^2\lambda}{2}.
\end{align*}

Consideramos a formulação fraca global do problema (\ref{PVC}) onde estamos procurando um mapeamento $\lbrace h,u \rbrace : J \rightarrow L^2 \times H_1$ que satisfaz:
\begin{equation}
\begin{aligned}
\left( \alpha u, v \right)_{L^2} - \left( h,v_x \right)_{L^2} + h(t,b) \cdot v(b) - h(t,a)\cdot v(a)= 0, && v\in L^2,  \\
\left( h_t,w \right)_{L^2} + \beta\left( (u)^2,w \right)_{L^2} + \left( u_x,w \right)_{L^2} = \left( \xi,w \right)_{L^2}, && w\in H_1,  \\
h(0) = h_0.
\end{aligned}\label{FFG}
\end{equation}

 Consideramos agora uma partição homog\^enea do intervalo $[a,b]$ em $m+1$ pontos  
 $a= x_1, x_2, \cdots, x_{m+1} = b$ e com $\Delta x = \frac{b-a}{m}$ e particionamos do domíno  $\Omega$ em $m$ sub-intervalos $I_j = [x_j, x_{j+1}]$. Podemos definir  o método de elementos finitos mistos substituindo os espaços $L^2$ e $H_1$ por espaços convenientes de dimensão finita $V^h$ e $W^h$. Exemplos destes espaços podem ser encontrados na literatura, e.g., \cite{BF91,Che05,BS08}. Nós vamos trabalhar com os espaços  $V^h = Span \{ \varphi_i \mid i=1,\cdots m \}$ e $W^h = Span \{ \psi \mid i=1,\cdots m \}$ onde,
\begin{equation}
 \begin{aligned}
 \varphi_j(x) = \left\{ \begin{array}{cc}
 (x - x_{j-1}) / \Delta x,  & x \in I_{j}, \\
 (x_{i+1}-x) / \Delta x,  & x \in I_{j+1}, \\
 0, & \hbox{ caso contr\'ario},
 \end{array} \right. &&
 \psi_i = \left\{ \begin{array}{cc}
 1, & x \in I_j, \\
 0, & \hbox{ caso contr\'ario}.
 \end{array} \right.
 \end{aligned}\label{BaseEFM}
\end{equation}

 Assim a aproximação de elementos finitos mistos em {\it tempo contínuo} consiste em encontrar $\{\bar{U},\bar{H}\}$, projeções da solução $\{u,h\}$ de (\ref{FFG}) sobre os espaços $V^h$ e $W^h$ que satisfaçam, 
\begin{equation}
\begin{aligned}
\left( \alpha \bar{U}, v \right)_{L^2} - \left( \bar{H},v_x \right)_{L^2} + \bar{H}(t,b) \cdot v(b) - \bar{H}(t,a)\cdot v(a)= 0, && v\in V^h,  \\
\left( \bar{H}_t,w \right)_{L^2} + \beta\left( (\bar{U})^2,w \right)_{L^2} + \left( \bar{U}_x,w \right)_{L^2} = \left( \xi,w \right)_{L^2}, && w\in W^h,  \\
\bar{H}(0) = h_0.
\end{aligned}\label{EFM_Fraca}
\end{equation}

\subsection{Elementos finitos mistos e híbridos}
  
A formulação clássica do método de elementos finitos mistos conduzirá em um problema de ponto sela, tipicamente existem métodos para tratar diretamente estes problemas como o bem conhecido método de Uzawa com suas diferentes versões que visam uma aceleração do método. Em vez disso, vamos aproveitar a mesma partição e considerar os intervalos $\Omega_j= I_j$ como subdomínios de $\Omega$ e a decomposição de (\ref{FFG}) sobre cada um deles, ou seja,
chamarmos $h_{j} = h\mid_{\Omega_j}$ e $u_{j} = u\mid_{\Omega_j}$ (a solução de (\ref{FFG}) restringida a cada subdomínio $\Omega_j$) para obter formulação fraca local em {\it tempo continuo} e {\it espaço continuo}. Aqui cada solução local pertence aos subespaços $H_1(\Omega_j)$ e $L^2(\Omega_j)$, 
\begin{equation}
\begin{aligned}
\left( \alpha u_j, v \right)_{L^2} - \left( h_j,v_x \right)_{L^2} + h_j(t,b) \cdot v(b) - h_j(t,a)\cdot v(a)= 0, && v\in L^2(I_j),  \\
\left( \partial_t h_j,w \right)_{L^2} + \beta\left( (u_j)^2,w \right)_{L^2} + 
\left( \partial_x u_j,w \right)_{L^2} = \left( \xi,w \right)_{L^2}, && w\in H_1(I_j),  \\
h_j(0) = h_0\mid_{I_j}.
\end{aligned} \label{FFLEC}  
\end{equation}

Agora vamos tirar o índice $h$ em (\ref{EFM_Fraca}) e procurar uma aproximação da solução de cada problema (\ref{FFLEC}) sobre cada espaço de dimensão finita $V_j = \{v\mid\Omega_j, v\in V^h \}$ e $W_j = \{w\mid\Omega_j, w\in W^h \}$, 
\begin{equation}
\begin{aligned}
\left( \alpha \bar{U}_j, v \right)_{L^2} - \left( \bar{H}_j,v_x \right)_{L^2} + \bar{H}_j(t,b) \cdot v(b) - \bar{H}_j(t,a)\cdot v(a)= 0 && v\in V_j,  \\
\left( \partial_t \bar{H}_j,w \right)_{L^2} + \beta\left( (\bar{U}_j)^2,w \right)_{L^2} + 
\left( \partial_x \bar{U}_j,w \right)_{L^2} = \left( \xi,w \right)_{L^2} && w\in W_j  \\
\bar{H}_j(0) = h_0\mid_{\Omega_j}.
\end{aligned} \label{FFLTC}  
\end{equation}

  Como as funções bases (\ref{BaseEFM}) estão definidas sobre cada subdomínio $\Omega_j$ podemos considerar as seguintes bases locais em cada um deles,
\begin{equation}
\begin{aligned}
\varphi_j^1(x) &= \frac{x_{j+1} - x}{\Delta x}, \\
\varphi_j^2(x) &= \frac{x -x_j}{\Delta x}, \\
\psi_j &= 1, 
\end{aligned}\label{BasesLocais}
\end{equation}
para garantir que a solução do sistema (\ref{FFLTC}) seja consistente com a solução de (\ref{EFM_Fraca}) devemos impor condições de acoplamento nas interfaces de cada subdomínio $\Omega_j$, 
\begin{subequations}
\begin{align}
\bar{U}_j(x_{j+1})  = \bar{U}_{j+1}(x_{j+1}), && \forall j \in [1,m], \\
\bar{H}_j(x_{j+1})  = \bar{H}_{j+1}(x_{j+1}), && \forall j \in [1,m]. 
\end{align}  
\end{subequations}\label{CondCons}

  A formulação  mista em {\it tempo discreto} consiste em aproximar a derivada temporal (via Euler recuado, por exemplo) em (\ref{FFLTC}) e procurar as aproximações $\lbrace \bar{U}_j^n,\bar{H}_j^n \rbrace$ nos espaços $ V_j$ \, e \, $W_j$, para cada tempo fixo $t_n$ com $n = 1,\cdots, N$, 
\begin{equation}\label{FFLTD}   
\begin{aligned}
\left( \alpha \bar{U}_j^n, v \right)_{L^2} - \left( \bar{H}_j^n,v_x \right)_{L^2} + \bar{H}_{j+1}^n \cdot v(x_{j+1}) - \bar{H}_j^n \cdot v(x_j)&= 0 && v\in L_{I_j},  \\
\left( \frac{\bar{H}_j^n-\bar{H}_j^{n-1}}{\Delta t},w \right)_{L^2} + \beta\left( (\bar{U}_j^n)^2,w \right)_{L^2} + \left( \partial_x \bar{U}_j^n,w \right)_{L^2} &= \left( \xi,w \right)_{L^2} && w\in H_{I_j},  \\
\bar{H}_j(0) = h_0\mid_{\Omega_j}. 
\end{aligned}
\end{equation}

Em várias famílias de elementos finitos mistos, as funções $w \in W^h$ podem ser descontínuas
em cada interface. Como uma consequência, a tentativa de impor as condições de consistência (\ref{CondCons}) poderia provocar um erro de conservação de fluxo, isto é, (\ref{CondCons}b) não seria satisfeita a menos que a solução aproximada $H$ for constante. Então vamos introduzir as variáveis $l^n_{j1} = h(t_n,x_j)$ e $l^n_{j2} = h(t_n,x_{j+1})$ (multiplicadores de Lagrange) as quais terão a função de conectar, ou acoplar, as soluções nas interfaces. As condições de consistência (\ref{CondCons}) se transformam em,
\begin{align*}
\bar{U}^n_j(x_j)    = \bar{U}^n_{j-1}(x_j),         &&   \bar{U}^n_j(x_{j+1}) &= \bar{U}^n_{j+1}(x_{j+1}),  \\
l^n_{j1}      = l^n_{(j-1),2},          &&   l^n_{j2}       &= l^n_{(j+1),1}.
\end{align*}
 
  Note que a aproximação de $h$ no interior de cada subdomínio $\Omega_j$ não fica mais com a ``responsabilidade'' da continuidade e, portanto, a solução $h$ pode ser aproximada por uma função descontínua, por exemplo, uma função constate por partes. Então estamos procurando as aproximações da solução $h_j$ e do fluxo $u_j$ em cada subdomínio $\Omega_j$ como sendo combinações lineares das bases (\ref{BasesLocais}) dos espaços $W_j$ e $V_j$, respetivamente,
\begin{align*}
\bar{U}_j^n(x) &= U^n_j\varphi_j^1(x) + U^n_{j+1}\varphi_j^2(x),  \\
\bar{H}_j^n(x) &= H^n_{j+\frac{1}{2}}\psi_j(x),
\end{align*}
onde,
\begin{align*}
U^n_j &= \bar{U}_j(t_n,x_j),  \\
U^n_{j+1} &= \bar{U}_j(t_n,x_{j+1}),  \\
H^n_{j+\frac{1}{2}}  &= \bar{H}_j(t_n,x_{j+\frac{1}{2}}).
\end{align*}

Para simplificar a notação e fazer mais claro o algoritmo vamos renomear os coeficientes das aproximações que pertençam a intervalos adjacentes a $I_j$.
\begin{align*}
U_{(j+1),1}^n &= \bar{U}_{j+1}(t_n,x_{j+1}) = U^n_R,  \\
U_{(j-1),2}^n &= \bar{U}_{j-1}(t_n,x_j)     = U^n_L,  \\
l_{(j+1),1}^n &= \bar{H}_{j+1}(t_n,x_{j+1}) = l^n_R,  \\
l_{(j-1),2}^n &= \bar{H}_{j-1}(t_n,x_j)     = l^n_L.  \\
\end{align*}

Seguindo \cite{DoH98,EAJV15}  utilizamos a regra dos trapézios para aproximar as integrais em \ref{FFLTD},
\begin{align*}
\left(\bar{U}_j^n,\varphi_j^1 \right) &= \int_{x_j}^{x_{j+1}} \left( U_{j1}^n\varphi_j^1(x) + U_{j2}^n\varphi_j^2(x) \right)\cdot \varphi_j^1 \approx \frac{\Delta x}{2}U_{j1}^n\cdot
\end{align*}
\begin{align*}
\left(\bar{U}_j^n,\varphi_j^2 \right) &= \int_{x_j}^{x_{j+1}} \left( U_{j1}^n\varphi_j^1(x) + U_{j2}^n\varphi_j^2(x) \right)\cdot \varphi_j^2 \approx \frac{\Delta x}{2}U_{j2}^n\cdot
\end{align*}
\begin{align*}
\left(\bar{H}_j^n,\partial_x\varphi_j^1 \right) &= \int_{x_j}^{x_{j+1}} 
H_{j+\frac{1}{2}}^n\psi_j(x)\cdot \partial_x\varphi_j^1 \approx -H_{j+\frac{1}{2}}^n
\end{align*}
\begin{align*}
\left(\bar{H}_j^n,\partial_x\varphi_j^2 \right) &= \int_{x_j}^{x_{j+1}} 
H_{j+\frac{1}{2}}^n\psi_j(x)\cdot \partial_x\varphi_j^2 \approx H_{j+\frac{1}{2}}^n
\end{align*}
\begin{align*}
\left(\partial_x \bar{U}_j^n,\psi_j \right) &= \int_{x_j}^{x_{j+1}} \partial_x \left( U_{j1}^n\varphi_j^1(x) + U_{j2}^n\varphi_j^2(x) \right)\cdot \psi_j \approx U^n_{j2} - U^n_{j1}\\
\end{align*}
\begin{align*}
\left(\left( \bar{U}_j^n \right)^2,\psi_j \right) &= \int_{x_j}^{x_{j+1}} \left( U_{j1}^n\varphi_j^1(x) + U_{j2}^n\varphi_j^2(x) \right)^2 \cdot \psi_j\approx \frac{\Delta x}{2} \left( (U^n_{j1})^2 + (U^n_{j2})^2\right)
\end{align*}
\begin{align*}
\left( \xi,\psi \right)_{L^2} &= \int_{x_j}^{x^{j+1}} \xi\cdot\psi \approx \frac{\Delta x}{2} \left( \xi_j + \xi_{j+1} \right)
\end{align*}
\begin{align*}
\left( \bar{H}_j^n,\psi \right)_{L^2} &= \int_{x_j}^{x^{j+1}} H_{j+\frac{1}{2}}^n\psi\cdot\psi \approx \Delta x H_{j+\frac{1}{2}}^n.
\end{align*} 

Obs: ``A regra de integração dos trapézios, em conjunto com o
espaço de Raviart-Thomas de menor índice (RT0), é
conveniente no contexto da formulação mista para ao
mesmo tempo acoplar o sistema discreto de equações na
interface e também para resultar em um algoritmo mais
simples de ser resolvido numericamente; para mais
detalhes, ver trabalhos (\cite{DPLW93,DoH98,EAF15})''. \\

  Substituindo as aproximações das integrais de volta em (\ref{FFLTD}) obtemos o seguinte sistema de equações algébrico,
\begin{equation}
\begin{aligned}
&\begin{array}{lr}
l_{j1}^n = l_L^n,  &  l_{j2}^n = l_R^n, \\
U_{j1}^n = U_L^n,  &  U_{j2}^n = U_R^n,  \\ 
\end{array} \\
&\alpha \frac{\Delta x}{2} U^n_{j1} + H^n_j - l^n_{j1} = 0,  \\
&\alpha \frac{\Delta x}{2} U^n_{j2} - H^n_j + l^n_{j2} = 0,  \\
&\frac{\Delta x}{\Delta t}\left( H_{j+\frac{1}{2}}^{n} - H_{j+\frac{1}{2}}^{n-1} \right)
+ \beta \frac{\Delta x}{2}\left(\left( U_{j1}^n  \right)^2 + \left( U_{j2}^n \right)^2\right)
+ U_{j2}^n - U_{j1}^n  = \frac{\Delta x}{2} \left( \xi_j + \xi_{j+1} \right).
\end{aligned}\label{AlgN}
\end{equation}
 
  Para definir um método iterativo para o problema algébrico gerado pela discretização de
elementos finitos mistos e híbridos é preciso substituir as condições de consistência pelas 
condições de transmissão de contorno de Robin (nas interfaces dos elementos)  \cite{Des90,Des91}, ver também \cite{JBBD97}, ou seja, 
\begin{align*}
l_{j1}^n & = l_L^n + \chi_1\left(U_L^n - U_{j1}^n \right),  \\
l_{j2}^n & = l_R^n - \chi_2\left(U_R^n - U_{j2}^n \right). \\
\end{align*}
Substituindo de volta em \ref{AlgN} temos:

\begin{equation}
\begin{aligned}
&l_{j1}^n  = l_L^n + \chi_1\left(U_L^n - U_{j1}^n \right), \\
&l_{j2}^n  = l_R^n - \chi_2\left(U_R^n - U_{j2}^n \right), \\
&\alpha \frac{\Delta x}{2} U^n_{j1} + H^n_j - l^n_{j1} = 0, \\
&\alpha \frac{\Delta x}{2} U^n_{j2} - H^n_j + l^n_{j2} = 0, \\
&\frac{\Delta x}{\Delta t}\left( H_{j+\frac{1}{2}}^{n} - H_{j+\frac{1}{2}}^{n-1} \right)
+ \beta \frac{\Delta x}{2}\left(\left( U_{j1}^n  \right)^2 + \left( U_{j2}^n \right)^2\right)
+ U_{j2}^n - U_{j1}^n  = \frac{\Delta x}{2} \left( \xi_j + \xi_{j+1} \right). 
\end{aligned}\label{SisAlg}
\end{equation}

  Neste momento é importante destacar que no paper {\it Mixed Finite Element Domain Decomposition for Nonlinear Parabolic Problems} os autores \cite{KPP00} provaram a convergência do algoritmo iterativo de elementos finitos mistos e híbridos para problemas da classe,
\begin{equation} 
\begin{aligned}
  u_t  -div(a(x,u) \nabla u) + f(x,t,u,\nabla u) = 0, && (x,t)\in \Omega \times J, \\
  u(x,0) = u_0(x), && x \in \Omega,  \\
  u(x,t) = 0, && x \in \partial\Omega, 
\end{aligned}    \label{Kim}
\end{equation}
com $J=[0,T]$ e $\Omega \in \mathbb{R}^n$ onde a função $a$ é de classe $C^2$ e a função $f$ seja de classe  $C^2$ na terceira e a quarta componente. Podemos ver que esta versão da equação KPZ (onde o termo fonte é uma função suave)  pode ser escrita na forma (\ref{Kim}) considerando $a(t,x)\equiv \nu$ \,  e \, $f(t,x,u,\nabla u) = (\nabla u)^2 + \eta(t,x)$. 

\section{Aproximação numérica}

\subsection{Algoritmo Iterativo}

Para definir o processo iterativo vamos supor conhecidas todas as variáveis no tempo 
$t_{n-1}$ que por sua sua vez servirão como aproximação inicial para o cálculo iterado 
das variáveis no tempo $t_{n}$.
Um grande número de índices está envolvido no sistema algébrico (\ref{SisAlg}) e ainda 
temos que acrescentar um outro índice para o nível de iteração. De aqui para 
frente supõe-se que estamos trabalhando no intervalo $I_j$ e no tempo $t_n$, portanto 
podemos suprimir o índice $j$ e o índice $n$ nas variáveis. A variável $H_{j+\frac{1}{2}}^{n-1}$ 
será chamada de $H_a$ onde o subscrito $a$ significa $anterior$. O sobrescrito $k$ 
será usado para numerar o nível de iteração.
Uma primeira estratégia muito simples para calcular a solução do sistema algébrico (\ref{SisAlg}) é 
o seguinte algoritmo iterativo de ponto fixo  onde as 
variáveis que trazem informação dos elementos adjacentes e os termos não lineares são 
mantidos no nível de iteração $k-1$, dessa forma o sistema de equações resultante será linear como é mostrado a seguir:
\begin{subequations}
\begin{align}
&l^k_1 + \chi_1 U^k_1   =  l^{k-1}_L + \chi_1 U^{k-1}_L  \\
&l^k_2 - \chi_2 U^k_2   =  l^{k-1}_R - \chi_2 U^{k-1}_R  \\
&\left(\alpha\frac{\Delta x}{2}+\chi_1\right) U^k_1 + H^k   =  l^{k-1}_L + \chi_1 U^{k-1}_L   \\
&\left(\alpha\frac{\Delta x}{2}+\chi_2\right) U^k_2 - H^k   =  -l^{k-1}_R + \chi_2 U^{k-1}_R  \\
&\frac{\Delta x}{\Delta t} H^k  + U^k_2  -  U^k_1 
=  \frac{\Delta x}{2} \left( \xi_1 + \xi_2 \right)  - \beta\frac{\Delta x}{2}\left((U^{k-1}_1)^2 
+ (U^{k-1}_2)^2\right) + \frac{\Delta x}{\Delta t}H_a.
\end{align}\label{AlgPF}
\end{subequations}
Na fronteira à esquerda temos a condição de Dirichlet que modifica as equações (\ref{AlgPF}a) e (\ref{AlgPF}c), 
\begin{align}
l^k_1 = h(t_n,a),  \\
\alpha\frac{\Delta x}{2} U^k_1 + H^k = h(t_n,a). 
\label{AlgPFL}
\end{align}  
Na fronteira à direita temos a condição de Dirichlet que modifica as equações (\ref{AlgPF}b) e (\ref{AlgPF}d), 
\begin{align}
l^k_2 = h(t_n,b), \\
\alpha\frac{\Delta x}{2} U^k_2 - H^k = - h(t_n,b).
\label{AlgPFR}
\end{align}  
Podemos escrever o sistema na forma simplificada,
\begin{align}
M \cdot X_j^k = F(X_L^{k-1},X_R^{k-1},\xi_{j1},\xi_{j2},H_a),
\label{EqIt}
\end{align}
onde $X_j = [l_1,l_2,U_1,U_2,H]^T$, $X_L$ e $X_R$ são os vetores das incógnitas em cada intervalo $I_j$, $I_{j-1}$ e $X_{j+1}$, e $F$ é a função vetorial (não linear) da parte direita de (\ref{AlgPF}) e $\xi_{j1} = \xi(t_n,x_j)$ e $\xi_{j2} = \xi(t_n,x_{j+1})$. Desse modo fica totalmente definido o sistema algébrico linear de equações que temos que resolver em cada intervalo $I_j$ para cada tempo $t_n$. Podemos adotar, ainda, uma estratégia Red-Black por atualizar as variáveis nos elementos de índice par, para imediatamente serem usadas para atualizar as variáveis nos elementos de índice ímpar \cite{JV15,EAJV15}. O algoritmo pode ser descrito da seguinte forma:

\subsubsection{Algoritmo}
\begin{enumerate}
\item Inicializar $X_0$ a partir do dado inicial 
\item LOOP1 $t = \Delta t \cdots N\Delta t$
\item LOOP2 até que $Erro \leq Tolerancia$
\begin{enumerate}
\item Atualizar $X_1^k$ e $X_m^k$ (Valores na fronteira) usando (\ref{EqIt}) com (\ref{AlgPFL} e \ref{AlgPFR}).
\item Atualizar as variáveis de sub-índice par ($X_{2j}^k$) usando (\ref{EqIt}).
\item Atualizar as variáveis de sub-índice ímpar ($X_{2j-1}^k$) usando (\ref{EqIt}).
\item Calcular o erro relativo $Erro = \frac{\parallel X^k-X^{k-1} \parallel}{\parallel X^k \parallel}$.
\end{enumerate}
\item inicializar o próximo processo iterativo $X^0 = X$.
\end{enumerate}

Todos os códigos foram programados em MATLAB.

\section{Experimentos numéricos}

Vamos testar o algoritmo computacional EFMH_KPZ para um exemplo simples. Consideremos o crescimento laminar para um {\it stomatrolite laminae} \cite{BBHW00}. Nesse modelo de equação KPZ, no lugar do ruído branco, aparece uma função constante, i.e.,
\begin{equation}
\partial_t h = \nu \partial^2_x h + \lambda \left( 1+\partial_x h  \right)^2 + v. \label{EQKPZ1}
\end{equation}   

Esta é uma equação determinística e a sua solução, para um dado inicial parabólico específico é
dado por,
\begin{equation}
  h(t,x) = A + (v+\lambda)t - \frac{\lambda}{\nu}\log(2\lambda t + B) 
  - \frac{(x-x_0)^2}{2\lambda t + B}. \label{Sol_Est}
\end{equation}  
onde $A$, $B$, $v$ e $x_0$, são constantes. Neste contexto, mantendo a força externa constante, vamos testar o código para diferentes dados iniciais esperando ver um crescimento lateral das interfaces.\\
\\
{\bf Experimento 1:} {\it Vamos resolver o problema de valor inicial e de contorno para a equação (\ref{EQKPZ1}) tomando $A=-1$, $B=1$, $x_0=0$ $v=1$, $\nu=1$ e $\lambda =1$. Consideramos condições de fronteira consistentes com a solução  (\ref{Sol_Est}) e  como dado inicial a parábola obtida ao avaliar (\ref{SisAlg}) em $t = 0$. Vamos rodar o código EFMH_KPZ.m para obter soluções aproximadas no tempo $T=0.1$, $T=0.5$ e  $T=1$. A partição do intervalo $[-1,1]$ será de $64$ pontos ($m=64$). A relação entre $\Delta x$ e $\Delta t$ será $\Delta x / \Delta y = C$, onde $C$ é uma constante que dependerá dos parâmetros $\nu$ e $\lambda$ mas no caso da seleção feita para este experimento particular observamos que a escolha $C=1/16$ é suficiente para garantir uma boa aproximação. Vamos usar \, $\chi_1 = \chi_2 =1$.}\\
\\
\begin{equation}
\left\{
\begin{aligned}
&\partial_t h = \nu \partial^2_x h + \lambda \left( 1+\partial_x h  \right)^2 + v, && \forall t\in[0,T],&\,\forall x\in[-1,1],\\
&h(0,x) = x(1-x),\\
&h(t,-1)= h(t,1) = t - log(2t+1)-\frac{1}{2t+1} -1.
\end{aligned}
\right.\label{Exp1}
\end{equation}

\begin{figure}[h!]
\includegraphics[scale=.6,trim = 4cm 0 0 0]{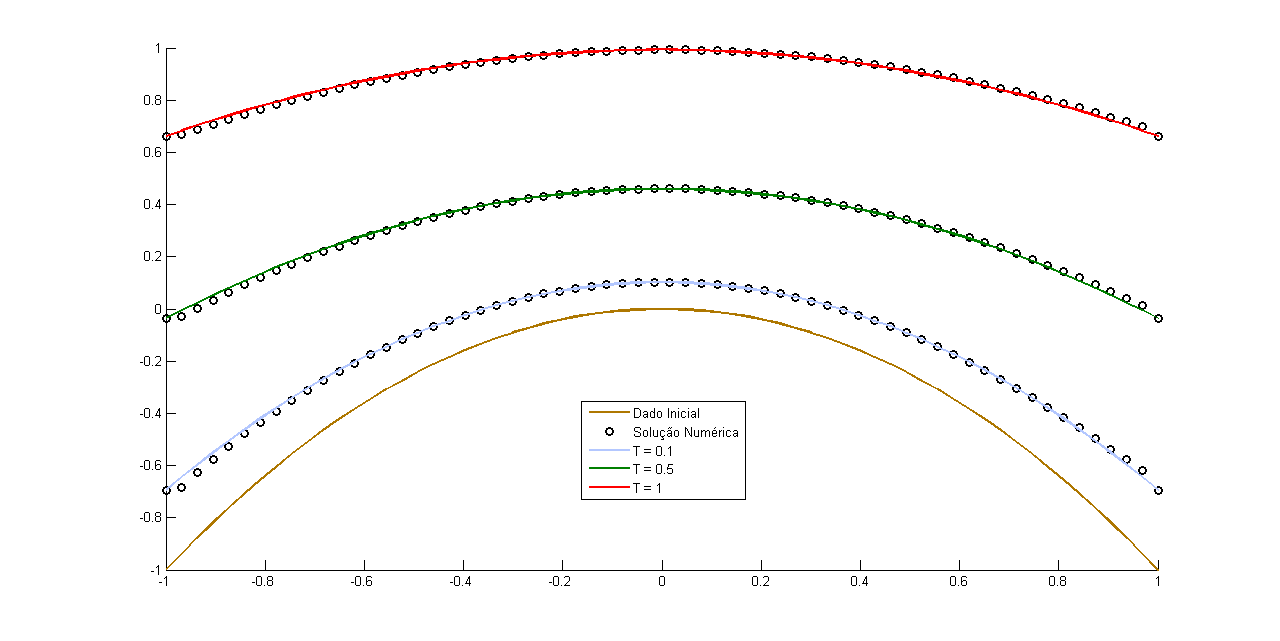} 
\caption{Solução aproximada da equação KPZ determinística para o crescimento laminar do estromatólito para $T=0.1$, $T=0.5$ e  $T=1$ com condições de contorno consistentes com a solução (conhecida). A solução numérica  está acompanhando a solução exata.} \label{Fig_Est_CFC}
\end{figure}

A motivação fundamental da escolha deste exemplo se dá pelo fato de ser possível ``testar o código'' em uma comparação direta com a solução exata (explícita) conhecida do problema (\ref{Sol_Est}). Assim, podemos realizar uma comparação simples, e efetiva, entre as soluções exata e aproximada nessa etapa preliminar de verificação da implementação do procedimento numérico proposto. Os experimentos numéricos mostram que a solução aproximada obtida está, de fato, em boa concordância com a solução exata, como mostrado na Figura \ref{Fig_Est_CFC}. O tempo real de computo foi $3min:06seg$ para $T=0.1$, $5min:12seg$ para $T=0.5$ e $25min:27seg$ para $T=1$ Cumpre ressaltar que temos ciência de que não se tratar de uma rigorosa demonstração matemática, mas sim de uma ``boa'' evidência numérica dos cálculos e aproximações realizadas. Porém, tais resultados motivam em nosso entendimento um estudo de análise numérica mais rigoroso. \\
   
\begin{figure}[h!]
\includegraphics[scale=.6,trim = 4cm 0 0 0]{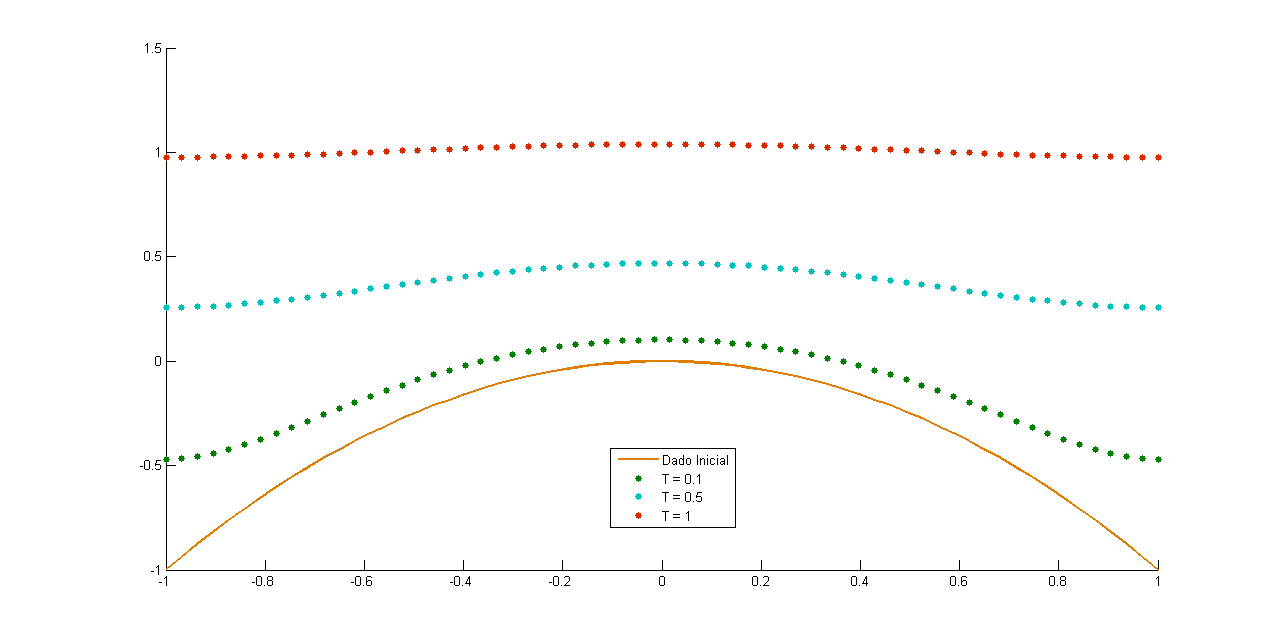} 
\caption{Solução aproximada da equação KPZ determinística para o crescimento laminar do Estromatólito para $T=0.1$, $T=0.5$ e  $T=1$ com condições de contorno periódicas.}
\label{Fig_Est_CFP}
\end{figure}

{\bf Experimento 2 :} {\it Repetimos o Experimento 1 usando os mesmos valores para os parâmetros, só que desta vez vamos impor condições de fronteira periódicas. Neste caso a solução exata não é conhecida. Na Figura \ref{Fig_Est_CFP} são apresentados os gráficos das soluções aproximadas calculadas pelo código EFMH_KPZ.m para tempos $T=0.1$, $T=0.5$ e  $T=1$.}\\

  Uma pequena variação do problema anterior é tomar condições de fronteira periódicas. Para este caso  não conhecemos a solução exata, mas, dado que os valores dos parâmetros do problema não variaram, esperamos que o comportamento da solução numérica obtida seja similar a aquela no exemplo anterior como, de fato, ocorre. Devemos observar que nos bordos o crescimento é mais rápido como produto das contribuições dos crescimentos à direita e à esquerda, além da contribuição da difusão. A solução numérica exibe este efeito como podemos ver na Figura \ref{Fig_Est_CFP}. Assim podemos concluir que os resultados neste estudo também são consistentes e bons\\  
   
{\bf Experimento 3 :} {\it Consideremos os mesmos parâmetros do Experimento 1 para $T = 1$. Rodamos o código EFMH_KPZ.m com $\chi_1 = \chi_2 = 0.1$ e fazemos um estudo de refinamento de malha variando a quantidade de pontos na partição, $m = 2^7,...,2^{10}$.} \\

Esperamos também que a consistência e a estabilidade do método estejam estreitamente associadas com a relação entre o tamanho de passo de tempo e a norma da partição no espaço. Além disso, sendo que o nosso método iterativo é basicamente um algoritmo de iteração de ponto fixo, a convergência para a solução desejada dependerá da aproximação inicial e de uma seleção dos parâmetros que façam com que a equação de iteração seja uma contração.\\

\begin{table}
\centering
\renewcommand{\arraystretch}{1}
\begin{tabular}{c|c|c} 
 $m$  & $\Delta x$ & $E$\\ 
\hline  $128$  & $ 0.0078125$     &  $2.6171040659\cdot 10^{-4}$   \\
\hline  $256$  &  $0.00390625$    &  $1.1188163686\cdot 10^{-4}$   \\
\hline  $512$  &  $0.001953125$   &  $5.2179146103\cdot 10^{-5}$   \\
\hline  1024   &  $0.0009765625$  &  $2.5335808615\cdot 10^{-5}$  \\
\end{tabular}
\label{Tab:Exp3}
\caption{Variação do erro relativo quando refinamos a malha. O experimento mostra que o procedimento numérico fornece aproximações com uma taxa de convergência numérica de primeira ordem.}
\end{table}

\begin{figure}[h!]
\centering
\includegraphics[scale=.4]{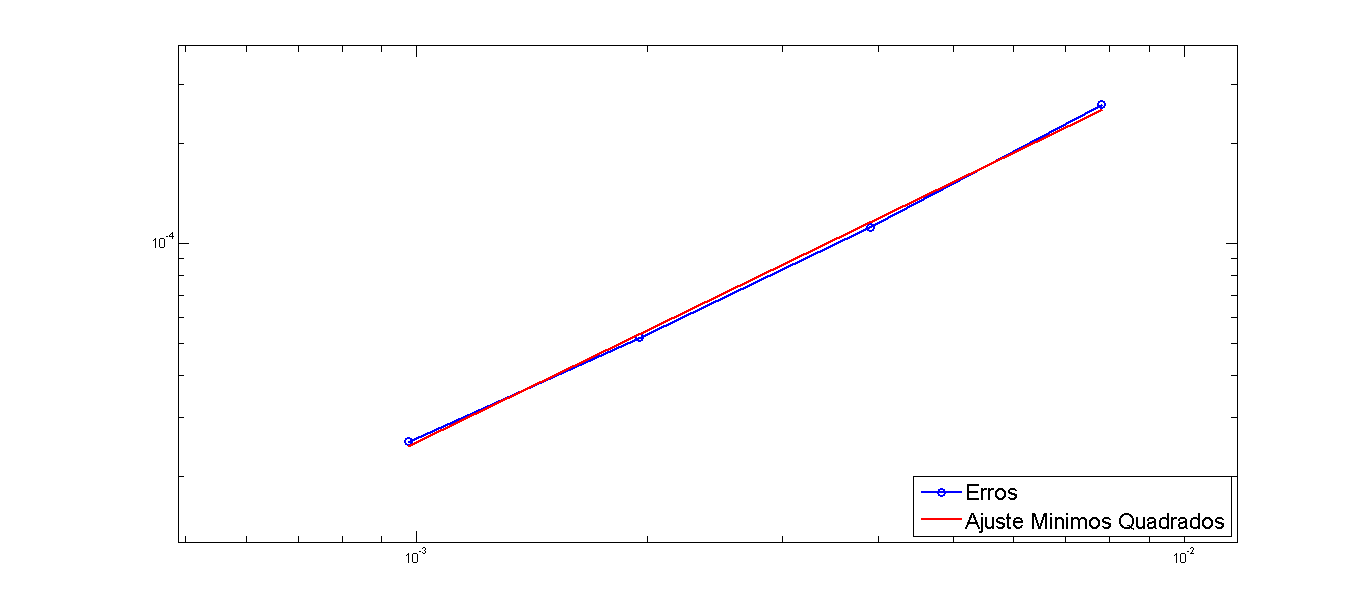}
\caption{Ajuste quadrático dos erro  obtidos no Experimento 3 ao refinar a malha variando a quantidade de pontos na partição, $m = 2^7,...,2^{10}$.} 
\label{Fig:ErrLL}
\end{figure}

  Os resultados do Experimento 3 fornecem alguns bons indicativos de ``convergência numérica\textbf{''}, pois o erro decresce com ordem $O(\Delta x)$, como pode ser visto na Tabela \ref{Tab:Exp3}3 e a Figura \ref{Fig:ErrLL}. Para ver resultados teóricos de convergência ver (\ref{Kim}).

\begin{figure}[h!]
\centering
\begin{tabular}{lll}
\hspace{-1.5 cm}
\subf{\includegraphics[scale=.3]{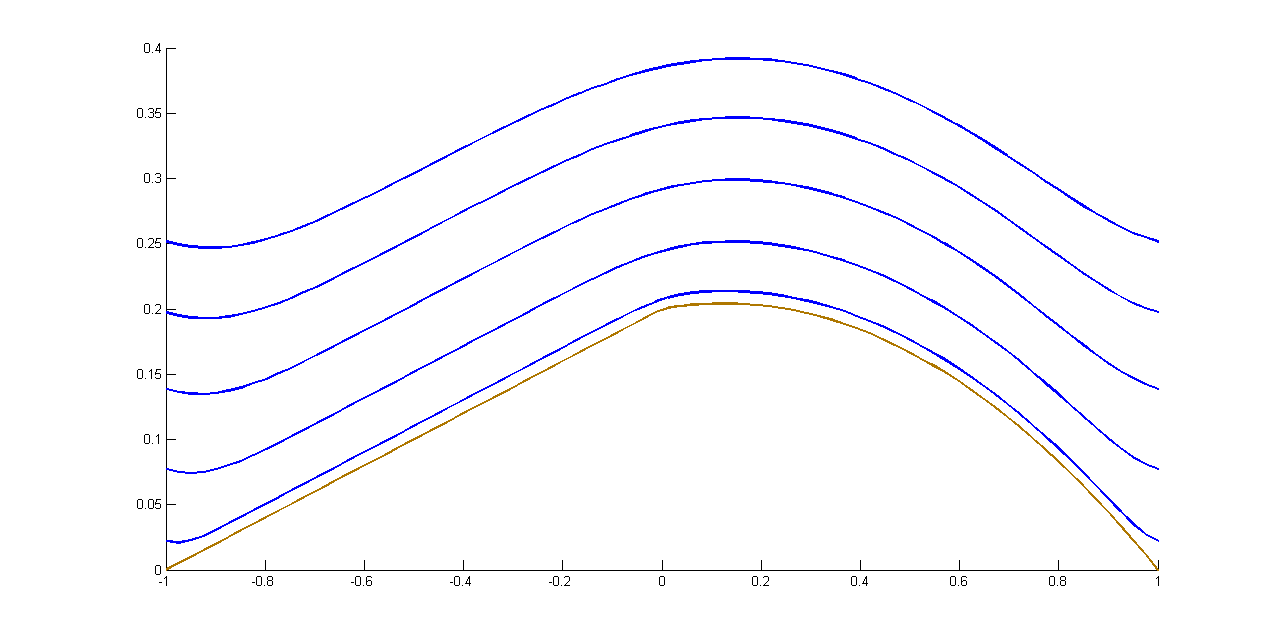} }
     {Dado 1}
&
\hspace{-1 cm}
\subf{\includegraphics[scale=.3]{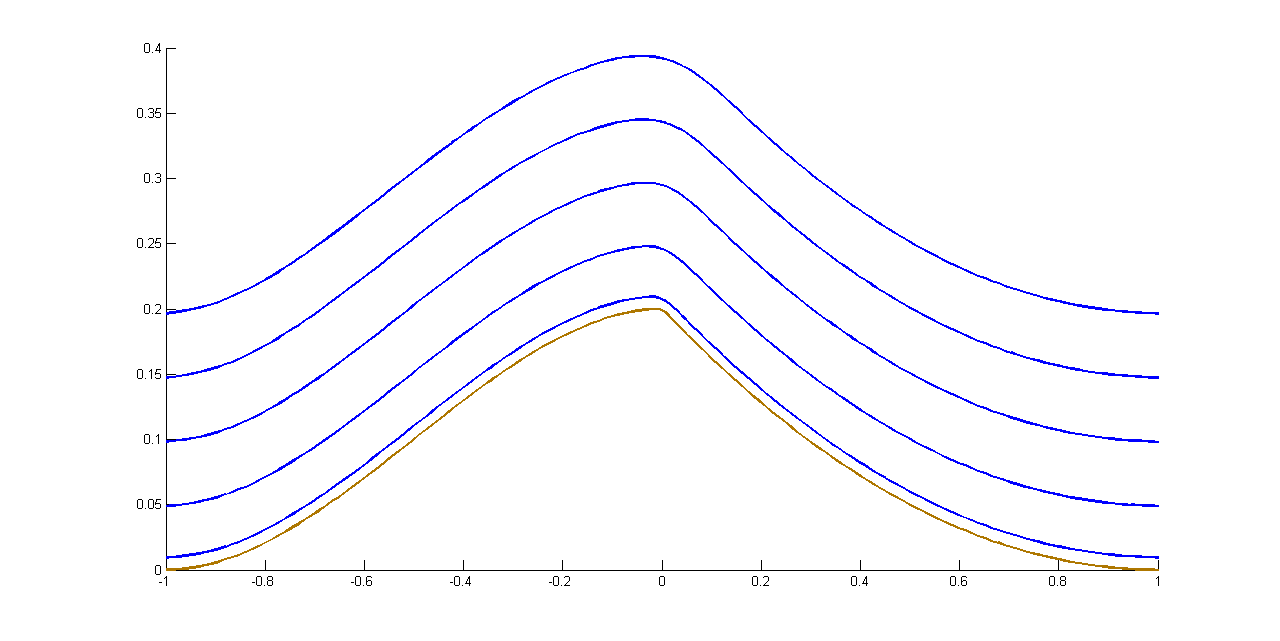} }
     {Dado 2} 
\\
\hspace{-1.5 cm}
\subf{\includegraphics[scale=.3]{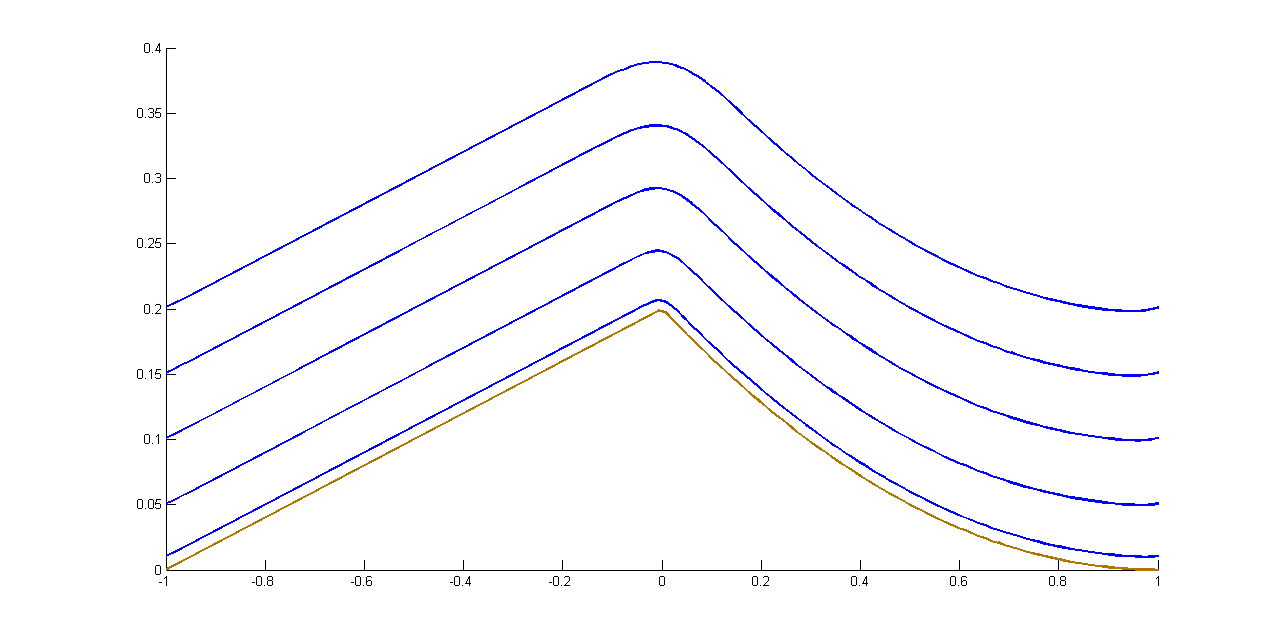} }
     {Dado 3}
&
\hspace{-1 cm}
\subf{\includegraphics[scale=.3]{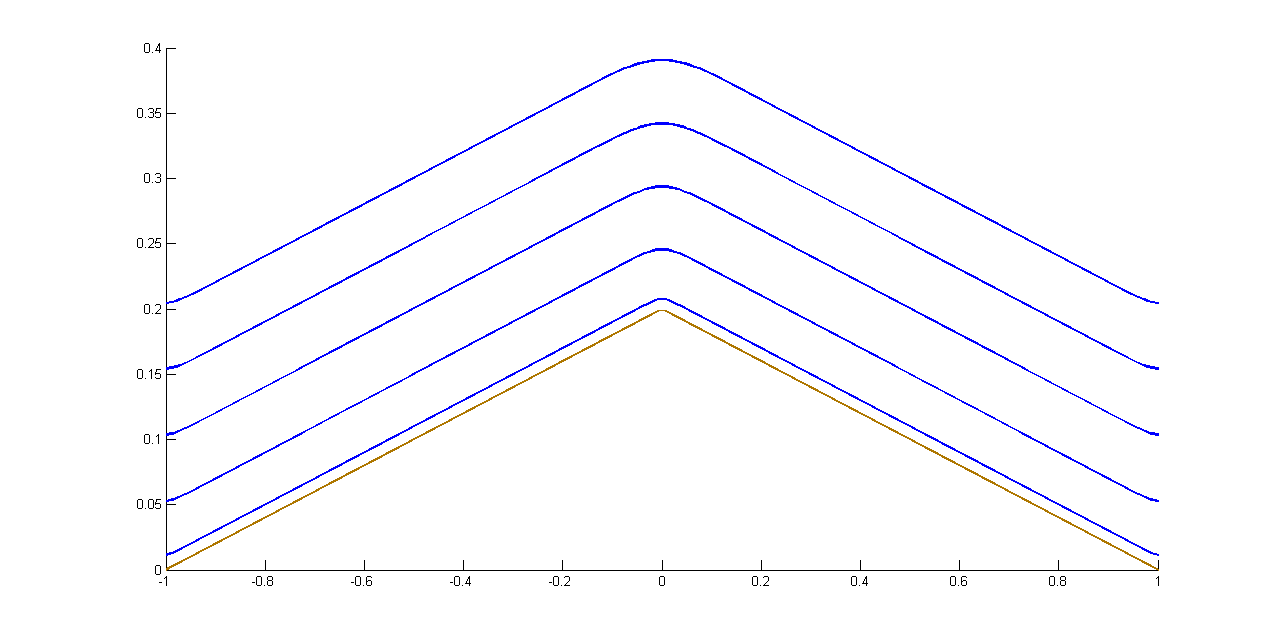} }
     {Dado 4} 
\\
\hspace{-1.5 cm}
\subf{\includegraphics[scale=.3]{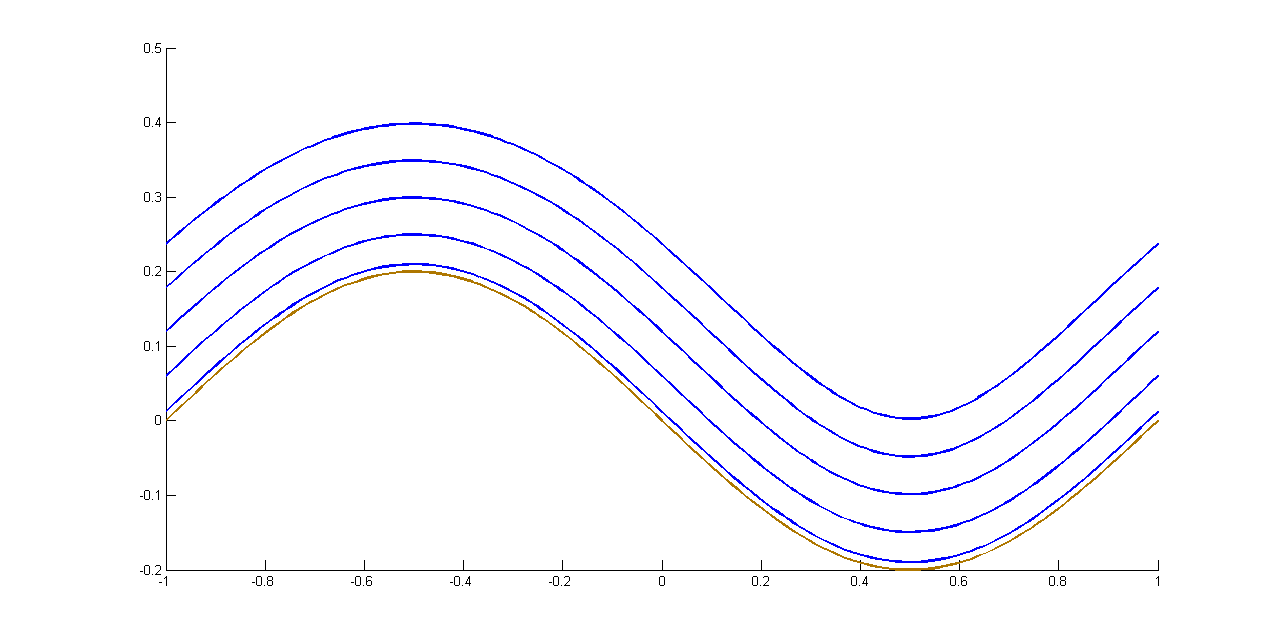} }
     {Dado 5} 
&
\hspace{-1 cm}
\subf{\includegraphics[scale=.3]{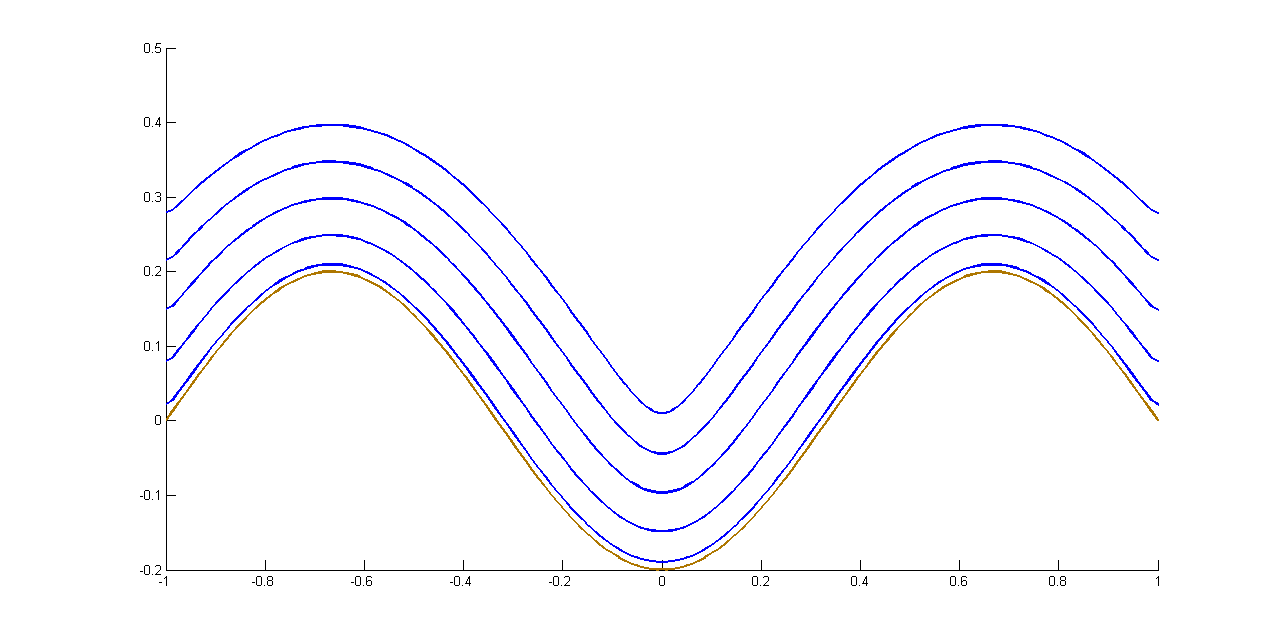} }
     {Dado 6} 
\\
\end{tabular}
\caption{Soluções aproximadas da equação KPZ determinística para seis dados iniciais diferentes. As soluções exibem crescimento lateral, como é esperado.}
\label{fig:6dados}
\end{figure}

{\bf Experimento 4 :} {\it Usando os mesmos valores do Experimento 2 para os parâmetros, calculamos as soluções aproximadas da KPZ para diferentes dados iniciais, que são mostrados na Figura \ref{fig:6dados}.}    \\

Com o objetivo de testar a efetividade do método utilizado, e do código implementado, avaliamos também o desempenho do procedimento na aproximação de uma equação KPZ para diferentes tipos de dado inicial. A Figura \ref{fig:6dados} mostra a estrutura dessas soluções obtidas para seis dados iniciais diferentes. Podemos notar que o comportamento das soluções numéricas obtidas indicam que a interface apresenta um crescimento lateral, o qual é um efeito esperado para a solução desse modelo KPZ.\\
  \\

\chapter{Um estudo de mollifiers e renormalização aplicado num modelo KPZ}\label{cap4}

  Neste capítulo vamos reforçar experimentalmente a idéia de que o problema de valor inicial e
   de contorno para a equação KPZ  é mal posto e que a mesma equação precisa ser renormalizada
   para fazer sentido como equação do modelo BD. Desta maneira, mostraremos o que foi sugerido na 
  secção (\ref{HC}), i.e. as  soluções  da equação KPZ renormalizada se aproximam à transformada de 
  Hopf-Cole da solução da equação estocástica do calor com ruído branco multiplicativo.

\section{Amolecimento do ruído branco}\label{SecAmol}
  
  Existem diferentes maneiras de amolecer o ruído branco bidimensional (espaço-tempo). Uma delas 
  seria amolecer apenas no espaço e a outra nas duas componentes. Veremos as duas variantes.\\
  
\begin{figure}[h!]
\centering
\begin{tabular}{lll}
\hspace{-3 cm}
\subf{\includegraphics[trim=2cm 2cm 3cm 4cm, scale=.4]{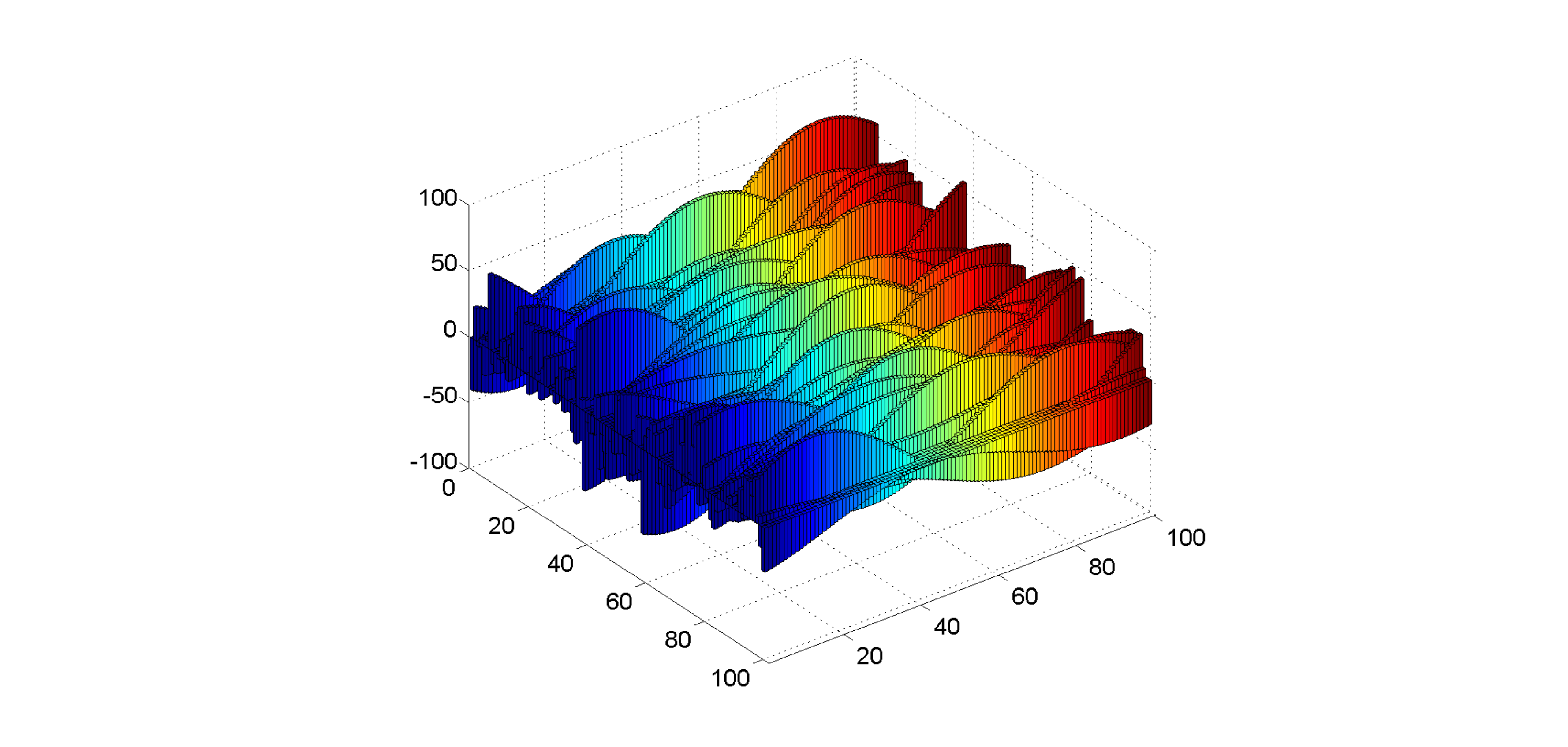}}
     {Ruído Amolecido $\kappa = 1/2$.}
&
\hspace{-4 cm}
\subf{\includegraphics[trim=-3cm 2cm 3cm 4cm, scale=.4]{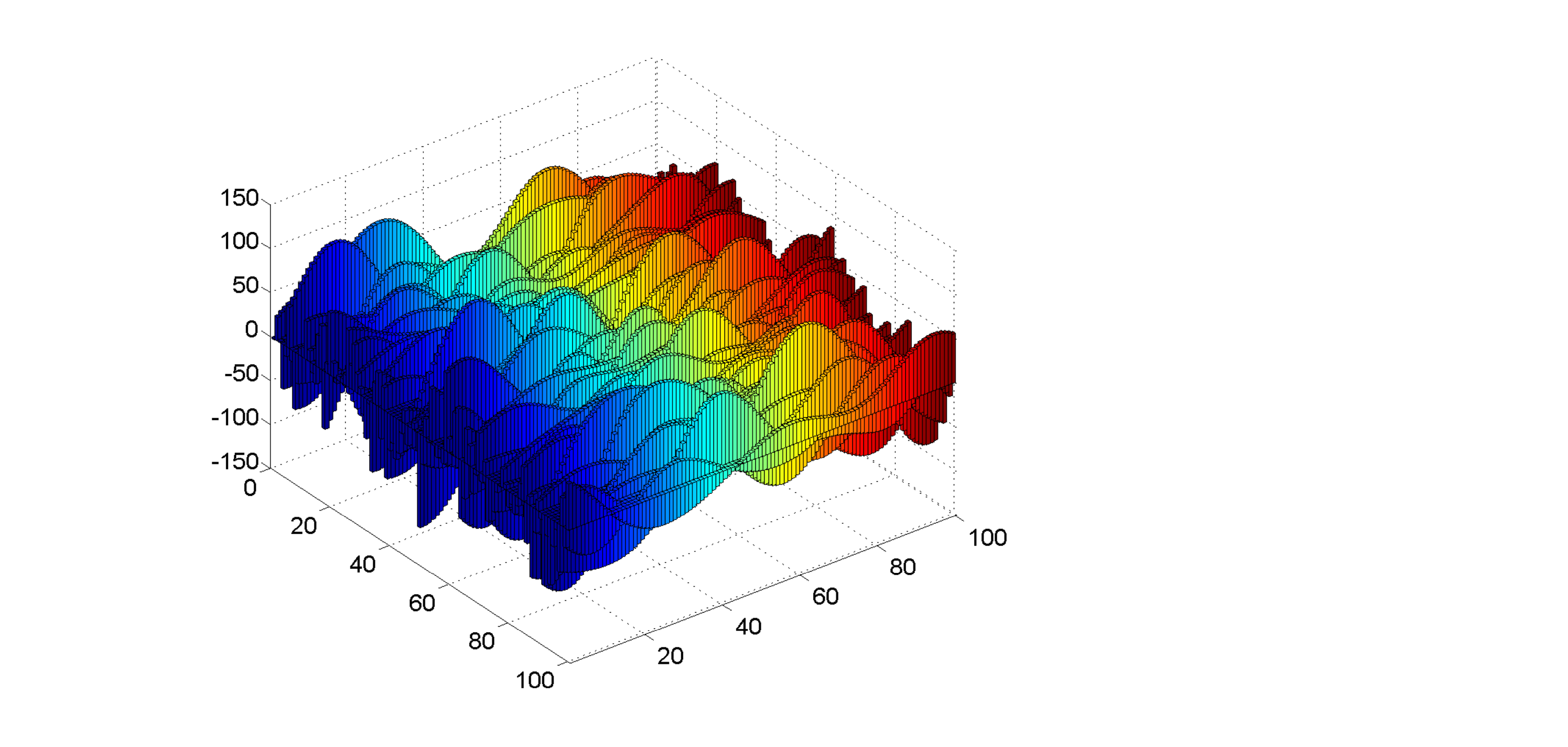}}
     {Ruído Amolecido $\kappa = 1/4$} 
\\
\end{tabular}
\caption{Versões amolecidas do ruído branco para dois valores do parâmetro $\kappa = 1/2, 1/4$ usando como \textit{mollifier} a função (\ref{mol1}).}
\label{AMOL1}
\end{figure}

  Para efetuar o amolecimento do ruído branco $\xi$ na componente espacial adotaremos a mesma 
  metodologia utilizada em \cite{Hai11} por meio da qual obtemos uma versão  $\xi_{\kappa}$ do 
  ruído branco $\xi$ (amolecido no espaço), usando como \textit{mollifier} uma função  
  $\varphi : \mathbb{R} \longrightarrow \mathbb{R}$, onde $\varphi$ é  par, suave, tem suporte 
  compacto e satisfaz $\varphi(0)=1$.\\
   Definindo 
  \begin{equation}
  \xi_{\kappa,k} = \xi_k \cdot \varphi(\kappa n), \label{MolR}
  \end{equation}
  onde $\xi_n$ é a k-ésima componente de Fourier de $\xi$ (\ref{RB2D}). Estes são ruídos 
  brancos complexos que satisfazem $\xi_k = \bar{\xi}_{-k}$. A Figura \ref{AMOL1} 
  mostra os perfis do ruído branco amolecido para diferentes valores do parâmetro $\kappa$ 
  utilizando o \textit{mollifier} (\ref{mol1}) no tempo e no espaço.\\
  
  Para o amolecimento  no espaço e no tempo seguimos a mesma ideia de Hairer. Primeiramente, 
  amolecemos cada componente de cada modo de Fourier de $\xi$ usando uma função $\varphi$ 
  que satisfaz as condições mencionadas anteriormente, ou seja, cada ruído branco escalar 
  no tempo (\ref{RB1D}) vai ser previamente amolecido. Seguidamente, repetimos o mesmo 
  procedimento no espaço usando a mesma função. Em princípio poderia ser qualquer outra função 
  satisfazendo essas condições porém, por mais simplicidade optamos por usar a mesma.
  
    A Figura \ref{AMOL2} mostra os perfis do ruído branco amolecido no espaço e no tempo 
    para diferentes valores do parâmetro $\kappa$ usando o \textit{mollifier} (\ref{mol1}):
    
\begin{figure}[h!]
\centering
\begin{tabular}{lll}
\hspace{-1.5 cm}

\subf{\includegraphics[scale=.3]{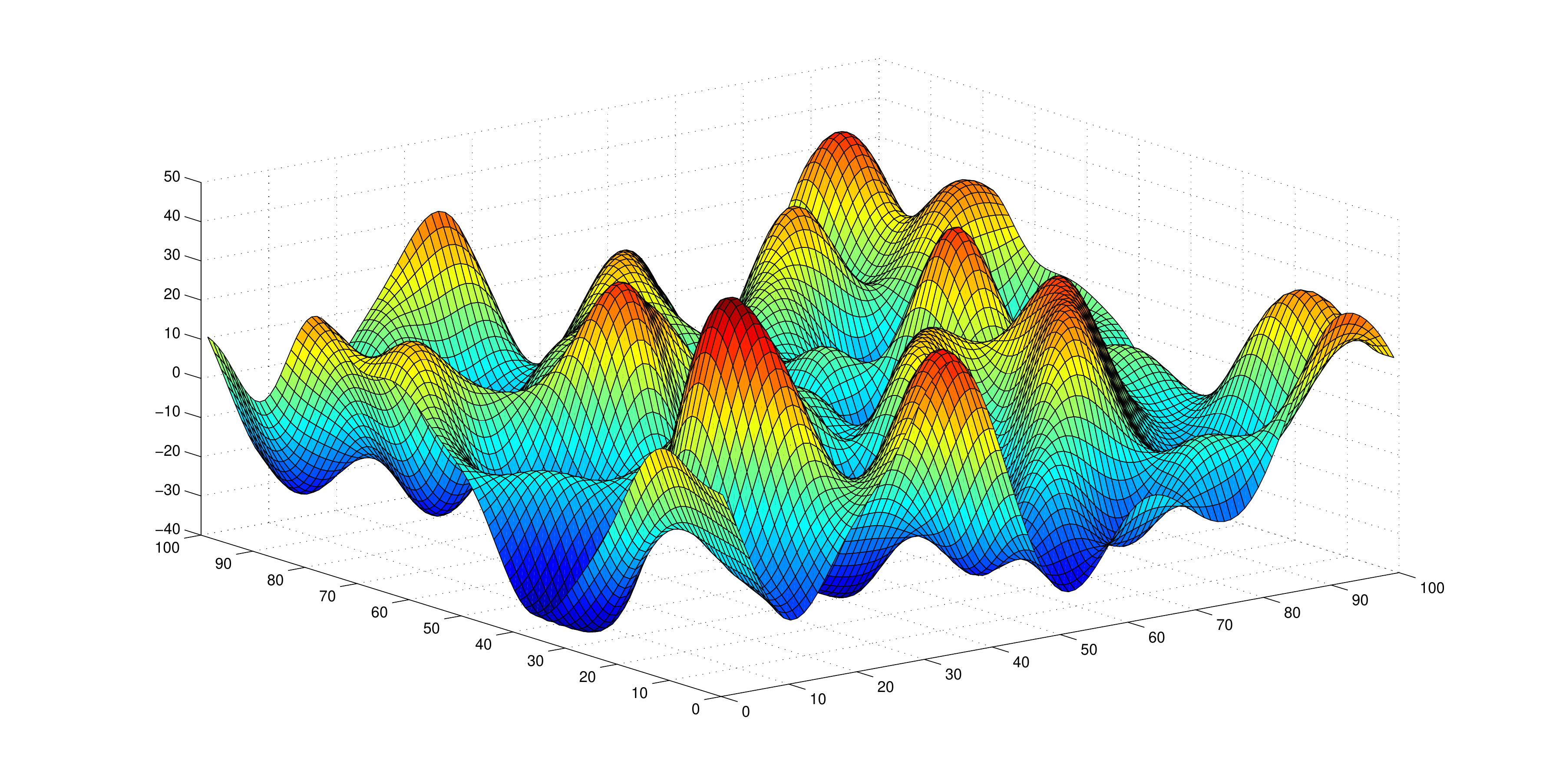} }
     {Ruído Amolecido $\kappa = 1/8$}
&
\hspace{-1 cm}
\subf{\includegraphics[scale=.3]{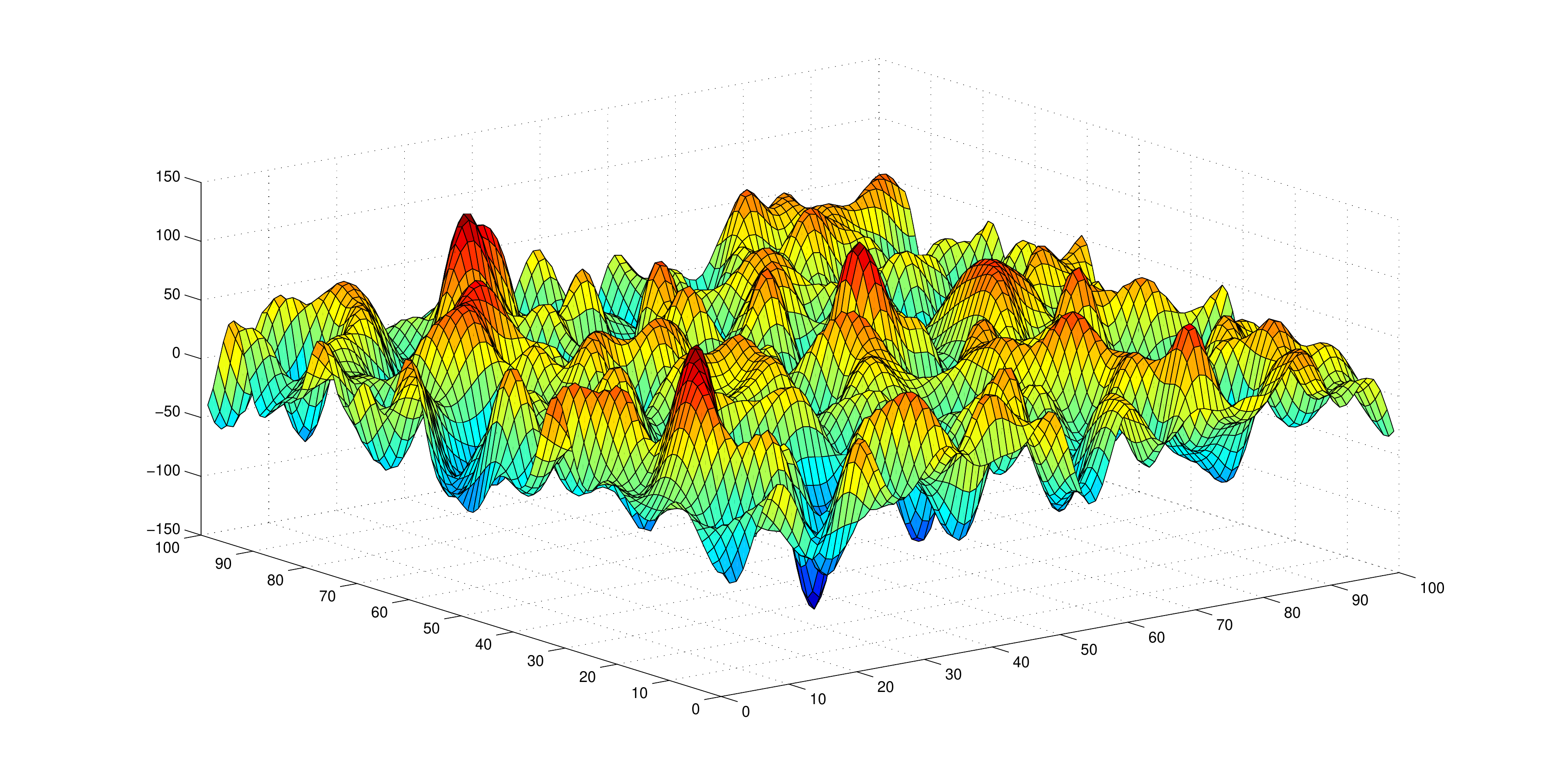} }
     {Ruído Amolecido $\kappa = 1/16$}
\\
\end{tabular}
\caption{Versões amolecidas de ruído branco para quatro valores do parâmetro $\kappa$ usando como \textit{molifier} a função \ref{mol1}.}
\label{AMOL2}
\end{figure}

\section{Divergência das constantes}  

\indent No Capitulo \ref{cap2} na Secção \ref{HC} foi sugerido que 
a equação KPZ precisa  ser renormalizada para fazer sentido. Em
outras palavras, problemas de valor inicial e condições de contorno
seriam mal postos quando considerada a forma clássica da equação 
KPZ. Nessa mesma secção obtemos algumas evid\^encias de que, ao
aplicar um amolecimento sobre o ruído branco, a solução da equação
poderia, de fato, estar próxima da transformada de Hopf-Cole da
solução da equação estocástica do calor se adicionamos uma constate
de renormalização. Assim, passamos a considerar a equação,

\begin{equation}\label{KPZMOL}
\partial_t h_{\kappa} = \nu \partial^2_x h_{\kappa} + \frac{\lambda}{2} \left[\left( \partial_x h_{\kappa} \right)^2 - C_{\kappa} \right] + \eta_{\kappa}.
\end{equation} 

\indent É fácil ver que a equação (\ref{KPZMOL}) é de fato a equação KPZ clássica se considerarmos a mudança de variáveis
\begin{equation}
  \hat{h}_{\kappa}(t,x) = h_{\kappa}(t,x) + \frac{\lambda}{2} \cdot C_{\kappa}\cdot t,
  \label{Cons}
\end{equation}
   o qual quer dizer que a solução da equação KPZ renormalizada vai coincidir com a solução da equação 
   KPZ clássica acrescido de um deslocamento temporal que depende do suporte do \textit{mollifier} 
   o qual diminui quando o parâmetro $\kappa \longrightarrow 0$. Também queremos verificar
    experimentalmente o resultado obtido por Hairer em \cite{Hai11} de que o processo limite quando 
    $\kappa \longrightarrow 0$ é independente da escolha do \textit{mollifier}.\\ 
\\
  Os passos a seguir são:
  
  \begin{itemize} \label{pasos} 
  \item[1)] Aproximar a solução da equação estocástica do calor com ruído 
           branco multiplicativo  utilizando o algoritmo  \textit{Milsteim.m}  introduzido no 
           Capítulo \ref{cap33}.
  \item[2)] Aplicar o método de elementos finitos mistos e  híbridos com decomposição 
           de domínio na equação KPZ clássica considerando a mesma realização do ruído branco usada no 
           passo anterior substituída por uma versão amolecida no espaço (veja \ref{AMOL1}), 
           e comparar as soluções obtidas por cada um destes algoritmos esperando obter trajetórias 
           com perfis similares acrescido de algum deslocamento temporal $C_{\kappa} \cdot t$.
  \item[3)] Diminuir progressivamente o valor do parâmetro de amolecimento do \textit{mollifier} 
           $\kappa \longrightarrow 0$ para tentar conferir que as constantes $C_{\kappa}$ crescem e 
           que os perfis das duas soluções ficam cada vez mais próximas.
  \item[4)] Repetir os passos $1,2,3$ variando o \textit{mollifier} e verificar que o comportamento não 
           varia qualitativamente. 
\end{itemize} 
  
\begin{figure}[h!]
  \centering
\includegraphics[scale=.4,trim = 5cm 0 0 0]{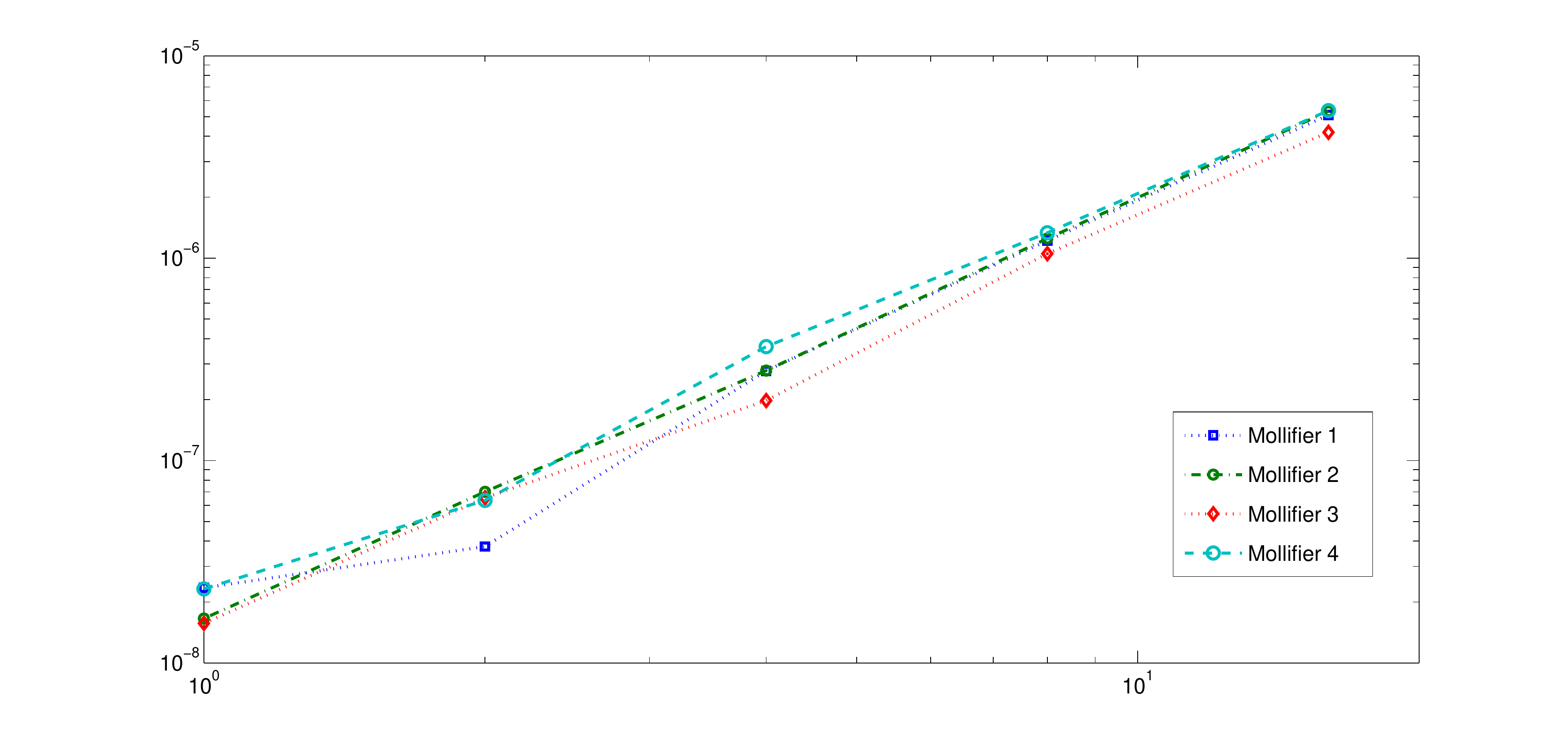} 
  \caption{Crescimento exponencial das constantes de renormalização obtidas para valores do parâmetro $\kappa = 1, 1/2, 1/4, 1/8, 1/16$ dos \textit{mollifier} mostrados na figura \ref{Molifiers}.}
\label{ConstMol} 
\end{figure}

\begin{figure}
  \centering
  \includegraphics[scale=.3,trim = 4cm 0 0 0]{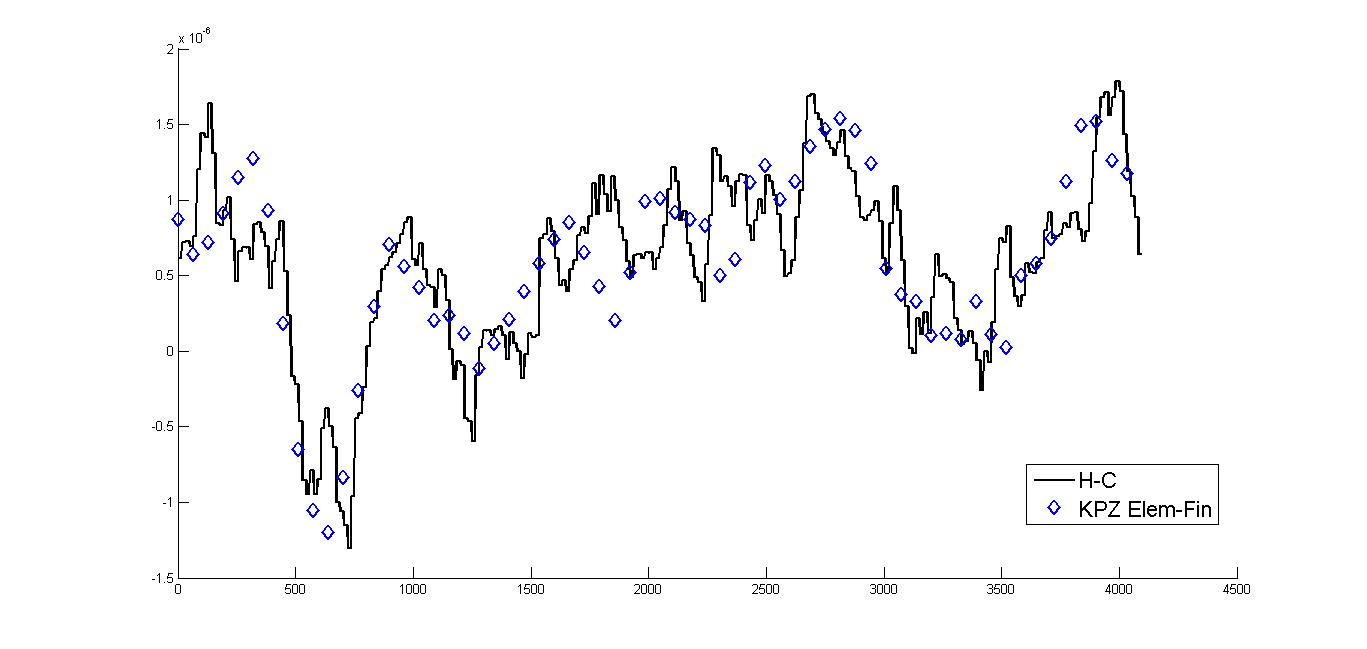}
  \includegraphics[scale=.3,trim = 4cm 0 0 0]{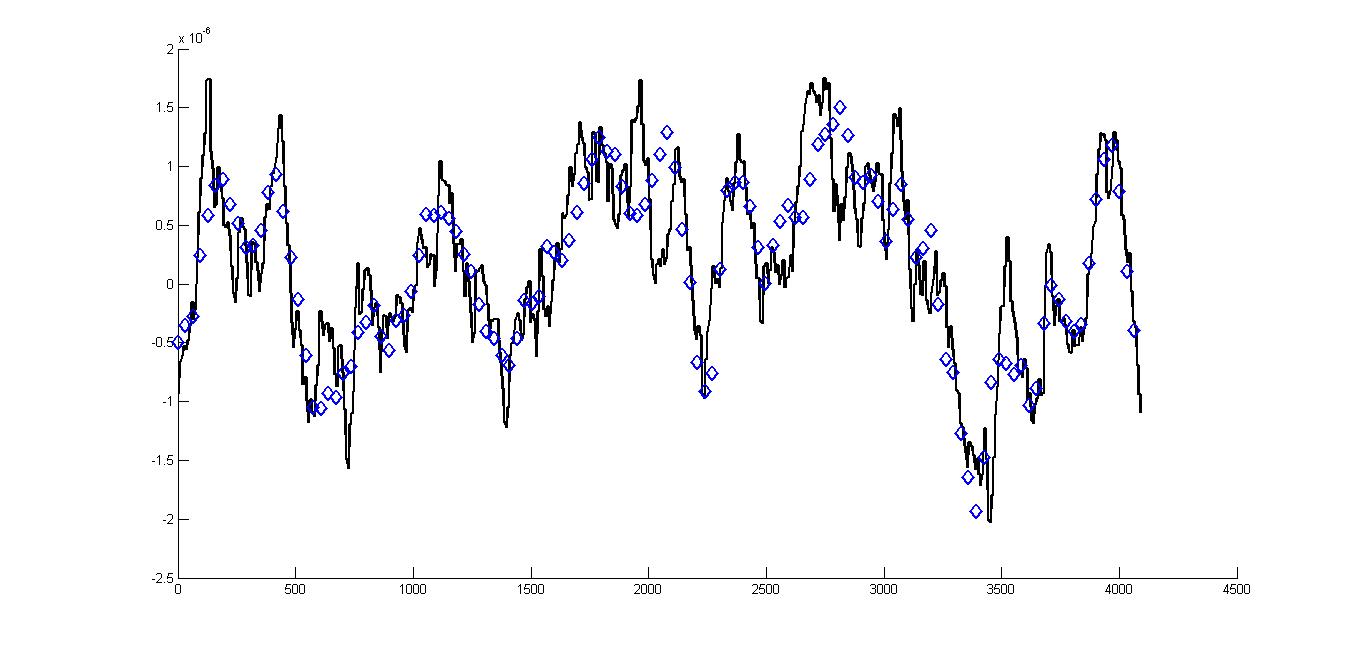} 
  \includegraphics[scale=.3,trim = 4cm 0 0 0]{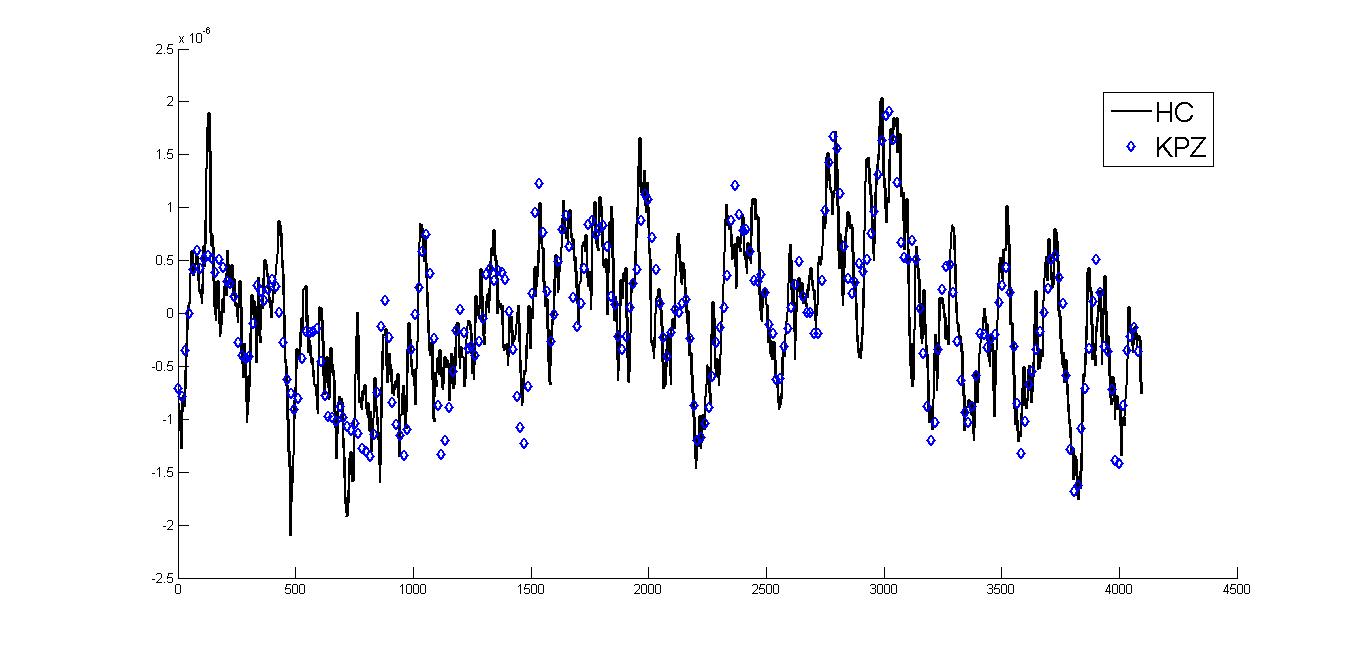}  
  \caption{Perfis das soluções da equação Estocástica do calor transformada via Hopf-Cole junto
           à solução da equação KPZ com ruído branco amolecido calculadas até os tempos
            $\frac{1}{256^2},\frac{1}{512^2},\frac{1}{1024^2}$}
  \label{HCKPZ}
  \end{figure}

\begin{figure}[h!]
  \centering
  \includegraphics[scale=.3,trim = 4cm 0 0 0]{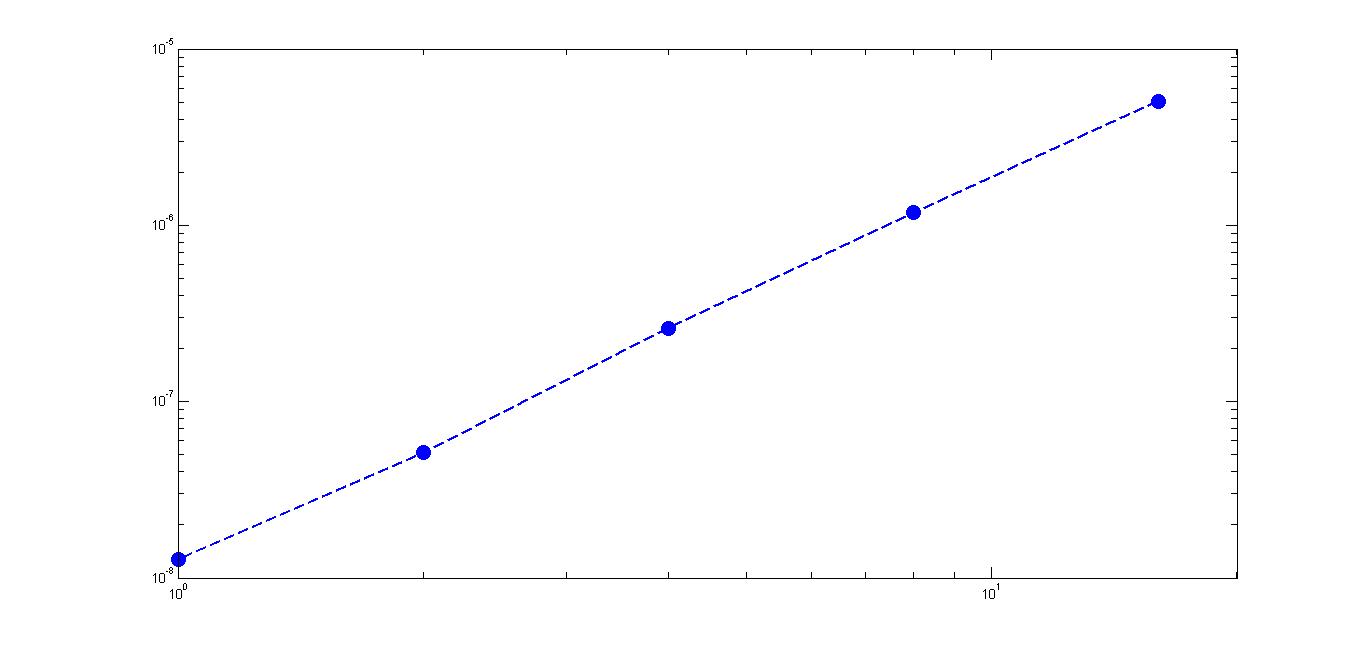} 
  \caption{Crescimento exponencial das constantes de renormalização obtidas para valores do parâmetro $\kappa = 1, 1/2, 1/4, 1/8, 1/16$ do \textit{mollifier} (\ref{mol1})}
  \label{Const} 
  \end{figure}

\begin{figure}[h!]
\centering
\begin{tabular}{lll}
\hspace{-1.5 cm}
\subf{\includegraphics[scale=.3]{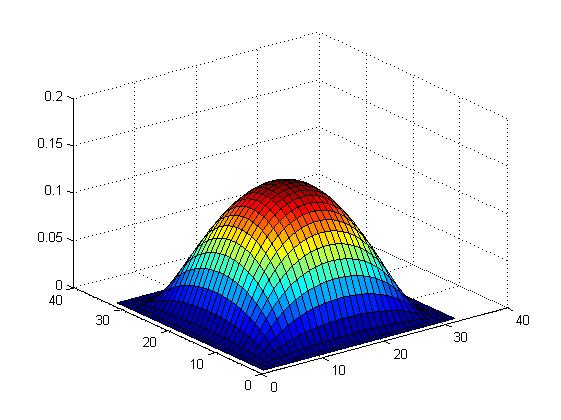} }{\textit{Mollifier} 1}
&
\hspace{-1 cm}
\subf{\includegraphics[scale=.3]{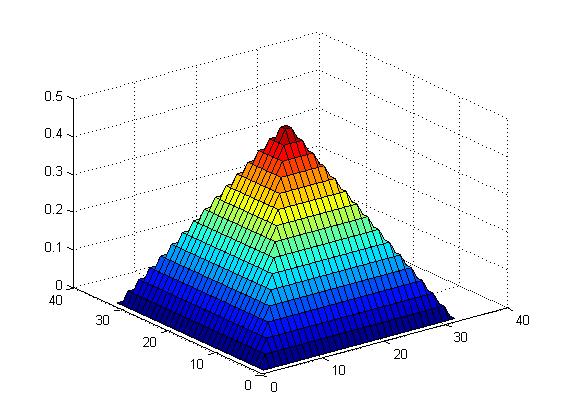}  }{\textit{Mollifier} 2} 
\\
\hspace{-1.5 cm}
\subf{\includegraphics[scale=.3]{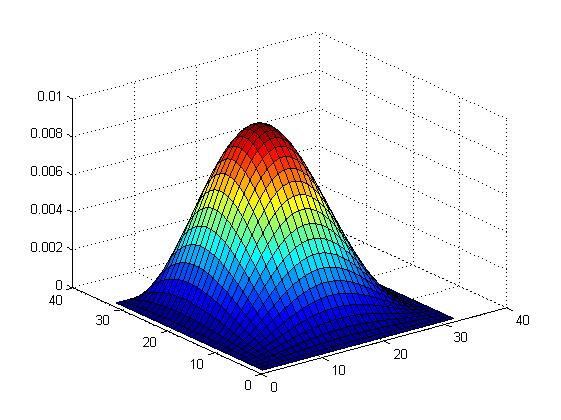}  }{\textit{Mollifier} 3}
&
\hspace{-1 cm}
\subf{\includegraphics[scale=.3]{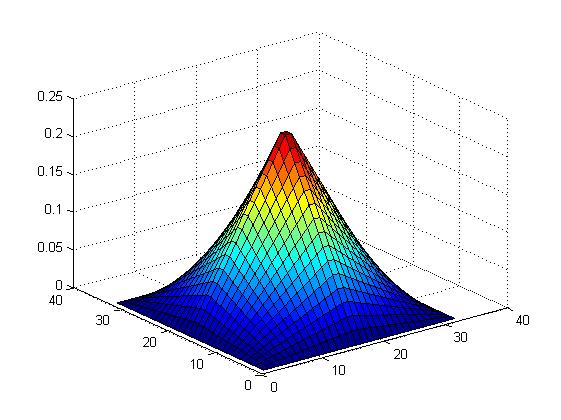} }{\textit{Mollifier} 4} 
\end{tabular}
\caption{\textit{Mollifiers} utilizados para amolecer o ruído Branco na equação KPZ}
\label{Molifiers}
\end{figure}

  No primeiro experimento vamos amolecer o ruído branco no espaço seguindo o procedimento da Secção (\ref{SecAmol})  usando o \textit{mollifier}
 \begin{equation}\label{mol1}
 \mathcal{G}_{\kappa}(t,x) = e^{-\left(\frac{1}{1-(\kappa\cdot t)^2} + \frac{1}{1-(\kappa\cdot x)^2}\right)}  
 \end{equation}
  onde o parâmetro $\kappa$ regula o tamanho do suporte. Agora variamos
   o parâmetro $\kappa$ tomando os valores $1/2,1/4,1/8,1/16$ e resolvemos, em cada caso, a 
   equação estocástica do calor usando o código \textit{Euler_Imp_Galerkin} implementado para o 
   método de Euler Galerkin semi-implícito. Logo, resolvemos a equação KPZ com ruído branco 
   amolecido usando o código \textit{EFMH_KPZ.m} implementado para o método de elementos finitos 
   mistos e híbridos adaptado para KPZ no Capitulo \ref{cap3}. Note, que novo problema cai dentro 
   da classe de problemas para a qual existem provas rigorosas de existência e unicidade 
   da solução (veja \cite{Den77}) e de convergência do método de elementos finitos mistos 
   híbrido com decomposição de domínio (veja \cite{KPP00}). O ruído branco é gerado pelo código
   \textit{Ruido.m} que pela sua vez utiliza a função de matlab \textit{randn} para gerar números 
   aleatórios com distribuição normal. A função\textit{ Ruido.m} possui um algoritmo para efetuar 
   o amolecimento com a flexibilidade de escolher o \textit{mollifier}.\\
\\     
   \indent No \textit{Experimento 1}(\ref{Exp1}) do Capitulo (\ref{cap3}) para a equação KPZ determinística, a escolha da relação $\Delta t/ \Delta x \approx 1/16$ foi bem sucedida na aproximação da solução. Tal escolha foi possível porque a solução era uma função muito regular e o tempo suficientemente curto. Porém, esta não é uma escolha adequada para o problema  estocástico onde a solução é muito mais áspera.  De fato, as observações dos experimentos onde o termo fonte é uma versão amolecida do ruído branco permitiram estimar uma relação $\Delta t = (\Delta x)^3$. Quer dizer  que para obter duas soluções até o tempo $T=1$ usando diferentes versões amolecidas de uma única realização do ruído branco, teríamos que armazenar matrizes da ordem  $N\times N^3$ pois a mesma realização não pode ser gerada duas vezes. O que vamos fazer é armazenar matrizes $N\times N$ e avançar até o tempo $T = 1/N^2$. O número de pontos na partição $N$ tem que variar junto com o tamanho do suporte do \textit{mollifier}  para conseguir capturar flutuações cada vez mais pequenas já que a solução se torna cada vez mais áspera.\\
   \\
\indent As duas soluções vão ser calculadas até os tempos 
   $\frac{1}{256^2},\frac{1}{512^2},\frac{1}{1024^2}$ para cada valor  do 
   parâmetro $\kappa$. A Figura \ref{HCKPZ} mostra como os perfis 
   das duas soluções vão se aproximando à medida que diminuímos o tamanho do suporte.\\
   \\
\indent Sabemos de (\ref{Cons}) que a diferença entre a transformada de Hopf-Cole e a solução amolecida 
  deve ser uma constate $C_{\kappa}$ vezes $t$ assim, podemos estimar cada uma destas constantes. 
  A Figura \ref{Const} mostra o crescimento exponencial das constantes. Estes resultados 
  apoiam a conjetura de que o processo de renormalização é necessário e que  
  as duas soluções coincidem no limite.   Repetimos este experimento para diferentes 
  \textit{mollifiers} cujos gráficos são mostrados na Figura \ref{Molifiers}. Na Figura 
  \ref{ConstMol} podemos ver o comportamento das constantes obtidas para cada um deles.
  
\newpage 
\section{Processo de renormalização da equação KPZ e experimentos numéricos}

   O processo de renormalização que vamos adotar neste trabalho é o mesmo que apresentado em \cite{Hai11}. 
   A idéia é definir um processo recursivo partindo de uma versão amolecida $\xi_{\kappa}$ do ruído 
   branco $\xi$ como em (\ref{MolR}). Vejamos qual foi a ideia desenvolvida por Hairer para
   chegar neste processo de renormalização.\\
   
   Chamemos de $\Pi_0$ a projeção ortogonal de $L^2([a,b])$ sobre o espaço das funções constates
   e definamos  $\Pi_0^{\perp} = 1 - \Pi_0$. Consideremos o conjunto de índices $\mathcal{T}_2$ 
   formado por todas as árvores binárias onde cada $\tau \in \mathcal{T}_2$ é escrito como 
   $\tau = [\tau_1,\tau_2]$ e é formado pela sua raiz e as árvores $\tau_1$ e $\tau_2$. 
   Agora vamos definir o processo
   
   \begin{equation*}
       \partial_t X^{\bullet}_\kappa = \partial^2_x X^{\bullet}_\kappa + \Pi_0^{\perp} \xi_\kappa,
   \end{equation*}
 \\
 e para cada árvore $\tau = [\tau_1,\tau_2]$ definimos $X^\tau$ como sendo a solução estacionária de
 
 \begin{equation*}
    \partial_t X^\tau_\kappa = \partial^2_x X^\tau_\kappa + \Pi_0^\perp \left( \partial_x X_\kappa^{\tau_1}, X_\kappa^{\tau_2}\right).
 \end{equation*}
 
 Definamos agora $Y_\kappa^\bullet = X^\bullet_\kappa + \sqrt{2} B(t)$ onde $B$ é um movimento Browniano. Então, pode ser demonstrado que existem constantes $C_\kappa^\tau$ tais que as soluções $Y_\kappa^\tau$ de 
 \begin{equation}
     \partial_t Y^\tau_\kappa = \partial^2_x Y^\tau_\kappa + \Pi_0^\perp \left( \partial_x Y_\kappa^{\tau_1}, Y_\kappa^{\tau_2}\right) - C_\kappa^\tau
 \end{equation}
 com condição inicial $Y_\kappa^\bullet(0) = X_\kappa^\bullet(0)$, tem limite quando $\kappa \longrightarrow 0$ e este limite é independente da escolha do \textit{mollifier} $\varphi$. A razão pela qual foi definido o processo $Y_\kappa^\tau$ é que, pelo menos no nível formal, se consideramos
 
 \begin{equation}
 h_\kappa (t) = \sum_\tau Y_\kappa^\tau,
 \end{equation}
\\   
   onde $\tau$ são todos os nós interiores (i.e. nós que não são folhas) então $h_\kappa^\tau$ resolve a equação
   
\begin{equation}\label{KPZR1}
   \partial_t h_\kappa= \partial^2_x h_\kappa + \left(\partial_x h_\kappa \right)^2 + \xi_\kappa - \sum_\tau C_\kappa^\tau.
\end{equation}

   Podemos notar a similaridade da equação (\ref{KPZR1}) com a equação KPZ. Ao longo do trabalho \cite{Hai11} é provado rigorosamente que o processo limite $ h = \lim_{\kappa \longrightarrow 0} h_\kappa $ existe e coincide com a transformada de Hopf-Cole da equação estocástica do calor. Além disso, são calculadas explicitamente as constantes $C_\kappa^\tau$ 
   
\begin{align*}
C_{\kappa}^{(1)} &= \frac{1}{\kappa} \int_{\mathbb{R}} \varphi^2(x) dx, \\
C_{\kappa}^{(2)} &= \frac{4\pi}{\sqrt{3}}\vert\log \kappa\vert - 8\int_{\mathbb{R_+}}\int_{\mathbb{R}} \frac{x\varphi'(y) \varphi(y) \varphi^2(y) \log\varphi}{x^2-xy+y^2} dx dy, \\
C_{\kappa}^{(3)} &= -\frac{C_2}{4}.
\end{align*}

Assim, a constante de renormalização fica determinada como $C_{\kappa} = C_{\kappa}^{(1)} + C_{\kappa}^{(2)} + C_{\kappa}^{(3)}$ e a versão da equação KPZ que devemos resolver é

\begin{equation}\label{KPZR}
\partial_t h_{\kappa} = \partial^2_x h_{\kappa} +  \left( \partial_x h_{\kappa} \right)^2   + \xi_{\kappa} - C_{\kappa}.
\end{equation}

\begin{figure}[h!]\label{KPZNN1}
\includegraphics[scale=.5]{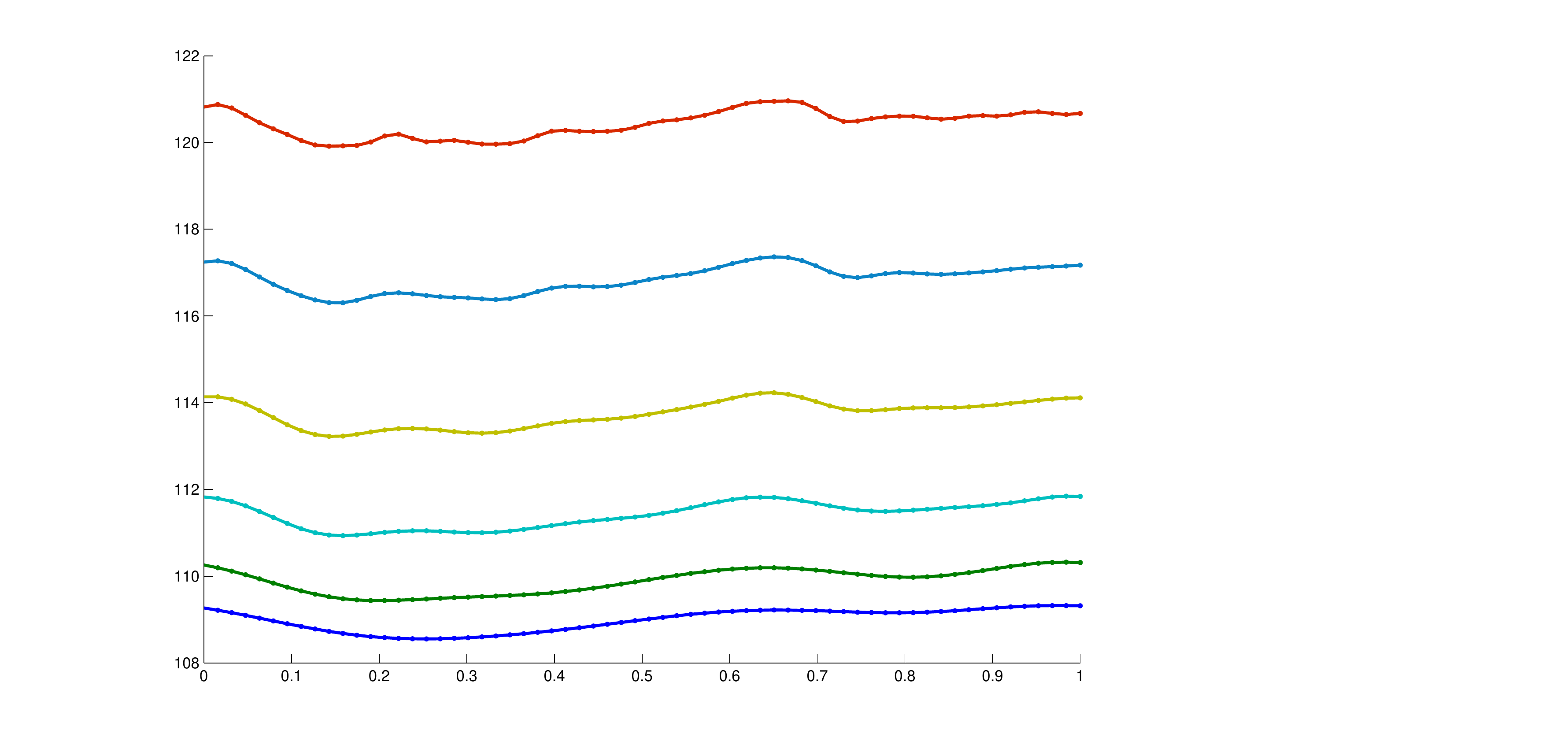} 
\caption{Soluções obtidas usando o código \textit{EFMH_KPZ} na equação (\ref{KPZR}) ignorando a constante de renormalização para valores do parâmetro $\kappa = 0.5000,  0.4000,  0.3333,  0.2857,  0.2500, 0.2222$. Pode-se observar que as soluções vão ficando mais distantes.}
\end{figure}

 Inicialmente, vejamos o que acontece no caso em que simplesmente ignoramos o processo de renormalização, ou seja, aplicar o método de elementos finitos mistos e híbridos na equação (\ref{KPZR}) sem a constante $C_{\kappa}$ para a mesma realização do ruído branco e diferentes graus de amolecimento. A Figura \ref{KPZNN1} mostra como, ao diminuir o parâmetro de amolecimento ($\kappa \longrightarrow 0$),  os perfis das soluções obtidas se aproximam cada vez uma das outras, mas, note que as alturas médias ficam cada vez mais distantes.
 
\begin{figure}[h!]
\hspace*{-1cm}   
\begin{tabular}{lll}
\hspace{-1.5 cm}

\subf{\includegraphics[scale=.3]{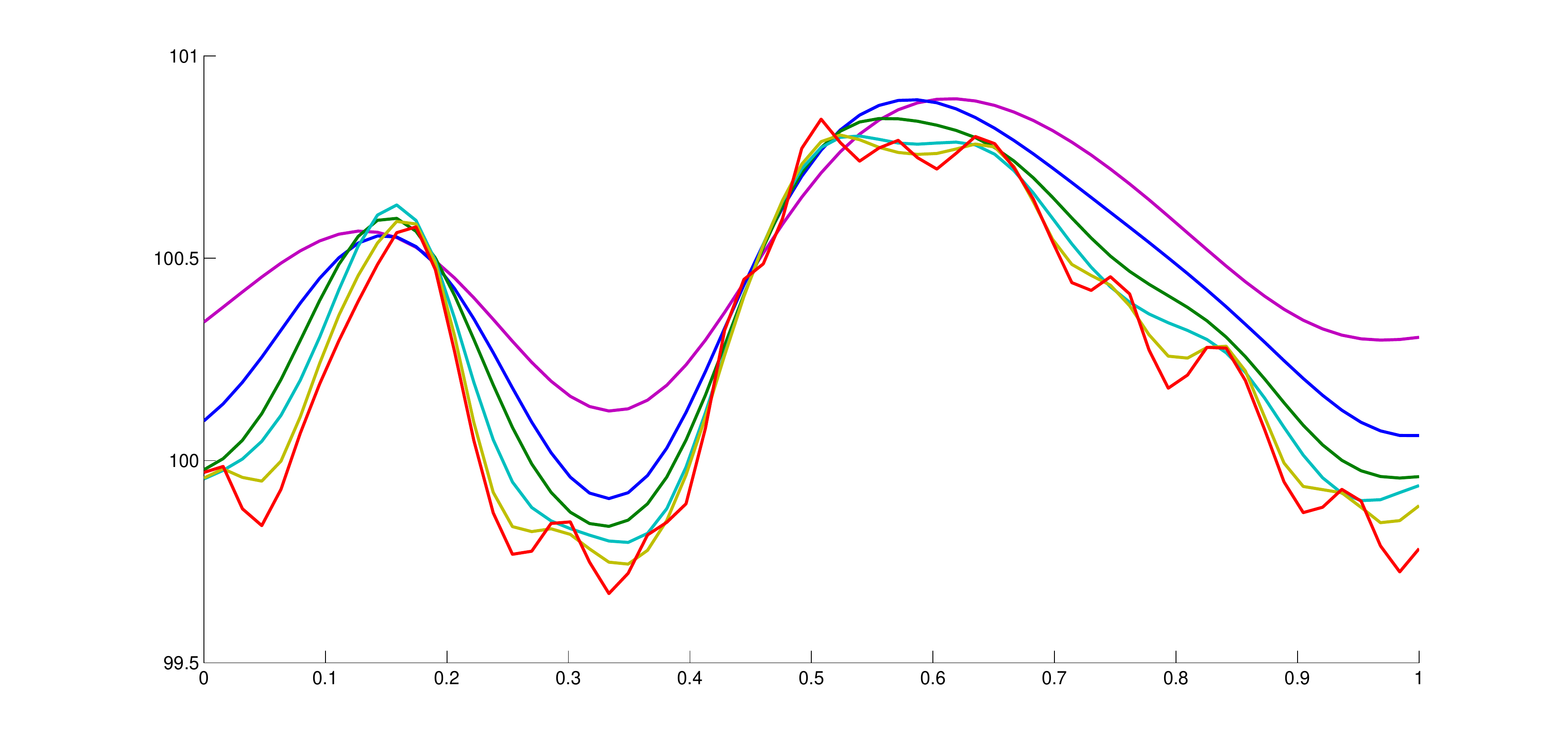} }
     {Ruído 1 Amolecido.}
&
\hspace{-1 cm}
\subf{\includegraphics[scale=.3]{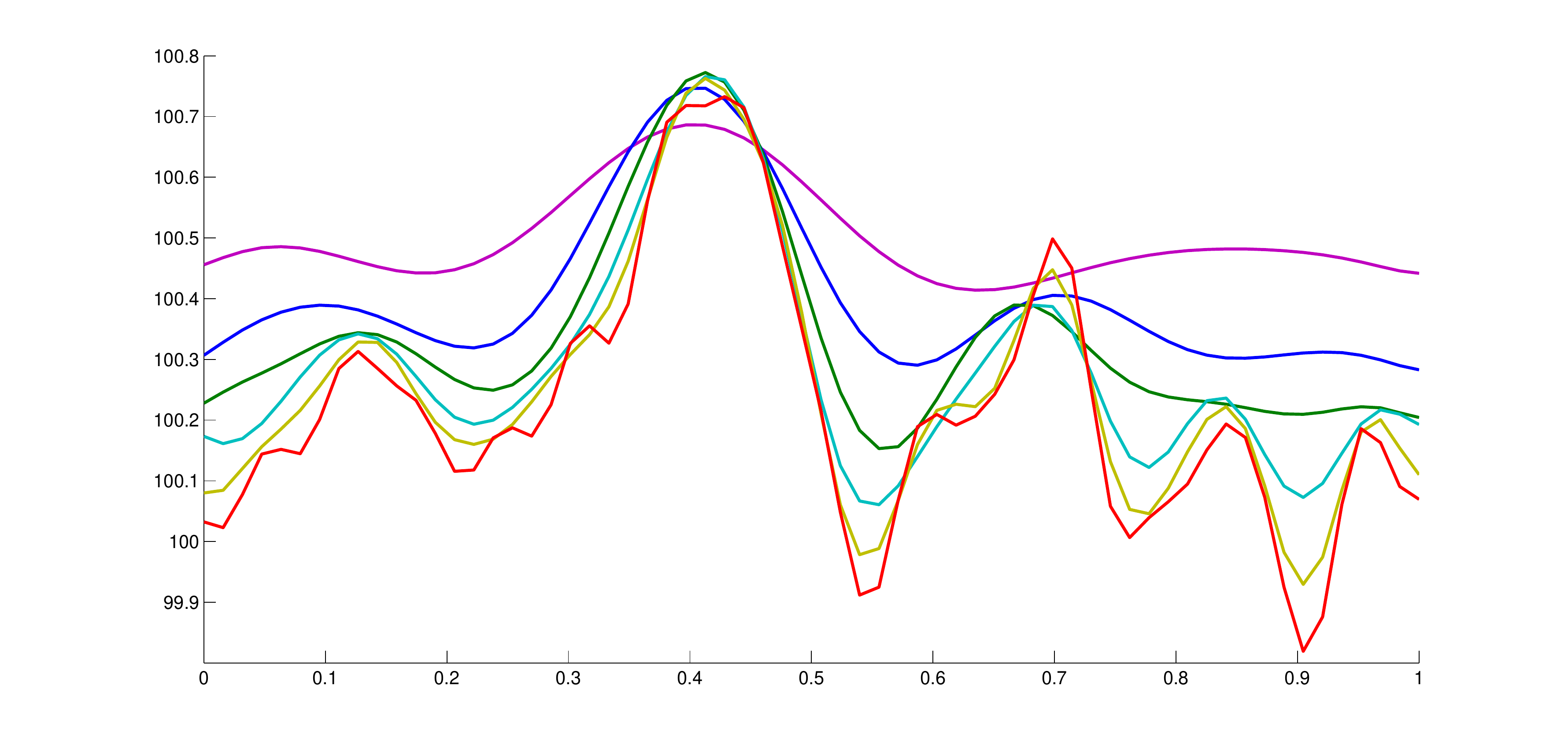} }
     {Ruído 2 Amolecido} 
\\
\hspace{-1.5 cm}
\subf{\includegraphics[scale=.3]{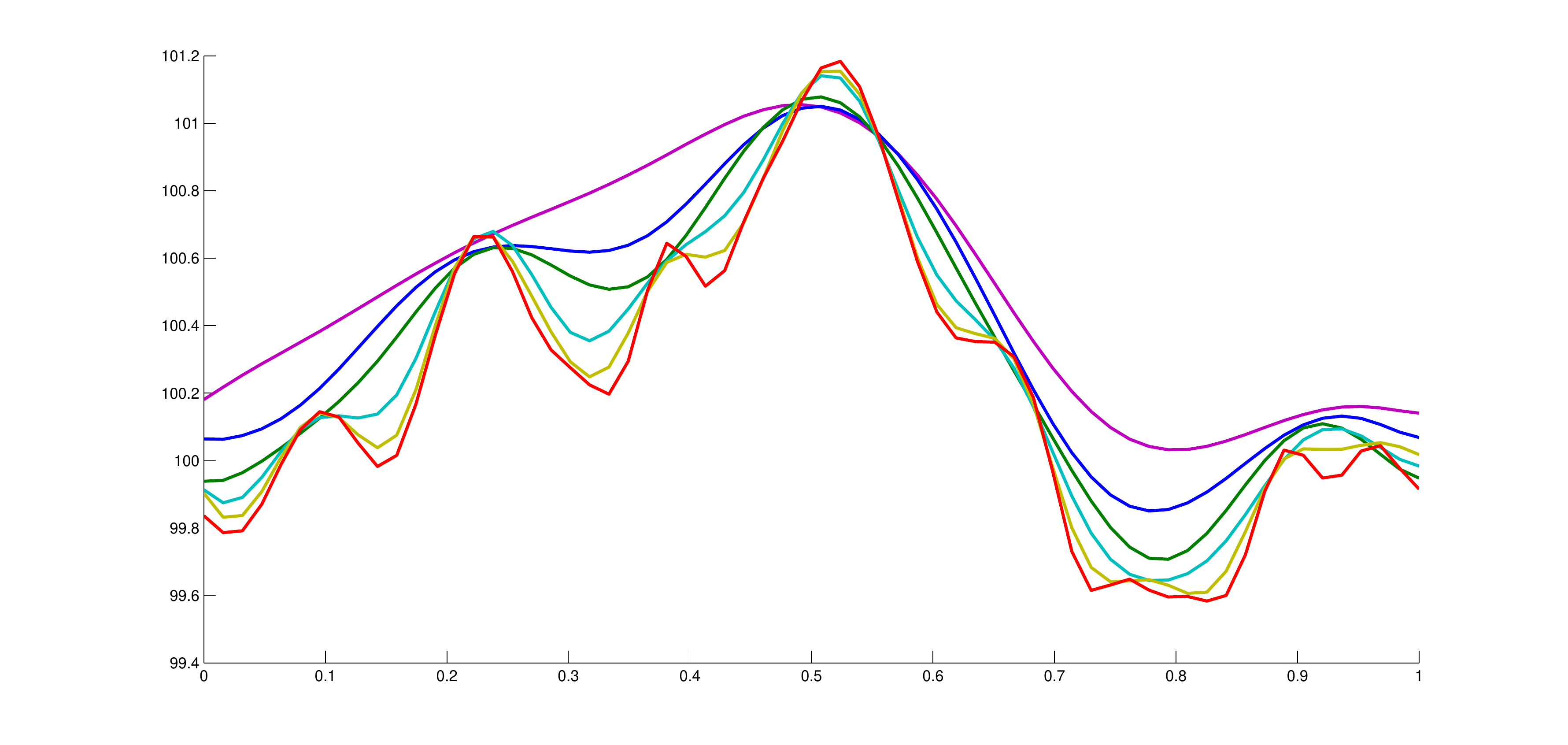} }
     {Ruído 3 Amolecido}
&
\hspace{-1 cm}
\subf{\includegraphics[scale=.3]{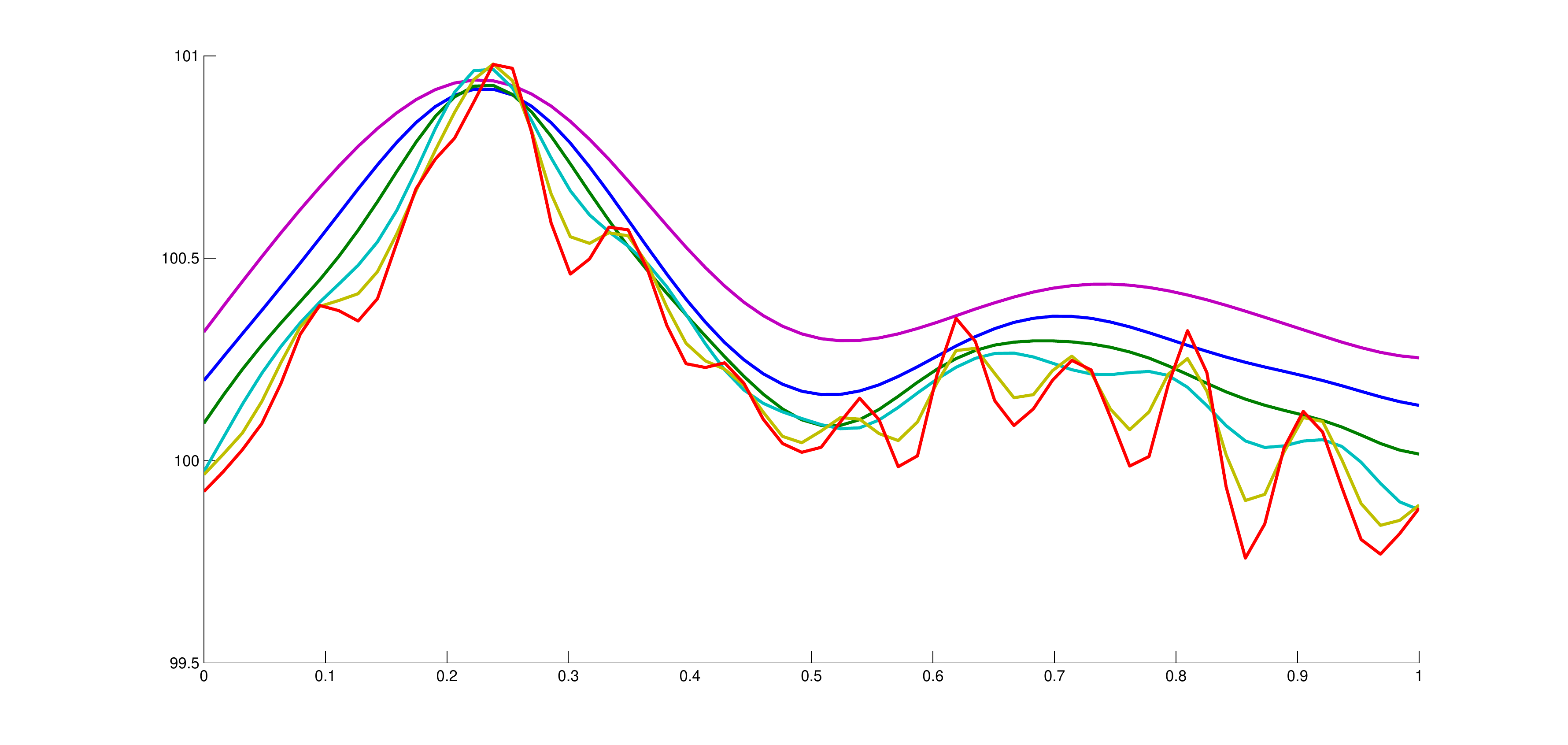} }
     {Ruído 4 Amolecido}
\\
\end{tabular}
\caption{Soluções obtidas usando o código \textit{EFMH_KPZ} na equação (\ref{KPZR}) incluindo a constante de renormalização para valores do parâmetro $\kappa = 0.5000,  0.4000,  0.3333,  0.2857,  0.2500, 0.2222$ sobre a mesma realização de 4 ruídos brancos. Pode-se observar que as soluções vão ficando próximas à medida que diminuímos o parâmetro de amolecimento $\kappa$}
\label{KPZSN}
\end{figure}

\newpage
   Agora vamos repetir o experimento descrito anteriormente mas desta vez incluindo a constante 
   de renormalização. A Figura \ref{KPZSN} mostra como tanto as alturas médias  quanto os perfis  
   vão ficando cada vez mais próximos à medida que diminuímos o parâmetro de amolecimento.

\newpage
   Vamos fazer um estudo de refinamento de malha para verificar que o erro, no sentido de (\ref{erro}), diminui à medida que fazemos $\kappa \longrightarrow 0$. Os experimentos serão feitos sobre $100$ realizações do ruído branco diminuindo o tamanho do parâmetro de amolecimento $\kappa$ junto com a norma da partição (uniforme) $\Delta x$. Esta variação é natural, uma vez que, quando aumenta a rugosidade da versão amolecida do ruído branco, precisamos diminuir a norma da partição para capturar as pequenas flutuações que aparecem. Desta maneira definimos $N_\kappa = \frac{N_0}{\kappa}$ onde $\Delta x = (b-a)/N_\kappa$. É claro que se $\kappa \longrightarrow 0$ então $N_\kappa \longrightarrow \infty$. Os erros são calculados pela formula

\begin{equation}\label{erro1}
\Vert h^{N_\kappa+1}_{\kappa,t} - h^{N_\kappa}_{\kappa,t}\Vert = \left(\mathbb{E}\left[\int_0^1 \vert h^{N_\kappa+1}_{\kappa,t}(x) - h^{N_\kappa}_{\kappa,t}(x)\vert^2 dx\right] \right)^{1/2}.
\end{equation}

\begin{figure}
\centering
\begin{minipage}[t]{.7\textwidth}
\centering
\vspace{0pt}
\hspace{-3cm}
\includegraphics[width=\textwidth, trim= 1cm 1cm 3cm 1cm,, clip=true]{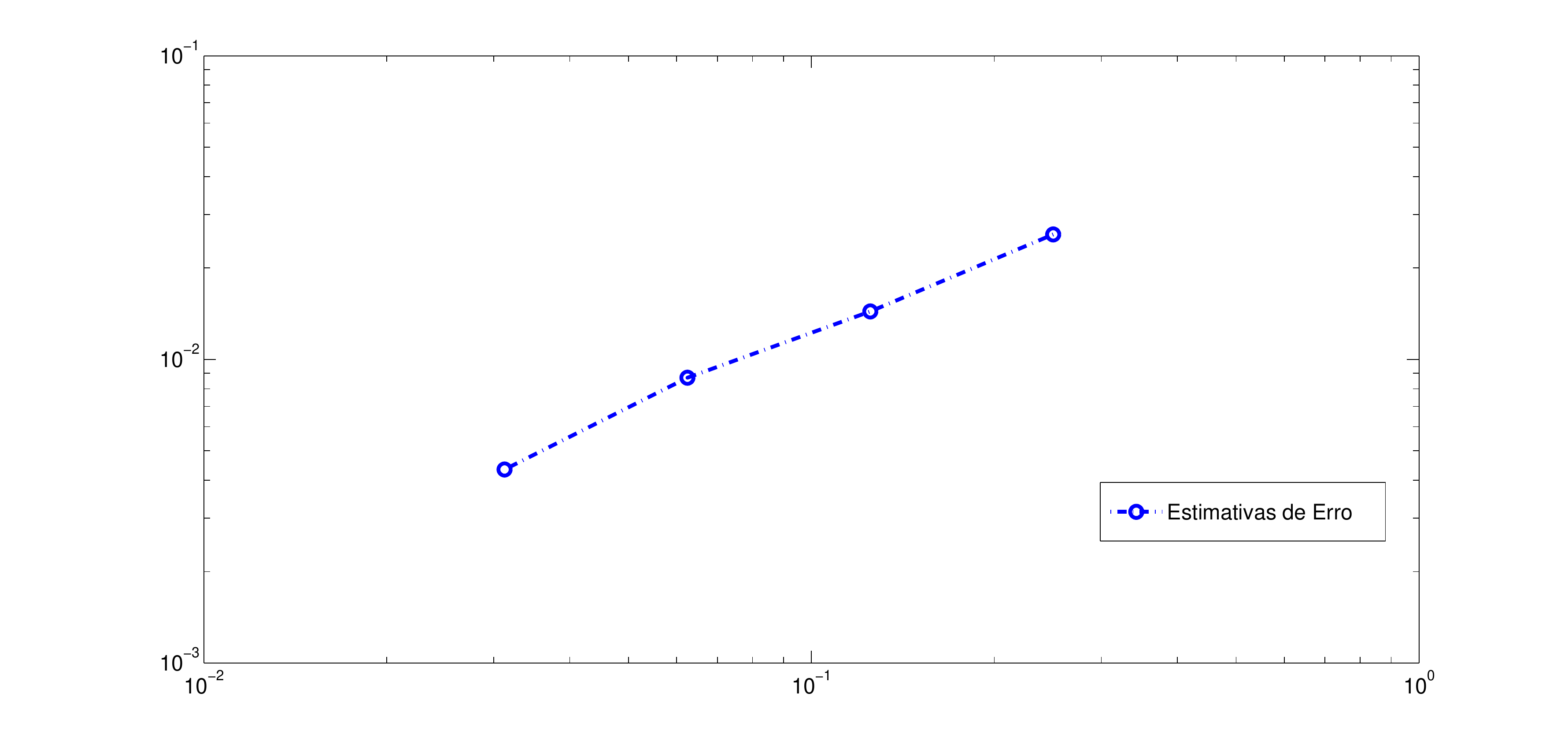}
\end{minipage}\hfill
\begin{minipage}[t]{.3\textwidth}
\centering
\vspace{1.5cm}
\hspace{0cm}
\renewcommand{\arraystretch}{1.2}
     	\begin{tabular}{c|c|c} 
       	    $N$  & $\kappa$ & $Erro$\\ 
			\hline  $8$   &  $0.2500$      &  $0.0258$  \\
			\hline  $16$  &  $0.1250$      &  $0.0144$  \\
			\hline  $32$  &  $0.0625$      &  $0.0087$  \\
		 	\hline  $64$  &  $0.0313$      &  $0.0043$  \\
	    \end{tabular}
\end{minipage}
\caption{Erros obtidos no sentido de (\ref{erro1}) sobre $100$ realizações do ruído branco variando o parâmetro de amolecimento $\kappa = 1, 1/2, 1/4, 1/8$ e o número de funções na expansão 
 $N_\kappa = 8/\kappa = 8, 16, 32, 64$.} \label{errotabla}
\end{figure}

Em nosso caso estamos trabalhando no intervalo $[a,b]=[0,1]$ até o tempo $T = 1$. Ao longo dos experimentos anteriores temos observado que a escolha da partição no tempo para garantir a convergência do algoritmo é da ordem de $\Delta t \approx (\Delta x)^3$. A Figura \ref{errotabla} mostra uma tabela e um gráfico log-log da variação do erro ao diminuir o par\^amento de amolecimento.

\chapter{Conclusões e perspectivas}\label{cap5}
\section{Conclusões}

  Nesta dissertação de mestrado trabalhamos com a aproximação numérica da solução de problemas de 
  valor inicial e de contorno associados à equação KPZ em dimensão 1. Durante o desenvolvimento 
  deste trabalho focamos no estudo dos métodos computacionais para equações diferenciais
  estocásticas semilineares e equações diferenciais não lineares. A ênfase se deu na parte 
  da modelagem computacional. Especificamente, foi apresentada uma coleção representativa de 
  experimentos que corroboraram recentes descobertas teóricas e permitiram testar nossa adaptação do 
  método de elementos
  finitos mistos e híbridos com decomposição de domínio à aproximação da solução da equação KPZ 
  acompanhada de um processo de renormalização introduzido em (\cite{Hai11}).\\

  Os métodos de \textit{Lord Rougemont}(\ref{LRM}), \textit{Euler Galerkin semi-implícito}(\ref{EGSI}) 
  e \textit{Milstein} (\ref{Milstein}) para SPDEs do tipo semilinear foram implementados 
  e os experimentos numéricos, que estão apresentados no Capitulo \ref{cap33}, mostraram boas evidências 
  de convergência 
  quando aplicados na equação estocástica do calor com ruído branco multiplicativo. Por exemplo, um
  estudo de erro mostrou que o erro de aproximação decai, em  cada um dos modelos supracitados. 
  Especificamente, o algoritmo 
  de \textit{Milstein} requer um esforço computacional menor do que o de \textit{Lord Rougemont} e 
  o de \textit{Euler Galerkin semi-implícito} (ver Figura \ref{ordem}). Os três algoritmos 
  produziram  soluções 
  muito próximas  para a mesma realização do ruído branco. Além disso, a transformada de Hopf-Cole 
  da solução da 
  equação estocástica do calor, obtida usando estes códigos, exibe as propriedades de crescimento 
  da rugosidade preditas teoricamente para o processo de deposição balística (ver Figura \ref{Rugo}).\\

 Por sua vez, a reformulação proposta via método de elementos finitos mistos e híbridos adaptado 
 para KPZ, foi bem sucedida na aproximação da solução do modelo determinístico. Isto 
 foi evidenciado quando aplicamos o método nos exemplos propostos. No caso determinístico, foi 
 possível comparar a solução numérica com a solução analítica, se consideramos condições de contorno 
 consistentes de Dirichlet para o problema. Variando as condições de fronteira de explícitas para periódicas, 
 percebemos que a solução numérica obtida preserva o comportamento das soluções do caso 
 determinístico, e que, além disso, as variações correspondem ao esperado no modelo físico.
 Portanto, temos evidências de convergência do método. No caso estocástico, ao aplicar o 
 código implementado para o algoritmo de elementos finitos mistos e híbridos com decomposição 
 de domínio na equação KPZ com ruído branco amolecido, observamos como o perfil da solução 
 obtida se aproxima do perfil da transformada de Hopf-Cole da solução da equação estocástica 
 do calor, ao diminuirmos o valor do parâmetro de amolecimento (\ref{HCKPZ}). Mais uma vez, 
 os resultados computacionais obtidos concordam com 
 os resultados preditos na teoria. Em outras palavras, as duas soluções são iguais no limite. \\
 
 Usando o processo de renormalização proposto em \cite{Hai11}, observou-se uma 
 compensação da divergência produzida pelo termo $(\partial_x h(t,x))^2$. Isto fica explícito 
 ao observarmos que as alturas médias e os perfis das soluções permaneceram cada vez mais próximos quando 
 $\kappa \longrightarrow 0$ (não importando a escolha do \textit{mollifier}). Observamos 
 que o tempo de cômputo empregado pelo algoritmo \textit{EHMH_DD.m} é maior do que o do algoritmo 
 de \textit{Milstein}. \\

\section{Perspectivas}

  Como continuação deste trabalho, pretendemos utilizar as ideias e conceitos apontados aqui a uma classe 
  mais ampla de SPDEs. Como exemplos de casos de estudo temos os problemas associados à dinâmica 
  de fluidos em meios porosos heterogêneos em dimensão maior do que 1, envolvendo sistemas de 
  equações diferenciais estocásticas, em linha com os trabalhos (\cite{FuP03,BFPS08,Abr15}). 
  Uma análise numérica mais profunda dos métodos de resolução associados aos novos problemas a serem
  estudados é também um dos objetivos da continuação deste trabalho. Com efeito, cumpre destacar que 
  uma teoria geral para a boa colocação de modelos estocásticos da KPZ ainda é um tema em aberto. 
  A implicação é óbvia, no que diz respeito à identificação dos espaços funcionais de aproximação
  adequados para modelos KPZ.\\

Algumas abordagens possíveis de pesquisa são

\begin{itemize}
   
  \item Estudar o desempenho do método de Galerkin espectral quando aplicado nas dimensões 
  temporal e espacial, diferentemente de seu uso combinado com diferenças finitas, como abordado 
  no Capitulo \ref{cap33} deste trabalho. Essa nova abordagem pode ser útil no desenvolvimento de 
  algoritmos com melhor desempenho computacional, uma vez que podemos aproveitar as vantagens 
  de cálculo que as bases exponenciais oferecem na discretização da dimensão temporal, da mesma forma que 
  foi feito na discretização da dimensão espacial. Além disso, esta abordagem pode simplificar
  a análise numérica já que a homogeneidade da discretização restringe o estudo à teoria espectral. 

  \item Considerar outros tipos de bases para a expansão das aproximações. Por exemplo, funções
  \textit{wavelets}, cuja propriedades de ortogonalidade e compacidade do suporte 
  permitem seu uso em métodos espectrais. Além disso, também permitem localizar os cálculos.
  Sendo assim, métodos espectrais que utilizem este tipo de base, podem ser melhor sucedidos 
  na tarefa de aproximar a solução de  problemas de natureza \textit{stiff}, os quais aparecem 
  usualmente na din\^amica de fluidos em meios porosos. 

  \item Explorar outras alternativas para obtenção de aproximações do integrando da integral estocática (\ref{semiL}). Tal abordagem é pertinente pois, neste trabalho, focamos nossa atenção em séries de Taylor para obter tais aproximações, as quais dependem da existência de derivadas de Frechet, que nem sempre são possíveis de obter. Outras aproximações polinomiais como interpolação ou projeção poderiam ser estudadas.

\item Empregar as ideias que foram utilizadas no Capítulo (\ref{cap33}) para propor alternativas de aprimoramento do desempenho computacional dos métodos espectrais, especificamente do método de elementos finitos mistos e híbridos aplicado na obtenção de aproximações de soluções de equações diferenciais estocásticas.
  
\end{itemize}


\phantomsection
\addcontentsline{toc}{chapter}{Bibliography}
\bibliographystyle{siam}
\bibliography{tese1}

\backmatter
\phantomsection
\addcontentsline{toc}{chapter}{Index}

\end{document}